\documentclass{amsart}

\usepackage{amsmath}
\usepackage{amsthm}
\usepackage{amsopn}
\usepackage{amssymb}
\usepackage{supertabular}
\usepackage[all]{xy}
\usepackage{graphicx}
\usepackage{lscape}

\newcommand{\abs}[1]{\lvert #1 \rvert}

\newcommand{\br}[1]{\overline{#1}}

\newcommand{\td}[1]{\widetilde{#1}}

\newcommand{\bra}[1]{\langle #1 \rangle}
\newcommand{\ZZ}{\mathbb{Z}}

\newcommand{\CC}{\mathbb{C}}
\newcommand{\FF}{\mathbb{F}}
\newcommand{\capeq}{\overset{\scriptscriptstyle{\cap}}{=}}
\newcommand{\dotin}{\overset{\scriptstyle{\centerdot}}{\in}}
\newcommand{\dotequiv}{\overset{\scriptstyle{\centerdot}}{\equiv}}
\newcommand{\customsizeA}[1]{\small{#1}}

\theoremstyle{definition}
 \newtheorem{thm}{Theorem}[section]
 \newtheorem{cor}[thm]{Corollary}
 \newtheorem{lem}[thm]{Lemma}
 \newtheorem{prop}[thm]{Proposition}
 \newtheorem{rmk}[thm]{Remark}
\newtheorem{defn}[thm]{Definition}

\newtheorem*{conventions}{Conventions}
\newtheorem{procedure}[thm]{Procedure}
\numberwithin{equation}{section}
\DeclareMathOperator{\ext}{Ext}
\DeclareMathOperator{\holim}{holim}

\title{Root invariants in the Adams spectral sequence}
\author[Mark Behrens]{Mark Behrens$\sp 1$}
\address{
Department of Mathematics \\
Massachusetts Institute of Technology \\
Cambridge, MA 02139, U.S.A.}

\begin{document}

\begin{abstract}

Let $E$ be a ring spectrum for which the $E$-Adams spectral
sequence converges.  We define a variant of Mahowald's root invariant
called the `filtered root invariant' which
takes values in the $E_1$ term of the $E$-Adams spectral sequence.  
The main theorems of this paper concern
when these filtered root invariants detect the actual root invariant, and
explain a relationship between filtered root invariants and
differentials and compositions in the $E$-Adams spectral sequence.
These theorems are compared to some known computations of root invariants at
the prime $2$. 
We use the filtered root invariants to compute some low dimensional root
invariants of $v_1$-periodic elements at the prime $3$.
We also compute the root invariants of
some infinite $v_1$-periodic families of elements at the prime $3$.

\end{abstract}

\maketitle

\tableofcontents

\footnotetext[1]{The author is partially supported by the NSF.}


\section{Introduction}

In his work on metastable homotopy groups \cite{Mahowald}, Mahowald
introduced an invariant that associates to every element $\alpha$ 
in the stable
stems a new element $R(\alpha)$ called the \emph{root invariant} of
$\alpha$.  
The construction has indeterminacy and so $R(\alpha)$ is in general only a
coset.
The main result of \cite{MahowaldRavenel} indicates a deep
relationship between elements of the stable homotopy groups of spheres
which are
root invariants and their behavior in the EHP spectral sequence.
Mahowald and Ravenel conjecture in \cite{MahowaldRavenel2} that, loosely
speaking, the root
invariant of a $v_n$-periodic element is $v_{n+1}$-periodic.  Thus the root
invariant is related simultaneously to unstable and chromatic phenomena.

The conjectural relationship between root invariants and the chromatic
filtration 
is based partly on computational evidence. 
For instance, at the prime $2$, we have \cite{MahowaldRavenel},
\cite{Johnson}
$$ R(2^i) = \begin{cases}
\alpha_{4t} & i = 4t \\
\alpha_{4t+1} & i = 4t+1 \\
\alpha_{4t+1}\alpha_1 & i = 4t+2 \\
\alpha_{4t+1}\alpha_1^2 & i = 4t+3 
\end{cases}
$$
while at odd primes we have \cite{Sadofsky}, \cite{MahowaldRavenel}
$$ R(p^i) = \alpha_i $$
demonstrating that the root invariant sends $v_0$-periodic families to
$v_1$-periodic families.  For $p \ge 5$ it is also known \cite{Sadofsky},
\cite{MahowaldRavenel} that
$$ \beta_i \in R(\alpha_i) $$
and 
$$ \beta_{p/2} \in R(\alpha_{p/2}). $$
These computations led the authors of \cite{MahowaldRavenel} to regard the
root invariants as \emph{defining} the $n^{\mathrm{th}}$ 
Greek letter elements as
the root invariants of the $(n-1)^\mathrm{st}$ 
Greek letter elements when the relevant
Smith-Toda complexes do not exist.  

Other evidence of the root invariant raising chromatic filtration is seen
in the cohomology of the Steenrod algebra.
Mahowald and Shick
define a chromatic filtration on $\ext_A(\FF_2, \FF_2)$
in \cite{MahowaldShick}, and Shick proves that
an algebraic version of the root invariant increases chromatic filtration
in this context in \cite{Shick}.  

The
Tate spectrum computations of
\cite{DavisMahowald}, \cite{tBP2}, \cite{MahowaldShick2},
\cite{GreenleesSadofsky}, \cite{AndoMoravaSadofsky}, and
\cite{GreenleesMay} indicate that the $\ZZ/p$-Tate spectrum of a
$v_n$-periodic cohomology theory is $v_n$-torsion.  The root invariant is
defined using the Tate spectrum of the sphere spectrum, so the results of
the papers listed above provide even more evidence that the root invariant
of a $v_n$-periodic element should be $v_n$-torsion.

The purpose of this paper is to introduce a new variant of the root
invariant called the \emph{filtered root invariant}.
We apply the filtered root invariant
towards the computation of some new root invariants at the prime~$3$ 
in low dimensions.  We also compute the root invariants of some
infinite $v_1$-families as infinite $v_2$-families at the prime $3$.  We
will also describe how the theory of filtered root invariants given in
this paper works out to give an alternative perspective on known
computations of root invariants at the prime~$2$.

Let X be a finite $p$-local spectrum, and let $E$ be a ring spectrum for
which the $E$-Adams spectral sequence ($E$-ASS) converges.
Given $\alpha \in \pi_t(X)$, we define its 
filtered root invariants $R^{[k]}_E(\alpha)$ to be a sequence of cosets of
$d_r$-cycles in the
$E^{k,*}_1$ term of the
$E$-ASS converging to $\pi_*(X)$, where $r$ depends on $k$. These invariants 
govern the
passage between the $E$-root invariant $R_E(\alpha)$, or the algebraic root
invariant $R_{alg}(\alpha)$, and the elements in the $E$-ASS that
detect $R(\alpha)$.  The passage of each of these filtered root invariants
to the next is governed by differentials and compositions in the
$E$-ASS.  Our method is to use algebraic and $E$-cohomology computations to
determine the first of these filtered root invariants and then iteratively
deduce the higher ones from differentials and compositions.  Under
favorable circumstances, the last of these filtered root invariants will
detect the root invariant.

In Section~\ref{sec:filt}, we discuss the $E$-Adams resolution, and
various filtered forms of the Tate spectrum which we shall be using to
define the the filtered root invariants.

In Section~\ref{sec:defs}, we  
define the filtered root invariants.  We also recall the definitions of 
the root
invariant, the $E$-root invariant, and the algebraic root invariant.  We
will sometimes refer to the root invariant as the ``homotopy root
invariant'' to distinguish it from all of the other variants being used.

Our main results are most conveniently stated using the language of Toda
brackets introduced in the appendix of \cite{Shipley}.  
In Section~\ref{sec:toda}, we define a variant
of the Toda bracket, the $K$-Toda bracket, which is taken to be an attaching 
map in a fixed finite CW-complex $K$. 
Some properties of $K$-Toda brackets are introduced.  We also define a
version on the $E_r$ term of the $E$-ASS.

In Section~\ref{sec:results}, we state the main results which relate the
filtered root invariants to root invariants, Adams differentials,
compositions, $E$-root invariants, and algebraic root invariants.  
If $R^{[k]}_E(\alpha)$ contains
a permanent cycle $x$, then $x$ detects an element of $R(\alpha)$ modulo
Adams filtration $k+1$.  If $R^{[k]}_E(\alpha)$ does not contain a
permanent
cycle, there is a formula relating its $E$-Adams differential to the next
filtered root invariant.  Thus if one knows a filtered root invariant, one
may sometimes deduce the next one from its $E$-Adams differential.  There
is a similar result concerning compositions.  The zeroth filtered root
invariant $R^{[0]}_E(\alpha)$ is the $E$-root invariant. The first filtered
root invariant $R^{[1]}_E(\alpha)$ is the $E \wedge \br{E}$-root invariant
($\br{E}$ is the fiber of the unit of $E$).
If $E=H\FF_p$,
then the first non-trivial filtered root invariant is the algebraic root
invariant $R_{alg}(\alpha)$.

We present the proofs of the main theorems of
Section~\ref{sec:results} in Section~\ref{sec:proofs}.  The
proofs are technical, and it is for this reason that they have been
relegated to their own section.

In Section~\ref{sec:bo} we give
examples of these theorems with $E=bo$, and show how the $bo$ resolution
computes the root invariants of the elements $2^k$ at the prime $2$.  This
was the motivating example for this project.

We
will need to make some extensive computations in an algebraic 
Atiyah-Hirzebruch spectral sequence (AAHSS).  This is the spectral sequence
obtained by applying $\ext_{BP_*BP}(BP_*, BP_*(-))$ to the cellular
filtration of projective space.
We introduce this spectral sequence in
Section~\ref{sec:aahss}.  We compute the $d_r$ differentials for small $r$
using formal group methods.  

Calculating homotopy root invariants from our filtered root invariants is a
delicate business.  In an effort to make our low dimensional calculations
easier to follow and less ad hoc, we spell out our methodology 
in the form of Procedure~\ref{procA} 
which we follow
throughout our low dimensional calculations.  The description of this
procedure is the subject of Section~\ref{sec:proc}.

In Section~\ref{sec:anss}, we compute some $BP$-filtered root
invariants of the elements $p^i$ as well as of the elements $\alpha_{i/j}$.
We find that at every odd prime $p$
\begin{align*}
\pm \alpha_i & \in R_{BP}^{[1]}(p^i) \\
\pm \beta_{i/j} & \in R_{BP}^{[2]}(\alpha_{i/j})
\end{align*}
These filtered root invariants hold at the prime $2$ modulo an
indeterminacy which is identified, but not computed, in this paper.  
The only exceptions are 
the cases $i=j=1$ and $i=j=2$ at the prime $2$ 
(these cases correspond to the existence of
the Hopf invariant $1$ elements $\nu$ and $\sigma$).
These filtered root invariants are computed by means of manipulation of
formulas in $BP_*BP$ arising from $p$-typical formal groups.

In Section~\ref{sec:beta1}
we compute the root invariants
$$ R(\beta_1) = \beta_1^p $$
at primes $p > 2$.
These root invariants were announced without proof in \cite{MahowaldRavenel}.
This computation is accomplished by applying our theorems to the Toda
differential in the Adams-Novikov spectral sequence (ANSS).


In Section~\ref{sec:anss3}, we apply the methods described to compute some root
invariants 
of the Greek letter elements $\alpha_{i/j}$
that lie within the $100$-stem at the prime $3$.  At the prime
$3$, $\beta_i$ is known to be a permanent cycle for $i \equiv 0,1,2,5,6 
\pmod 9$ 
\cite{BehrensPemmaraju} and is conjectured to exist for $i \equiv 3 \pmod
9$.  The element $\beta_3$ is a permanent cycle.  One might expect from the
previous work for
$p \ge 5$ that $\beta_i \in R(\alpha_i)$ when $\beta_i$ exists.
Surprisingly, there is at least one instance where $\beta_i$ exists, yet is not
contained in $R(\alpha_i)$.  Our low dimensional computations of
$R(\alpha_{i/j})$ at $p=3$ are
summarized in Table~\ref{tab:3roots}.  
\begin{table}[ht]\label{tab:3roots}
\caption{Low dimensional root invariants of $\alpha_{i/j}$ at $p=3$}
\begin{center}
\begin{tabular}{c|c}
\hline
Element & Root Invariant \\
\hline
$\alpha_1$ & $\beta_1$ \\
$\alpha_2$ & $\pm \beta_1^2 \alpha_1$ \\
$\alpha_{3/2}$ & $-\beta_{3/2}$ \\
$\alpha_3$ & $\beta_3$ \\
$\alpha_4$ & $\pm \beta_1^5$ \\
$\alpha_5$ & $\beta_5$ \\
$\alpha_{6/2}$ & $\beta_{6/2}$ \\
$\alpha_{6}$ & $-\beta_6$ \\
\hline
\end{tabular}
\end{center}
\end{table}

All of these root invariants are $v_2$-periodic in the sense that they are
detected in $\pi_*(L_2 S^0)$ \cite{ShimomuraWang}.
A similar phenomenon happens at the prime $2$ with the root invariants of
$2^i$: they are not all given by the Greek letter elements $\alpha_{i}$, 
but the elements $R(2^i)$ are nevertheless $v_1$-periodic
\cite{MahowaldRavenel}.
The lesson we learn is that if one believes that the homotopy Greek letter
elements should be determined by iterated root invariants (as suggested in 
\cite{MahowaldRavenel}) then they will not always agree with the 
algebraic
Greek letter elements, even when the latter are permanent cycles.  

These
results will be partly generalized to 
compute $R(\alpha_i)$ for $i \equiv 0,1,5
\pmod 9$ at the prime $3$ in Section~\ref{sec:modcomp}.  
The remainder of this paper is
devoted to providing the machinery necessary for this computation.

In the ANSS the $\beta$ family lies in low Adams-Novikov filtration, but in
the ASS this family is in high filtration.  For the purposes of infinite
chromatic families, it is often useful to take both spectral sequences into
account simultaneously.  In Section~\ref{sec:mss} we explain how our
framework can be applied to the Mahowald spectral sequence to compute
algebraic root invariants.  As mentioned earlier, algebraic root invariants
are the first non-trivial $H\FF_p$ filtered root invariants.

Infinite families of Greek letter elements are constructed as homotopy
classes through the use of Smith-Toda complexes.  In their computation of
$R(\alpha_i)$ for $p \ge 5$, Mahowald and Ravenel introduce modified root
invariants \cite{MahowaldRavenel} which take values in the homotopy groups
of certain Smith-Toda complexes.
In Section~\ref{sec:mod} we adapt our results to modified root invariants.

Our modified root invariant methods are applied in
Section~\ref{sec:modcomp} to make some new computations of the root invariants
of some infinite $v_1$-periodic families at the prime $3$.  Specifically,
we are able to show that
$$ (-1)^{i+1} \beta_i \in R(\alpha_i) $$
for $i \equiv 0,1,5 \pmod 9$.

This paper represents the author's dissertation work.
The author would like to extend his heartfelt gratitude to his adviser, 
J.~Peter May, for his guidance and encouragement, and to Mark Mahowald, for
many enlightening conversations regarding the contents of this paper.


\begin{conventions}
Throughout this paper we will be working in the stable homotopy category
localized at some prime $p$.  We will always denote the quantity $q =
2(p-1)$, as usual.  All ordinary homology will be taken with $\FF_p$
coefficients.
If $p=2$,
let $P_N^M$ denote the stunted projective space with bottom cell in
dimension $N$ and top cell in dimension $M$.  Here $M$ and $N$ may be
infinite or negative.  
See \cite{MahowaldRavenel} for details.
If $p$ is odd, then projective space is replaced by
$B\Sigma_p$.  The complex $B\Sigma_p$ has a stable cell in every positive 
dimension congruent to $0$ or $-1$ mod $q$, and we will use the notation
$P_N^M$ to indicate the stunted complex with cells in dimensions between
$N$ and $M$.  When $M = \infty$, or when $N=-\infty$ 
the superscript or subscript may be omitted.

Given a spectrum $E$, the Tate spectrum
$$ \Sigma (E \wedge P)_{-\infty} = 
\Sigma \mathrm{holim}(E \wedge P_{-n}) $$
will be denoted $tE$.  To relate this notation to that in
\cite{GreenleesMay}, we have
$$ tE = t(\iota_* E)^{\ZZ/p}. $$
There is a unit $S^0 \rightarrow tE$.  For $E = S^0$, this is the inclusion
of the $0$-cell.  If $X$ is a finite complex, then the Segal conjecture for
the group $\ZZ/p$ \cite{SegalI}, \cite{SegalII}, \cite{RavenelSegal}
(also known as Lin's Theorem \cite{Lin}, \cite{LinExt} at $p=2$ and
Gunawardena's Theorem \cite{Gunawardena} for $p > 2$) implies that the map
\begin{equation}\label{eq:segalconj}
X = X \wedge S^0 \rightarrow X \wedge tS^0 = tX
\end{equation}
is $p$-completion (the last equality requires $X$ to be finite).

Suppose $A$ and $B$ are two subsets of a set $C$.  We shall write $A \capeq
B$ to indicate that $A \cap B$ is nonempty.  This is useful notation when
dealing with operations with indeterminacy.  If we are working over a ring
$R$, we shall use the notation $\doteq$ to indicate that two quantities are
equal modulo multiplication by a unit in $R^\times$.  We shall similarly
use the notation $\dotin$ for containment up to multiplication by a unit.

We will denote the regular ideal $(p, v_1, v_2, \ldots, v_{n-1}) \subseteq
BP_*$ by $I_n$.  

Finally, we will be using the following abbreviations for spectral
sequences.
\begin{description}

\item[ASS] The classical Adams spectral sequence.

\item[ANSS] The Adams-Novikov spectral sequence derived from $BP$.

\item[$E$-ASS] The generalized Adams spectral sequence derived from a ring
spectrum $E$.

\item[AHSS] The Atiyah-Hirzebruch spectral sequence.  We will
be using the form that computes stable homotopy groups from homology.

\item[AAHSS] The algebraic Atiyah-Hirzebruch spectral sequence, which
uses the cellular filtration to compute $\ext(X)$.

\item[MSS] The Mahowald spectral sequence, which computes $\ext(X)$ by
applying $\ext(-)$ to an Adams resolution of $X$.

\end{description}

\end{conventions}


\section{Filtered Tate spectra}\label{sec:filt}

Given a ring spectrum $E$, we will establish some notation for dealing with
the $E$-Adams resolution.  We will then mix the Adams filtration with the
skeletal filtration in the Tate spectrum $tS^0$.  These filtered Tate
spectra will carry the filtered root invariants defined in
Section~\ref{sec:defs}.  Our treatment of the Adams resolution follows
closely that of Bruner in \cite[IV.3]{Hinfty}.

For $E$ a ring spectrum, let $\br{E}$ be the fiber of the unit, so there is
a cofiber sequence
$$ \br{E} \rightarrow S^0 \xrightarrow{\eta} E. $$ 
For $X$ a spectrum, let
$W_k(X)$ denote the $k$-fold smash power $\br{E}^{(k)} \wedge X$.  
We shall also
use the notation $W_{k}^l(X)$ to denote the cofiber
$$ W_{k+l+1}(X) \rightarrow W_{k}(X) \rightarrow W_k^l(X). $$
We may drop the $X$ from the notation when $X=S^0$.
Note that with our definitions $W_k^{k-1}(X) \simeq \ast$.
The $E$-Adams resolution of $X$ now takes the form
$$ \xymatrix{
X \ar@{=}[r] & W_0(X) \ar[d] & W_1(X) \ar[d] \ar[l] & W_2(X) \ar[d] 
  \ar[l] & W_3(X) \ar[d] \ar[l] & \cdots \ar[l] \\
& W_0^0(X) & W_1^1(X) & W_2^2(X) & 
W_3^3(X)
} $$
The notation $W_k^l(X)$ is used because the $E$-ASS for $W_k^l(X)$ is
obtained from the $E$-ASS for $X$ by setting $E_1^{s,t} = 0$ for $s < k$
and $s > l$, and adjusting the differentials accordingly.  
If the resolution converges to the $p$-completion ($p$-localization) of
$X$, then $W_\infty(X) \simeq \ast$ in the $p$-complete ($p$-local) stable
homotopy category, and $W_s^\infty(X) \simeq W_s(X)$.

We shall denote
$E_r(X)$ for the $E_r$ term of the $E$-ASS for $X$.  An element of
$E_1^{s,t}(X) = \pi_{t-s}(W_s^s(X))$ is a $d_r$-cycle if and only if it
lifts to an element of $\pi_{t-s}(W_s^{s+r})$.  Given an element $\alpha
\in \pi_n(X)$ we shall let $filt_E(\alpha)$ denote its $E$-Adams filtration.

In what follows, the reader may find it helpful to assume that 
the spectra $W_s$ are CW
spectra and the maps
$$ W_s \rightarrow W_{s-1} $$
are the inclusions of subcomplexes.  If this is the case, then in what follows
the homotopy colimits may simply regarded as unions.  
As Bruner points out \cite[IV.3.1]{Hinfty}, this assumption represents
no loss of generality,
since any infinite tower may be replaced with a tower of inclusions of
CW-spectra through the use of CW approximation and mapping telescopes.
The Tate spectrum $\Sigma^{-1} tS^0$ is bifiltered, as depicted in the
following diagram.
$$ \xymatrix{
\vdots & \vdots & \vdots
\\
W_0(P^{N+1}) \ar[u] & W_1(P^{N+1}) \ar[u] \ar[l] & W_2(P^{N+1}) \ar[u]
\ar[l] & \cdots \ar[l]
\\
W_0(P^{N}) \ar[u] & W_1(P^{N}) \ar[u] \ar[l] & W_2(P^{N}) \ar[u]
\ar[l] & \cdots \ar[l]
\\
W_0(P^{N-1}) \ar[u] & W_1(P^{N-1}) \ar[u] \ar[l] & W_2(P^{N-1}) \ar[u]
\ar[l] & \cdots \ar[l]
\\
\vdots \ar[u] & \vdots \ar[u] & \vdots \ar[u]
} $$
In the above diagram, $W_k(P^N)$ is the
spectrum $W_k(P^N)_{-\infty} = \holim_M W_k(P^N_{-M})$, 
where the homotopy inverse limit is taken
\emph{after} smashing with $W_k$.  We emphasize that this is in general 
quite different from
what one obtains if one smashes with $W_k$ after taking the homotopy
inverse limit.

Given increasing sequences of integers
\begin{align*}
I & = \{ k_1 < k_2 < \cdots < k_l \} \\
J & = \{ N_1 < N_2 < \cdots < N_l \}
\end{align*}
we define subsets $S(I,J)$ of $\ZZ \times \ZZ$ by
$$ S(I,J) = \bigcup_{i=1}^{l} \{ (a,b) \: : \: a \ge k_i, b \le N_i \}. $$
We give the set of all multi-indices $(I,J)$ the structure of a poset by
declaring $(I,J) \le (I',J')$ if and only if $S(I,J) \subseteq S(I',J')$.

\begin{defn}[Filtered Tate spectrum]\label{defn:filttate}
Given sequences 
\begin{align*}
I & = \{ k_1 < k_2 < \cdots < k_l \} \\
J & = \{ N_1 < N_2 < \cdots < N_l \}
\end{align*}
with $k_i \ge 0$,
we define the filtered Tate spectrum (of the sphere) as the homotopy colimit
$$ W_I(P^J) = \bigsqcup_i W_{k_i}(P^{N_i}). $$
We allow for the possibility of $N_l = \infty$.
More generally, given another pair of sequences $(I',J') \le (I,J)$, we
define spectra
$$ W_I^{I'}(P^J_{J'}) = \mathrm{cofiber}\left( W_{I'+1}(P^{J'-1})
\rightarrow W_I(P^J) \right) $$
where $I'+1$ (respectively $J'-1$) is the sequence obtained by increasing 
(decreasing)
every element of the sequence by $1$.
\end{defn}

\begin{figure}
\begin{center}
\includegraphics[width=4in]{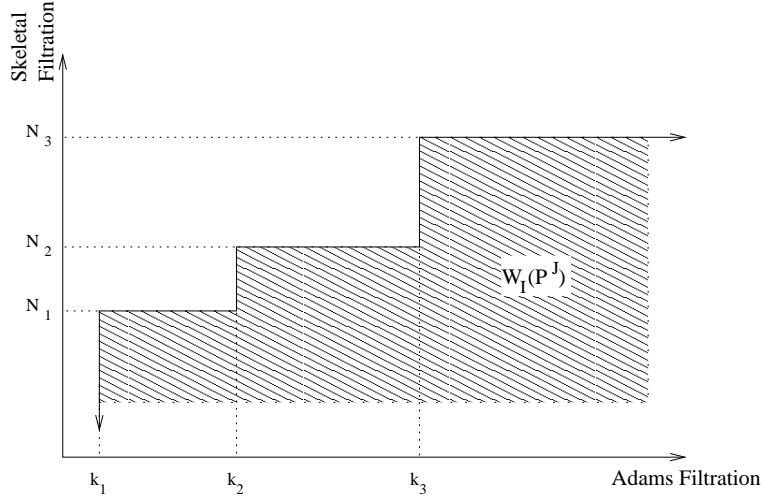}
\end{center}
\caption{The filtered Tate spectrum}
\label{fig:filt}
\end{figure}

Figure~\ref{fig:filt} displays a diagram of $S(I,J)$ 
intended to help the reader
visualize the filtered Tate spectrum.  The entire Tate spectrum is
represented by the right half-plane.  The shaded region $S(I,J)$ 
is the portion
represented by the filtered Tate spectrum $W_I(P^J)$.


\section{Definitions of various forms of the root invariant}\label{sec:defs}

In this section we will recall the definitions of the
Mahowald root invariant, the $E$-root invariant, and the algebraic root
invariant.  We will then define the filtered root
invariants.

\begin{defn}[{\it Root invariant} \rm]
Let $X$ be a finite complex, and $\alpha \in
\pi_t(X)$.  The \emph{root invariant} (also called the \emph{Mahowald
invariant}) 
of $\alpha$ is the coset 
of all dotted
arrows making the following diagram commute.
$$ \xymatrix{
S^{t} \ar@{.>}[r] \ar[d]^\alpha & \Sigma^{-N+1} X \ar[dd]  
\\
X \ar[d]
\\
tX \ar[r] & \Sigma P_{-N} \wedge X 
}$$
This coset is denoted $R(\alpha)$.  Here the map $X \rightarrow tX$
$p$-completion by the Segal conjecture (\ref{eq:segalconj}), 
and $N$ is chosen to be minimal
such that the composite $S^{t} \rightarrow \Sigma P_{-N} \wedge X$ 
is non-trivial.
\end{defn}

The root invariant is difficult to compute because it involves knowing the
homotopy groups of the finite complex $X$ and of $X \wedge P_{-N}$.  For
this reason, Mahowald and Ravenel \cite{MahowaldRavenel} 
introduced a calculable approximation
to the root invariant called the $E$-root invariant, for $E$ a ring
spectrum.  We find it useful to
generalize to arbitrary spectra and to elements of $\pi_*(E)$.

\begin{defn}[{\it $E$-root invariant} \rm]
Let $E$ be a spectrum.  Let $x \in \pi_t(E)$.  Define the 
\emph{$E$-root invariant} of $x$ to be the coset $R_E(\alpha)$ 
of dotted arrows making the following
diagram commute.
$$ \xymatrix{
S^{t} \ar@{.>}[r] \ar[d]_x & \Sigma^{-N+1} E \ar[dd] 
\\
E \ar[d] 
\\
tE \ar[r] & \Sigma E \wedge P_{-N}   
} $$
It is quite possible that the composite $S^{t} \rightarrow tE$ is trivial.  
If this is the case the
$E$-root invariant is said to be \emph{trivial}.  Otherwise, in the
diagram above we choose $N$ be minimal
such that the composite $S^t \rightarrow \Sigma E \wedge P_{-N}$ is
non-trivial.

If $E$ is a ring spectrum, and $\alpha \in \pi_*(X)$ for $X$ a finite
spectrum, let $x = h(\alpha) \in E_*(X)$ be the Hurewicz image of $\alpha$.
We will then refer to the $E \wedge X$-root invariant of $h(\alpha)$ simply 
as the 
$E$-root invariant of $\alpha$.  Therefore, by abuse of notation, we have
$$ R_E(\alpha) = R_{E \wedge X}(h(\alpha)). $$
\end{defn}

Finally, one may define root invariants in $\ext$.  These are called
algebraic root invariants.

\begin{defn}[{\it Algebraic root invariant} \rm]
Let $\alpha$ be an element of
$\ext^{s,t}(H_*X)$.
We have the following diagram of $\ext$ groups which defines the algebraic
root invariant $R_{alg}(\alpha)$.
$$\xymatrix@C+1em{ 
\ext^{s,t}(H_*X) \ar[d]_{f} \ar@{~>}[r]^{R_{alg}(-)} &
\ext^{s, t+N-1}(H_*X) \ar[d]_{\iota_N}
\\
\ext^{s,t-1}(H_*P_{-\infty} \wedge X) \ar[r]_{\nu_{N}} &
\ext^{s,t-1}(H_*P_{-N} \wedge X)
}$$
Here $f$ is induced by the inclusion of the $-1$-cell of
$P_{-\infty}$, $\nu_N$ is the 
projection onto the $-N$-coskeleton,
$\iota_N$ is inclusion of the $-N$-cell,
and $N$ is minimal with respect to the property that
$\nu_{N} \circ f(\alpha)$ is non-zero.  Then the algebraic root invariant
$R_{alg}(\alpha)$ is defined to be the coset of lifts 
$\gamma \in \ext^{s, t+N-1}(H_*X)$ of the
element $\nu_{N} \circ f (\alpha)$.
\end{defn}

We wish to extend these definitions to a sequence of filtered root
invariants that appear in the $E$-Adams resolution.  
Suppose that $X$ is a finite complex and $\alpha \in \pi_t(X)$.
We want to lift $\alpha$ over the smallest possible filtered Tate spectrum
(Definition~\ref{defn:filttate}).
To this end,
we shall describe a pair of sequences 
\begin{align*}
I & = \{ k_1 < k_2 < \cdots < k_l \} \\
J & = \{ -N_1 < -N_2 < \cdots < -N_l \}
\end{align*}
associated to $\alpha$, which we define inductively.  
Let $k_1 \ge 0$ be
maximal such that the composite
$$ S^{t-1} \xrightarrow{\alpha} \Sigma^{-1} X \rightarrow \Sigma^{-1} tX
\rightarrow W_0^{k_1-1}(P \wedge X)_{-\infty} $$
is trivial.
Next, choose 
$N_1$ to be maximal such that the composite
$$ S^{t-1} \xrightarrow{\alpha} \Sigma^{-1} X \rightarrow \Sigma^{-1} tX
\rightarrow W_0^{(k_1-1, k_1)}(P_{(-N_1+1, \infty)} \wedge X) $$
is trivial.  Inductively, given 
\begin{align*}
I' & = (k_1, k_2, \ldots, k_i) \\
J' & = (-N_1, -N_2, \ldots, -N_i)
\end{align*}
let $k_{i+1}$ be maximal so that the composite
$$ S^{t-1} \xrightarrow{\alpha} \Sigma^{-1} X \rightarrow \Sigma^{-1} tX
\rightarrow W_0^{(I'-1,k_{i+1}-1)}(P_{(J'+1,\infty)} \wedge X) $$
is trivial.  If there is no such maximal $k_{i+1}$, we declare that
$k_{i+1} = \infty$ and we are finished.  Otherwise, choose $N_{i+1}$ to be
maximal such that the composite
$$ S^{t-1} \xrightarrow{\alpha} \Sigma^{-1} X \rightarrow \Sigma^{-1} tX
\rightarrow W_0^{(I'-1,k_{i+1}-1, k_{i+1})}(P_{(J'+1,-N_{i+1}+1,\infty)} 
\wedge X) $$
is trivial, and continue the inductive procedure.  We shall refer to the
pair $(I,J)$ as the $E$-bifiltration of $\alpha$.

Observe that there is an exact sequence
$$ \pi_{t-1}(W_I(P^J \wedge X)) \rightarrow \pi_t(tX) \rightarrow
\pi_{t-1}(W^{I-1}(P_{J+1} \wedge X)). $$
Our choice of $(I,J)$ ensures that the image of $\alpha$ in 
$\pi_{t-1}(W^{I-1}(P_{J+1} \wedge X))$ is trivial.  Thus $\alpha$ lifts to an
element $f^\alpha \in \pi_{t-1}(W_I(P^J \wedge X))$.

\begin{defn}[{\it Filtered root invariants} \rm]
Let $X$ be a finite complex, let $E$ a ring spectrum such that the $E$-Adams
resolution converges, and let $\alpha$ be an element of $\pi_t(X)$ of
$E$-bifiltration $(I,J)$.  Given a lift 
$f^\alpha \in \pi_{t-1}(W_I(P^J \wedge X))$, the $k^\mathrm{th}$ filtered
root invariant is said to be trivial if $k \ne k_i$ for any $k_i \in I$.
Otherwise, if $k=k_i$ for some $i$, 
we say that the image $\beta$ of $f^\alpha$ under the collapse map
$$ \pi_{t-1}(W_I(P^J \wedge X)) \rightarrow
\pi_{t-1}(W_{k_i}^{k_i}(\Sigma^{-N_{i}} X)) $$
is an element of the \emph{$k^\mathrm{th}$ filtered root invariant} of 
$\alpha$.
The $k^\mathrm{th}$ filtered root invariant is the coset $R^{[k]}_E(\alpha)$ of
$E_1^{k,t+k+N_i-1}(X)$ of all such 
$\beta$ as we vary the 
lift $f^\alpha$.
\end{defn}

\begin{rmk}
Let $r_i$ denote the difference $k_{i+1}-k_i$.  Then there is a
factorization 
$$  \pi_{t-1}(W_I(P^J \wedge X)) \rightarrow
\pi_{t-1}(W_{k_i}^{k_{i+1}-1}(\Sigma^{-N_{k_i}} X)) \rightarrow
\pi_{t-1}(W_{k_i}^{k_i}(\Sigma^{-N_{k_i}} X)) $$
Thus any such $\beta \in R^{[k]}_E(\alpha)$ is actually a $d_r$-cycle 
for $r < r_i$.
\end{rmk}

\begin{figure}
\begin{center}
\includegraphics[width=4in]{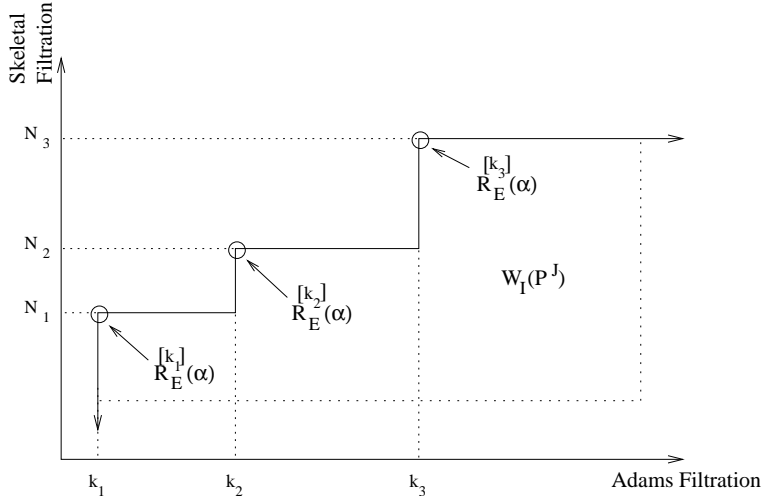}
\end{center}
\caption{The the bifiltration of the filtered root invariants in the
filtered Tate spectrum.}
\label{fig:filtroot}
\end{figure}

Figure~\ref{fig:filtroot} gives a companion visualization to
Figure~\ref{fig:filt}, by displaying the bifiltrations of the filtered root
invariants in the filtered Tate spectrum.


\section{The Toda bracket associated to a complex}\label{sec:toda}

In this section we introduce a variant of the Toda bracket.  This treatment
is essentially a specialization of the treatment of Toda brackets given in
the appendix of \cite{Shipley}.    
Suppose $K$ is a finite CW-spectrum with one bottom dimensional cell and
one top dimensional cell.  Suspend $K$ accordingly so that it
is connective and
$n$-dimensional with
one cell in dimension zero and 
one cell in dimension $n$.  For some other spectrum $X$ we shall define 
the \emph{$K$-Toda bracket}
to be an operator which, when defined, 
takes an element of $\pi_t(X)$ to a coset 
of $\pi_{t+n-1}(X)$.  
We shall present a dual definition, and 
show this dual definition is equivalent to
our original definition.
We will also define a variant on the $E_r$ term of the
$E$-ASS. 

In what follows, we let $K^j$ be the $j$-skeleton of $K$, and let $K_i^j$
be the quotient $K^j/K^{i-1}$.  We shall omit the top index  for the 
$i$-coskeleton $K_i = K/K^{i-1}$.

\begin{defn}[{\it $K$-Toda bracket} \rm]
Let
$$ f: \Sigma^{-1} K_1 \rightarrow S^0 $$
be the attaching map of the $1$-coskeleton of $K$ to the $0$-cell, 
so that the cofiber of $f$ is $K$.
Let $\nu : K_1 \rightarrow S^n$ be the projection onto the top cell.
Suppose $\alpha$ is an element of $\pi_t(X)$.
We have
$$ \pi_t(X) \xleftarrow{\nu_*} \pi_{t+n}(X \wedge K_1)
\xrightarrow{f_*} \pi_{t+n-1}(X). $$
We say the $K$-Toda bracket 
$$ \bra{K}(\alpha) \subseteq \pi_{t+n-1}(X) $$
is \emph{defined} if $\alpha$ is in the image of $\nu_*$.  Then the
$K$-Toda bracket is the 
collection of all $f_*(\gamma) \in \pi_{t+n-1}(X)$ where $\gamma \in 
\pi_{t+n}(X \wedge K_1)$ is any element satisfying $\nu_*(\gamma) =
\alpha$.
If $\bra{K}(\alpha)$ contains only 
one element, we will say that
it is \emph{strictly defined}.
\end{defn}

\begin{rmk}
If the $(n-1)$-skeleton $K^{n-1}$ is coreducible, then the attaching map $f$
factors as a composite
$$ f: \Sigma^{-1} K_1 \xrightarrow{\nu} S^{n-1} \xrightarrow{\xi} S^0. $$
Then the product
$\xi \cdot \alpha$ is in $\bra{K}(\alpha)$.
\end{rmk}

We remark on the relationship to Shipley's treatment of Toda brackets in
triangulated categories given in the appendix of \cite{Shipley}. 
In Shipley's terminology, if $K$ is an $m$-filtered object in 
$\{ f_1, f_2, \ldots, f_{m-1} \}$ 
with $F_1K \simeq S^0$ and 
$F_mK/F_{m-1}K \simeq S^n$, then
we have the equality (modulo indeterminacy)
$$ \bra{K}(\alpha) \capeq \bra{f_1, f_2, \ldots f_{m-1}, \alpha}. $$
In particular, if we take $K$ to be the $3$-filtered object
$$ \ast \subseteq S^0 \subseteq K^{n-1} \subseteq K $$
with attaching maps 
\begin{align*}
f_1: & \Sigma^{-1} K_1^{n-1} \rightarrow S^0 \\
f_2: & S^{n-1} \rightarrow K^{n-1}_1
\end{align*}
then we have the equality (modulo indeterminacy)
$$ \bra{K}(\alpha) \capeq \bra{f_1, f_2, \alpha}. $$
The reason we don't have exact equalities is that with the Toda bracket, you
can get more indeterminacy by varying the $m$-filtered object, whereas the
$m$-filtered object is fixed in the case of the $K$-Toda bracket.  However,
the $K$-Toda bracket can also get indeterminacy that is not seen by the
Toda bracket if $X \ne S^0$, because there are potentially 
more choices of lift
$\gamma$ if you smash with $X$.

In some of our applications the following dual variant will be more
natural to work with.  Our use of the word `dual' stems from the use the 
skeletal instead of
coskeletal filtration of $K$.

\begin{defn}[{\it Dual definition of the $K$-Toda bracket} \rm]
Let
$$ g: S^{n-1} \rightarrow K^{n-1} $$
be the attaching map of the $n$-cell of $K$ to the $(n-1)$-skeleton,
so that the cofiber of $g$ is $K$.
Let $\iota : S^0 \rightarrow K^{n-1}$ be the inclusion of the bottom cell.
Suppose $\alpha$ is an element of $\pi_t(X)$.
We have
$$ \pi_t(X) \xrightarrow{g_*} \pi_{t+n-1}(X \wedge K^{n-1})
\xleftarrow{\iota_*} \pi_{t+n-1}(X) $$
We say the (dual) $K$-Toda bracket
$$ \bra{K}(\alpha) \subseteq \pi_{t+n-1}(X) $$
is \emph{defined} if $g_*(\alpha)$ is in the image of $\iota_*$.  Then the
$K$-Toda bracket is the 
collection of all $\gamma \in \pi_{t+n-1}(X)$ where $\iota_*(\gamma) =
g_*(\alpha)$.
\end{defn}

We reconcile our use of the same notation in both definitions with the
following lemma.

\begin{lem}~\label{lem:todas}
The $K$-Toda bracket and the dual $K$-Toda bracket are equal.
\end{lem}

\begin{proof}
Let $\alpha$, $f$, $g$, $\iota$, and $\nu$ be the maps given in our two
definitions of the $K$-Toda bracket.
For the purposes of this proof, we shall refer to the dual $K$-Toda bracket
as $\bra{K}^d(\alpha)$.  The map of cofiber sequences makes the
following diagram commute.
$$ \xymatrix{
S^0 \ar[r] \ar[d]_{\iota} & K \ar[r] \ar@{=}[d] & 
K_1 \ar[r]^{f} \ar[d]_\nu & S^1 \ar[d]^\iota
\\
K^{n-1} \ar[r] & K \ar[r] & S^n \ar[r]_-{g} & \Sigma K^{n-1}
}$$
Taking the last square in the above triangle, extending it to a map of
cofiber sequences the other way, and smashing with $X$, we get the
following commutative diagram whose columns are exact.
$$ \xymatrix{
\pi_{t+n}(X \wedge K^{n-1}_1) \ar[d]_{\iota'_*} \ar@{=}[r] & 
\pi_{t+n}(X \wedge K_1^{n-1}) \ar[d]^{f'_*} 
\\
\pi_{t+n}(X \wedge K_1) \ar[d]_{\nu_*} \ar[r]^{f_*} & 
\pi_{t+n-1}(X) \ar[d]^{\iota_*} 
\\
\pi_t(X) \ar[d]_{g'_*} \ar[r]_-{g_*}
  &
\pi_{t+n-1}(X \wedge K^{n-1}) \ar[d]^{\nu'_*}
\\
\pi_{t+n-1}(X \wedge K_1^{n-1}) \ar@{=}[r] &
\pi_{t+n-1}(X \wedge K_1^{n-1})
} $$
Here $\iota'$ is the inclusion and $\nu'$ is the projection, and $f'$ and
$g'$ are determined from the diagram.
The bracket $\bra{K}(\alpha)$ is computed by taking the preimage under 
$\nu_*$ and
then applying $f_*$.  The bracket $\bra{K}^d(\alpha)$ is computed applying
$g_*$ and then taking the preimage over $\iota_*$.
We see that for any lift $\gamma \in \pi_{t+n}(X \wedge K_1)$ of $\alpha$,
$f_*(\gamma)$ is a lift of $g_*(\alpha)$.  It follows that we have
$\bra{K}(\alpha) \subseteq \bra{K}^d(\alpha)$.

For the reverse containment, suppose that $\beta \in \pi_{t+n-1}(X)$ is a
lift of $g_*(\alpha)$.  
Then $\beta$ is an arbitrary element of $\bra{K}^d(\alpha)$.
Let $\gamma \in \pi_{t+n}(X \wedge K_1)$ be any
lift if $\alpha$.  Such a lift exists since 
$$ g'_*(\alpha) = \nu'_* g_*(\alpha) = \nu'_* \iota_*(\beta) = 0. $$
We must add a correction term to $\gamma$ since it is
not necessarily the case that $f_*(\gamma) = \beta$.  Let $\delta$ be the
difference $f_*(\gamma) - \beta$.  Then $\iota_*(\delta) = 0$, so there is
a lift to $\td{\delta} \in \pi_{t+n}(X \wedge K_1^{n-1})$.  
Then $\gamma' = \gamma - \iota'_* (\td{\delta})$ 
is another lift of $\alpha$, and $f_*(\gamma) =
\beta$.  We have therefore proven that $\bra{K}^d(\alpha) \subseteq
\bra{K}(\alpha)$.
\end{proof}

\begin{rmk}
Perhaps a more 
conceptual way to prove Lemma~\ref{lem:todas} would be to regard the
$K$-Toda bracket as a $d_n$ in the AHSS for $K$ with respect to the
coskeletal filtration, and the dual $K$-Toda bracket as a $d_n$ in the AHSS
with respect to the skeletal filtration.  There is a comparison of these
spectral sequences which is an isomorphism on $E_1$-terms, hence the
spectral sequences are isomorphic, and in particular they have the same 
differentials $d_n$.
See, for instance, Appendix B of \cite{GreenleesMay}.
\end{rmk}

We need to extend our $K$-Toda brackets to operations on the $E_r$ term of
the $E$-ASS.  Suppose that the attaching map
$f: \Sigma^{-1} K_1 \rightarrow S^0$ has $E$-Adams filtration $d$.

\begin{defn}[{\it $K$-Toda brackets on $E_r(X)$} \rm]
Suppose that $\alpha$ is an element of $E^{s,t}_r(X)$.  Choose a lift of
$\alpha$ to $\td{\alpha} \in \pi_{t-s}(W_s^{s+r-1}(X))$.  Let $\td{f}:
\Sigma^{-1} K_1 \rightarrow W_d(S^0)$ be a lift of $f$.  Then we have
$$ \pi_{t-s}(W_s^{s+r-1}(X)) \xleftarrow{\nu_*} \pi_{t-s+n}(W_s^{s+r-1}(X 
\wedge K_1))
\xrightarrow{\td{f}_*} \pi_{t+n-1}(W_{s+d}^{s+d+r-1}(X)). $$
If there exists a lift $\td{\alpha}$ which is in the image of $\nu_*$, then
we say that the $K$-Toda bracket is \emph{defined}.  Take $\gamma \in 
\pi_{t-s+n}(W_s^{s+r-1}(X))$ such that $\nu_*(\gamma) = \td{\alpha}$.  Then
the image of $\td{f_*}(\gamma)$ in $E_r^{s+d,t+d+n-1}(X)$ is in the
$K$-Toda bracket $\bra{K}(\alpha)$.  We define the $K$-Toda bracket to be 
the set of all
such images for various choices of $\td{\alpha}$, $\gamma$, and $\td{f}$.
We shall say that the $K$-Toda bracket has \emph{$E$-Adams degree $d$}.
\end{defn}

\begin{rmk}
If the $K$-Toda bracket of $\alpha$ is defined, where $\alpha \in
E^{s,t}_\infty(X)$, then if $\alpha$ detects $\br{\alpha}$, every element of 
$\bra{K}(\alpha)$ detects an element of $\bra{K}(\br{\alpha})$.  In
particular, if $\bra{K}(\alpha)$ contains $0$, then there is an element of
$\bra{K}(\td{\alpha})$ of $E$-Adams filtration greater than $k+d$.
\end{rmk}


\section{Statement of results}\label{sec:results}

In this section we will state our main results concerning filtered root
invariants.  The proofs of many of the theorems are given in
Section~\ref{sec:proofs}.
Throughout this section let $E$ be a ring spectrum, 
let $X$ be a finite complex, and suppose 
$\alpha$ is an element of $\pi_t(X)$ of $E$-bifiltration $(I,J)$ with
\begin{align*}
I & = (k_1, k_2, \ldots, k_l) \\
J & = (-N_1, -N_2, \ldots, -N_l)
\end{align*}

\begin{thm}[{\it Relationship to homotopy root invariant} \rm]\label{thmA}
Suppose that $R^{[k_i]}_E(\alpha)$ contains a permanent cycle $\beta$. Then
there exists an element
$\br{\beta} \in \pi_*(X)$ which $\beta$ detects such that 
the following diagram commutes up to elements of $E$-Adams
filtration greater than or equal to $k_{i+1}$.
$$
\xymatrix{
S^{t} \ar[d]_{\alpha} \ar[r]_-{\br{\beta}} & \Sigma^{-N_i+1} X \ar[dd]  \\
X \ar[d] \\
tX \ar[r] & \Sigma P_{-N_i} \wedge X  
} $$
\end{thm}

The proof of Theorem~\ref{thmA} is deferred to Section~\ref{sec:proofs}.
If $k_{i+1} = \infty$, then $i$ is equal to $l$ (the maximal index of the
bifiltration), and
we get the following corollary.

\begin{cor}\label{corA}
The top filtered root invariant $R^{[k_l]}_E(\alpha)$
is contained in $E_\infty(X)$.  Let $\br{R}^{[k_l]}_E(\alpha)$ be the set of
elements of $\pi_*(X)$ detected by $R^{[k_l]}_E(\alpha)$.
There are two possibilities:
\begin{enumerate}
\item The image of the elements of $\br{R}^{[k_l]}_E(\alpha)$
in $\pi_*(P_{-N_l}\wedge X)$ under the inclusion of the
bottom cell is zero, and the homotopy root invariant lies in a higher stem
than the $k_l^\mathrm{th}$ filtered root invariant.

\item There is an equality (modulo indeterminacy)
$$ \br{R}^{[k_l]}_E(\alpha) \capeq R(\alpha) $$
\end{enumerate}
\end{cor}

\begin{proof}
Since $k_l = \infty$, the diagram in Theorem~\ref{thmA} commutes modulo
elements of infinite Adams filtration.
The $E$-Adams resolution was assumed to converge, so this means 
the diagram actually
commutes.  Let $\br{\beta}$ be the element described in
Theorem~\ref{thmA}.  If the image of $\br{\beta}$ in $\pi_{t-1}(P_{-N_i})$
is zero, then the homotopy root invariant lies in a higher stem than the
$i^\mathrm{th}$ filtered root invariant.  Otherwise, $\br{\beta}$ is an
element of $R(\alpha)$.
\end{proof}

Unfortunately, Corollary~\ref{corA} is difficult to invoke in practice.
This is because given a filtered root invariant, one usually does not know 
whether it is the highest one or not.
In practice we are in the situation where we have a permanent cycle in a
filtered root invariant and 
we would like to show that the diagram of Theorem~\ref{thmA}
commutes on the nose.  One strategy is to write out the
Atiyah-Hirzebruch spectral sequence for $\pi_*(P_{-N_i})$ and try to
show there are no elements of higher $E$-Adams filtration which
could be the difference between the image of $\alpha$ and the image of
$\br{\beta}$ in $\pi_*(P_{-N_i})$.  This method is outlined in 
Procedure~\ref{procA}.

We now present some theorems which relate filtered root invariants to
differentials and compositions in the Adams spectral sequence.  Recall
that $r_j = k_{j+1}-k_j$.

\begin{thm}[{\it Relationship to Adams differentials} \rm]\label{thmB}
Suppose that the $P_{-N_i}^{-N_{i+1}}$-Toda bracket has $E$-Adams 
degree $d$
and that $d \le r_{i+1}$.
Then the following is true.
\begin{enumerate}

\item
$\bra{P_{-N_i}^{-N_{i+1}}}(R^{[k_{i+1}]}_E(\alpha))$
is defined and contains a permanent cycle.

\item
$R^{[k_i]}_E (\alpha)$ consists of elements
which are $d_r$ cycles for $r < r_i+d$.

\item There is a containment
$$ d_{r_i+d} R^{[k_i]}_E(\alpha) \subseteq 
\bra{P_{-N_i}^{-N_{i+1}}}(R^{[k_{i+1}]}_E(\alpha))$$
where both elements are thought of as elements of $E^{*,*}_{r_i+d}(X)$.

\end{enumerate}
\end{thm}

The proof of Theorem~\ref{thmB} is deferred to Section~\ref{sec:proofs}.
We intend to use Theorem~\ref{thmB} in reverse: given the $i^\mathrm{th}$
filtered root invariant, we would like to deduce the $(i+1)^\mathrm{st}$
filtered root invariant from the presence of an Adams differential.
If the differential in Theorem~\ref{thmB} is zero, we still may be able to
glean some information from the presence of a non-trivial composition.

\begin{thm}[{\it Relationship to compositions} \rm]\label{thmC}
Suppose that 
$R^{[k_i]}_E(\alpha)$ is a coset of permanent cycles.  Let
$\td{R}^{[k_i]}_E(\alpha)$ denote the coset of all lifts of elements of
$R^{[k_i]}_E(\alpha)$ to $\pi_{t+N_i-1}(W_{k_i}(X))$.  When defined,
there are lifts
of the Toda brackets
$\bra{P_{-m}^{-N_i}}(\td{R}_E^{[k_i]}(\alpha))$
to
$$ \td{\bra{P_{-m}^{-N_i}}}(\td{R}_E^{[k_i]}(\alpha)) \subseteq
\pi_{t+m-2}(W_{k_{i+1}}(X)). $$
Let $M$ be the minimal such $m > N_i$ with the property that 
$\td{\bra{P_{-m}^{-N_i}}}(\td{R}_E^{[k_i]}(\alpha))$ contains a non-trivial
element.  Let $d$ be the $E$-Adams degree of the Toda bracket 
$\bra{P_{-N_i}^{-N_{i+1}}}(-)$,
and suppose that $d \le r_{i+1}$.  Then the following is true.
\begin{enumerate}

\item  Let
$$ \br{R}^{[k_i]}_E(\alpha) \subseteq \pi_{t+N_i-1}(X) $$
denote the elements which are detected by the permanent cycles 
of $R^{[k_i]}_E(\alpha)$.  Then 
the Toda bracket 
$$ \bra{P_{-M}^{-N_i}}(\br{R}^{[k_i]}_E(\alpha)) \subseteq \pi_{t+M-2}(X)$$
is defined.

\item The Toda bracket $\bra{P_{-M}^{-N_{i+1}}}(R^{[k_{i+1}]}_E(\alpha))$ 
is defined in
$E_{r_{i+1}}(X)$, and contains a permanent cycle.  We shall denote the
collection of all elements which are detected by these permanent cycles by
$$ \br{\bra{P_{-M}^{-N_{i+1}}}(R^{[k_{i+1}]}_E(\alpha))} 
\subseteq \pi_{t+M-2}(X). $$

\item There is an equality (modulo indeterminacy)
$$ \bra{P_{-M}^{-N_i}}(\br{R}^{[k_i]}_E(\alpha)) \capeq 
\br{\bra{P_{-M}^{-N_{i+1}}}(R^{[k_{i+1}]}_E(\alpha))}. $$

\end{enumerate}
\end{thm}

\begin{rmk}
The hypothesis $d \le r_{i+1}$ in Theorems~\ref{thmB} and \ref{thmC}
is a necessary technical hypothesis to make
the proofs work.  In practice, $d$ is often equal to $1$.  Since $r_j$ is
always positive, the hypothesis is satisfied in this case.
\end{rmk}

Finally, we give a partial description of the first filtered
root invariant.  We begin with a simple observation.

\begin{lem}\label{filtlem}
If $filt_E(\alpha) = k$, then $R^{[s]}_E(\alpha)$ is trivial for $s < k$.
\end{lem}

\begin{proof}
We must show that in the $E$-bifiltration of $\alpha$, $k_1 \ge k$.
Consider the following diagram.
$$ \xymatrix{
\pi_t(X) \ar[r] \ar[d] & \pi_t(W_0^{k-1}(X)) \ar[d] \\
\pi_t(tX) \ar[r] & \pi_{t-1}(W_0^{k-1}(P_{-\infty}))
} $$
Since $filt_E(\alpha) = k$, the image of $\alpha$ in $\pi_t(W_0^{k-1}(X))$
is trivial.  Therefore the image of $\alpha$ in 
$\pi_{t-1}(W_0^{k-1}(P_{-\infty}))$ is trivial.  By the maximality of
$k_1$, we have $k_1 \ge k$.
\end{proof}

If $filt_E(\alpha) = 0$, (i.e. when $\alpha$ has a non-trivial Hurewicz
image) then we can sometimes identify the first
non-trivial filtered root invariant with the $E$-root invariant.

\begin{prop}[{\it Relationship to 
the $E$-root invariant} \rm]\label{prop:Eroot}
The $E$-root invariant $R_E(\alpha)$ is non-trivial if and only if 
$R^{[0]}_E(\alpha)$
is non-trivial.  If this is the case, regarding $R_E(\alpha)$ as
being contained in $E^{0,*}_1(X)$, we have
$$ R^{[0]}_E(\alpha) \subseteq R_E(\alpha). $$
\end{prop}

\begin{proof}
This is immediate from the definitions.
\end{proof}

If $filt_E(\alpha) = 1$, then the $E$-root invariant is trivial and 
Lemma~\ref{filtlem} implies that the zeroth
filtered root invariant is trivial.  We can however sometimes compute the
first filtered root invariant using the $E\wedge \br{E}$-root invariant.

\begin{prop}[\it{Relationship to the $E\wedge \br{E}$-root invariant} 
\rm]\label{prop:E^Eroot}
Suppose that $\alpha$ has $E$-Adams filtration $1$.
Then there exists an element $\td{\alpha} \in \pi_{t}(E \wedge \br{E} 
\wedge X^{\wedge}_p) = E^{1,t+1}_1(X^{\wedge}_p)$ 
which detects $\alpha$ in the $E$-ASS, 
and such such that 
$R_{E \wedge \br{E}}(\td \alpha)$ is trivial if and
only if $R^{[1]}_E(\alpha)$ is trivial.  
There is a containment
$$ R^{[1]}_{E}(\alpha) \subseteq R_{E\wedge \br{E}}(\td{\alpha}). $$
\end{prop}

In practice, we will not know what choice of detecting element $\td{\alpha}$
to choose, so the following corollary will prove useful.

\begin{cor}\label{cor:E^Eroot}
Suppose that $filt_E(\alpha) = 1$.
If $\td{\alpha} \in \pi_t(E \wedge \br{E} \wedge X)$ 
is any element which detects $\alpha$ in the $E$-ASS, then
$$ R^{[1]}_{E}(\alpha) \subseteq R_{E \wedge \br{E}}(\td{\alpha}) + A. $$
Here $A$ is the image of the map
$$ \pi_{t}(W_0^{(0,1)}(P_{(-N,-N+1)}\wedge X)) \xrightarrow{\partial} 
\pi_{t-1}(W_1^1(\Sigma^{-N} X)) $$
where $\partial$ is the boundary homomorphism associated to the cofiber
sequence
$$ W_1^1(\Sigma^{-N}X) \xrightarrow{\iota} W_0^1(P_{-N} \wedge X) \rightarrow 
W_0^{(0,1)}(P_{(-N,-N+1)}\wedge X) $$
and the $E\wedge\br{E}$-root invariant of $\td{\alpha}$ is carried by the
$-N$-cell of $P_{-\infty}$.
\end{cor}

The proof of Proposition~\ref{prop:E^Eroot} and Corollary~\ref{cor:E^Eroot}
is deferred to
Section~\ref{sec:proofs}.
When $E = H\FF_p = H$, 
we can identify the first
filtered root invariant with the algebraic root invariant.

\begin{thm}[{\it Relationship to the algebraic root invariant} \rm]\label{thmD}
If $E$ is the Eilenberg-MacLane spectrum $H\FF_p = H$ and $\alpha$ has 
Adams filtration $k$, then 
$k_1 = k$.  Furthermore, the filtered root invariant 
$R^{[k]}_{H}(\alpha)$ consists of $d_1$ cycles
which detect a coset of non-trivial elements 
$\br{R}^{[k]}_{H}(\alpha) \subseteq E_2^{k, t+k+N_1-1}(X)$, and there
exists a choice of $\td{\alpha} \in E_2^{k, t+k}(X)$ which 
detects $\alpha$ in the
ASS such that  
$$ \br{R}^{[k]}_{H}(\alpha) \subseteq R_{alg}(\td{\alpha}). $$
\end{thm}


The proof of Theorem~\ref{thmD} is deferred to Section~\ref{sec:proofs}.
We have given a partial scenario as to how the filtered root
invariants can be used to calculate root invariants using the $E$-ASS.  
One first calculates the zeroth filtered root invariant as an $E$-root
invariant or the first non-trivial root invariant as an 
algebraic root invariant.
Then the idea, while
only sometimes correct, is that ``if a filtered root invariant does
not detect the root invariant, then it either supports a
differential or a composition that points to the next 
filtered root invariant.''  The last filtered root invariant then has a
chance of detecting the homotopy root invariant.  
We stress that many things can interfere
with this actually happening.  


\section{Proofs of the main theorems}\label{sec:proofs}

In all of the proofs below, we shall assume that our finite complex $X$
is actually $S^0$.  The general case is no different, but smashing
everything with $X$ complicates the notation.

\begin{proof}[Proof of Theorem~\ref{thmA}]
Let $\beta$ be an element of $R^{[k_i]}_E(\alpha)$.  Then there exists a
lift $f^\alpha$ such that $\beta$ is the image of $f^\alpha$ under the
collapse map
$$ \pi_{t-1}(W_I(P^J)) \rightarrow \pi_{t-1}(W_{k_i}^{k_{i+1}-1}(S^{-N_i})).
$$
Consider the following diagram.
$$ \xymatrix@R-1em@C-1em{
& S^{t-1} \ar@{.>}[ddddl]_{\br{\beta}} \ar[d]_\alpha
\ar@{.>}[ddddr]^{\td{\beta}}
\ar[ddrr]^{f^\alpha} 
\\
& S^{-1} \ar[d] 
\\
& P_{-\infty} \ar[d]_{\br{\nu}} & & W_I(P^J) 
\ar[d]_{\nu} 
\\
& P_{-N_i} & & W_{k_i}(P_{-N_i}) \ar[ll] \ar[dd]_\eta 
& W_{k_{i+1}}(P_{-N_i}) \ar[l]
\\
S^{-N_i} \ar[ur]_{\br{\iota}} & & W_{k_i}(S^{-N_i}) \ar[ll] 
\ar[dd]_\eta \ar[ur]_{\iota}
\\
& & & W_{k_i}^{k_{i+1}-1}(P_{-N_i}) 
\\
& & W_{k_i}^{k_{i+1}-1}(S^{-N_i}) \ar[ur]_\iota
} $$
In the above diagram, $\nu$ is induced from the composite
$$ W_I(P^J) \rightarrow W_{I'}(P^{J'}_{-N_i}) \rightarrow W_{k_i}(P_{-N_i})
$$
where $I' = (k_i, \ldots, k_l)$ and $J'=(-N_i, \ldots, -N_l)$.
Since $R^{[k_i]}_{E}(\alpha)$ contains a permanent cycle $\beta$, 
there exists a map $\td{\beta}$ (as above) 
such that $\eta \td{\beta}$ projects to $\beta$, and such
that
$\iota \eta \td{\beta} = \eta \nu f^\alpha$.  
If $\delta = \iota \td{\beta} - \nu f^\alpha$, then $\delta$ lifts
to $\pi_{t-1}(W_{k_{i+1}}(P_{-N_i}))$.
Let $\br{\beta}$ be the map induced by $\td{\beta}$, and denote by 
$\br{\delta} \in
\pi_{t-1}(P_{-N_i})$ the image of $\delta$.  Then $filt_E(\br{\delta}) \ge
k_{i+1}$ and have the following formula.
$$ \br{\iota} \circ \br{\beta} = \br{\nu} \circ \alpha + \br{\delta} $$
This is precisely what we wanted to prove.
\end{proof}

\begin{proof}[Proof of Theorem~\ref{thmB}]
Fix a lift $f^\alpha \in \pi_{t-1}(W_I(P^J))$ of $\alpha$.  Define a
spectrum
$$ U = W_{(k_i, k_{i+1}, k_{i+2})}(P_{-N_i}^{(-N_i, -N_{i+1}, \infty)}). $$
There is a natural map
$$ W_I(P^J) \rightarrow U $$
and let $\gamma \in \pi_{t-1}(U)$ be the image of $f^\alpha$ under this
map.

Since we have assumed that the $P_{-N_i}^{-N_{i+1}}$-Toda bracket has
$E$-Adams degree $d$, there is a
lift of the attaching map 
$$ f: \Sigma^{-1} P_{-N_i+1}^{-N_{i+1}} \rightarrow S^{-N_i} $$
to a map
$$ \td{f}: \Sigma^{-1} P_{-N_i+1}^{-N_{i+1}} \rightarrow W_d(S^{-N_i}) $$
Define a filtered stunted projective space
$(P_{-N_i}^{-N_{i+1}})_{[d]}$ 
by the following cofiber sequence.
$$ \Sigma^{-1} P_{-N_i+1}^{-N_{i+1}} \xrightarrow{\td{f}} 
W_d(S^{-N_i}) \rightarrow (P_{-N_i}^{-N_{i+1}})_{[d]} $$
The spectrum $U$ is given by the homotopy pushout 
$$ \left( W_{k_{i+2}}(P_{-N_i})
\sqcup_{W_{k_{i+2}}(P_{-N_i}^{-N_{i+1}})}
W_{k_{i+1}}(P_{-N_i}^{-N_{i+1}}) \right) 
\sqcup_{W_{k_{i+1}}(S^{-N_i})} 
W_{k_i}(S^{-N_i}). $$
Since we have assumed that $d \le r_{i+1} = k_{i+2} - k_{i+1}$, the spectrum 
$U$ admits the 
equivalent description as
$$ \left( W_{k_{i+2}}(P_{-N_i})
\sqcup_{W_{k_{i+2}}(P_{-N_i}^{-N_{i+1}})}
W_{k_{i+1}}((P_{-N_i}^{-N_{i+1}})_{[d]}) \right) 
\sqcup_{W_{k_{i+1}+d}(S^{-N_i})} 
W_{k_i}(S^{-N_i}). $$
 
Consider the following commutative diagram.
\customsizeA{
$$ \xymatrix@C-4em@R+.5em{
& \pi_{t-2}(W_{k_{i+1}+d}^{k_{i+2}+d-1}(S^{-N_i})) 
\\ 
& &
\pi_{t-2}(W_{k_{i+1}+d}(S^{-N_i})) \ar[lu]_{\eta} \ar[rd] 
\\
& \pi_{t-1}(W_{k_{i+1}}^{k_{i+2}-1}(P_{-N_i+1}^{-N_{i+1}})) 
\ar[uu]_{\td{f}} \ar[dl]^\nu & &
\pi_{t-2}(W_{k_{i+1}}(S^{-N_i})) 
\\
\pi_{t-1}(W_{k_{i+1}}^{k_{i+2}-1}(S^{-N_{i+1}})) 
\ar@{~>}@/^2pc/[ruuu]^{\bra{P_{-N_i}^{-N_{i+1}}}(-)} & &
\pi_{t-1}(U) \ar[ll] \ar[ul]^{g_2} \ar[uu]_{\partial_{tot}}
\ar[r]_-{g_3} &
\pi_{t-1}(W_{k_i}^{k_{i+1}-1}(S^{-N_i})) \ar[u]_\partial
} $$}
Here $g_3$ is obtained by collapsing out the first and second factors of
the homotopy pushout $U$, and similarly $g_2$ is obtained by
collapsing out the first and third factors.  The map $\partial_{tot}$ is the
boundary homomorphism of the Meyer-Vietoris sequence, and may be thought
of as collapsing out all of the factors of the homotopy pushout.  
The wavy arrow indicates that the
$P_{-N_i}^{-N_{i+1}}$-Toda bracket is taken by taking the inverse image
under $\nu$, followed by application of $\td{f}$.

Our element $\gamma \in \pi_{t-1}(U)$ has compatible images in all of the
other groups in the diagram.  The image of $\gamma$ under $g_3$ projects to
an element of $R^{[k_i]}_E(\alpha)$ in $E_{r_i}$.  Following the outside of
the diagram from the lower right-hand corner to the top corner 
counter-clockwise 
amounts to taking
$d_{r_i+d}$ in the $E$-ASS.  Following from the lower right-hand corner to
the top clockwise applies the Toda bracket to the $k_{i+1}^\mathrm{st}$
filtered root invariant.  Thus $d_{r_i+d}R^{[k_i]}_E(\alpha)$ and
$\bra{P_{-N_i}^{-N_{i+1}}}(R^{[k{i+1}]}_E(\alpha))$ have a common element.
\end{proof}

\begin{proof}[Proof of Theorem~\ref{thmC}]
Fix a lift $f^\alpha \in \pi_{t-1}(W_I(P^J))$ of $\alpha$.  Let $V_m$ be
defined by
$$ V_m = W_{(0,k_i,k_{i+1},k_{i+2})}(P^{(-N_{i-1}, -N_i,
-N_{i+1},\infty)}_{-m}) $$ 
Then $V_m$ is defined by the following homotopy pushout.
\begin{gather*}
P_{-m}^{-N_{i-1}} 
\sqcup_{W_{k_i}(P_{-m}^{-N_{i-1}})}
 W_{k_i}(P_{-m}^{-N_i})
\sqcup_{W_{k_{i+1}}(P_{-m}^{-N_i})}
\\
 W_{k_{i+1}}(P_{-m}^{-N_{i+1}}) 
\sqcup_{W_{k_{i+2}}(P_{-m}^{-N_{i+1}})}
W_{k_{i+2}}(P_{-m}) 
\end{gather*}
There is a natural map
$$ W_I(P^J) \rightarrow V_m $$
and we define $\gamma_m \in \pi_{t-1}(V_m)$ to be the image of $f^\alpha$
under this map.
We need to lift $\gamma_m$ a little more.  The composite
$$
\pi_{t-1}(V_m) \xrightarrow{\partial} 
\pi_{t-2}(W_{k_{i+1}}(P_{-m}^{-N_i})) \rightarrow 
\pi_{t-2} (W_{k_{i+1}}(S^{-N_i}))
$$
sends $\gamma_m$ to zero, since it carries the $E$-Adams differential of
an element in $R^{[k_i]}_E(\alpha)$,
which was hypothesized to be zero.  
Here $\partial$ is a Meyer-Vietoris boundary
homomorphism.
Thus our element $\gamma_m$ lifts to
an element $\td{\gamma}_m \in \pi_{t-1}(\td{V}_m)$ where $\td{V}_m$ is
the following spectrum.
\begin{gather*} 
P_{-m}^{-N_{i-1}} 
\sqcup_{W_{k_i}(P_{-m}^{-N_{i-1}})}
W_{k_i}(P_{-m}^{-N_i}) 
\sqcup_{W_{k_{i+1}}(P_{-m}^{-N_i-1})}
\\
W_{k_{i+1}}(P_{-m}^{-N_{i+1}}) 
\sqcup_{W_{k_{i+2}}(P_{-m}^{-N_{i+1}})}
W_{k_{i+2}}(P_{-m})
\end{gather*}
The following diagram explains the lifted bracket 
$\td{\bra{P_{-m}^{-N_i}}}(\td{R}_E^{[k_i]}(\alpha))$.
\begin{equation}\label{diag:lifted_toda}
\xymatrix{
\pi_{t-1}(\td{V}_m) \ar[r]^-{\partial} \ar[d]_{p_2} &
\pi_{t-2}(W_{k_{i+1}}(P_{-m}^{-N_i-1})) \ar[d] &
\pi_{t-2}(W_{k_{i+1}}(S^{-m})) \ar[l] \ar[d]
\\
\pi_{t-1}(W_{k_i}(S^{-N_i})) \ar[r]^g
\ar@/_2pc/@{~>}[rr]_{\bra{P_{-m}^{-N_i-1}}(-)} &
\pi_{t-2}(W_{k_i}(P_{-m}^{-N_i-1})) &
\pi_{t-2}(W_{k_i}(S^{-m})) \ar[l]
}
\end{equation}
Here $p_2$ is projection onto the second factor of the homotopy pushout
$\td{V}_m$, $\partial$ is a Meyer-Vietoris boundary, and $g$ is the attaching
map of the $-N_i$-cell to $P^{-N_i-1}_{-m}$.
Assuming (inductively) that for every $N_i < m' < m$, 
$$ \td{\bra{P_{-m'}^{-N_i}}}(\td{R}_E^{[k_i]}(\alpha)) = 0 $$
we there is a lift of $\partial(\td{\gamma}_m)$ to
$\pi_{t-2}(W_{k_{i+1}}(S^{-m})$.  The set of all such lifts is defined to
be
$$ \td{\bra{P_{-m}^{-N_i}}}(\td{R}_E^{[k_i]}(\alpha)) \subseteq 
\pi_{t-2}(W_{k_{i+1}}(S^{-m})) $$
Diagram~\ref{diag:lifted_toda} implies that this set of lifts is indeed a 
lift of the Toda bracket
$$ \bra{P_{-m}^{-N_i}}(\td{R}_E^{[k_i]}(\alpha)) \subseteq 
\pi_{t-2}(W_{k_i}(S^{-m})). $$
Let $M$ be the first $m$ such that 
$\td{\bra{P_{-m}^{-N_i}}}(\td{R}_E^{[k_i]}(\alpha))$ contains a non-trivial
element.  
Let $\td{\gamma} = \td{\gamma}_M$ and let $\td{V} = \td{V}_M$.
The image of $\td{\gamma}$ under the composite
$$ \pi_{t-1}(\td{V}) \xrightarrow{\partial}
\pi_{t-2}(W_{k_{i+1}}(P_{-M}^{-N_i-1}))
\rightarrow \pi_{t-2}(W_{k_{i+1}}(P_{-M+1}^{-N_i-1})) $$
is zero, and $\td{\gamma}$ lifts even further, to an element
$\td{\td{\gamma}} \in
\pi_{t-1}(\td{\td{V}})$,
where $\td{\td{V}}$ is defined to be the following spectrum.
\begin{gather*}
P_{-M}^{-N_{i-1}} 
\sqcup_{W_{k_i}(P_{-M}^{-N_{i-1}})}
W_{k_i}(P_{-M}^{-N_i}) 
\sqcup_{W_{k_{i+1}}(S^{-M})}
\\
W_{k_{i+1}}(P_{-M}^{-N_{i+1}}) 
\sqcup_{W_{k_{i+2}}(P_{-M}^{-N_{i+1}})}
W_{k_{i+2}}(P_{-M})
\end{gather*}
Let $f$ be the attaching map
$$ f: \Sigma^{-1} P_{-M+1}^{-N_{i+1}} \rightarrow S^{-M}. $$
Since $\bra{P_{-M}^{-N_{i+1}}}$ has $E$-Adams degree $d$, there is a lift
$$ \td{f} : \Sigma^{-1} P_{-M+1}^{-N_{i+1}} \rightarrow W_d(S^{-M}). $$

\begin{landscape}
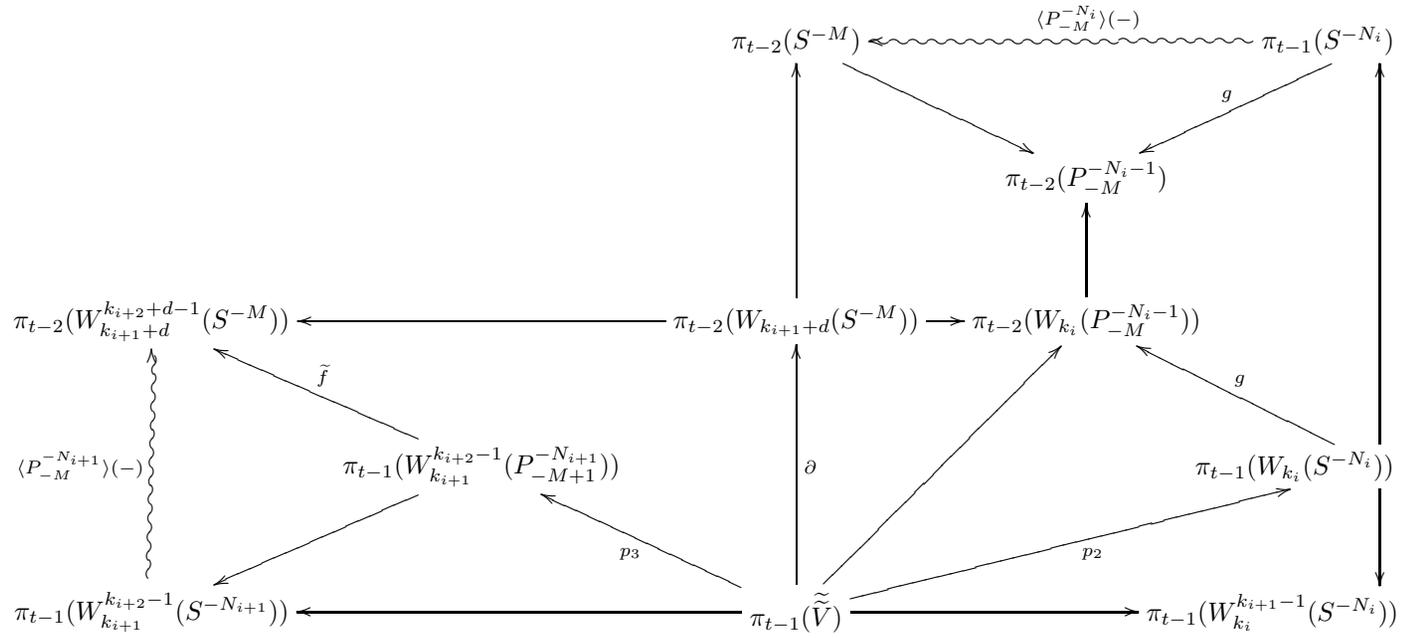
\begin{figure}
$$ \xymatrix@C-1em@R+1em{
\\
&& 
\pi_{t-2}(S^{-M}) \ar[dr] &&
*+[l]{\pi_{t-1}(S^{-N_i})} \ar@{~>}[ll]_{\bra{P_{-M}^{-N_i}}(-)} \ar[ld]_g
\\
&&&
\pi_{t-2}(P_{-M}^{-N_i-1})
\\
\pi_{t-2}(W_{k_{i+1}+d}^{k_{i+2}+d-1}(S^{-M})) &&
\pi_{t-2}(W_{k_{i+1}+d}(S^{-M})) \ar[ll] \ar[uu] \ar[r] &
\pi_{t-2}(W_{k_i}(P_{-M}^{-N_i-1})) \ar[u]
\\
&
\pi_{t-1}(W_{k_{i+1}}^{k_{i+2}-1}(P_{-M+1}^{-N_{i+1}})) \ar[ul]_{\td{f}} 
\ar[dl] &&&
*+[l]{\pi_{t-1}(W_{k_i}(S^{-N_i}))} \ar[ul]_g \ar[uuu] \ar[d]
\\
\pi_{t-1}(W_{k_{i+1}}^{k_{i+2}-1}(S^{-N_{i+1}}))
\ar@{~>}[uu]^{\bra{P_{-M}^{-N_{i+1}}}(-)} &&
\pi_{t-1}(\td{\td{V}}) \ar[ll] \ar[ul]^{p_3} \ar[uu]_{\partial} \ar[uur]
\ar[rru]_{p_2} \ar[rr] &&
*+[l]{\pi_{t-1}(W_{k_i}^{k_{i+1}-1}(S^{-N_i}))}
} $$
\caption{The main diagram for the proof of Theorem~\ref{thmC}}
\label{fig:thmC}
\end{figure}
\end{landscape}

Let $(P_{-M}^{-N_{i+1}})_{[d]}$ be the cofiber of $\td{f}$.
Our hypothesis that $d \le r_{i+1}$ implies that $\td{\td{V}}$ has the
following equivalent description as a homotopy pushout.
\begin{gather*}
P_{-M}^{-N_{i-1}} 
\sqcup_{W_{k_i}(P_{-M}^{-N_{i-1}})}
W_{k_i}(P_{-M}^{-N_i}) 
\sqcup_{W_{k_{i+1}+d}(S^{-M})}
\\
W_{k_{i+1}}((P_{-M}^{-N_{i+1}})_{[d]}) 
\sqcup_{W_{k_{i+2}}(P_{-M}^{-N_{i+1}})}
W_{k_{i+2}}(P_{-M}) 
\end{gather*}
Our proof is reduced to chasing the diagram displayed in
Figure~\ref{fig:thmC}.

In Figure~\ref{fig:thmC}, $\td{f}$ and $g$ are attaching maps, $p_2$ is
projection onto the second factor of the homotopy pushout $\td{\td{V}}$,
$p_3$ is projection onto the third factor followed by projection onto the
$-N_i$-cell, and $\partial$ is the Meyer-Vietoris boundary.
The wavy arrows represent Toda brackets.
The element $\td{\td{\gamma}} \in \td{\td{V}}$ maps to compatible elements
of every group displayed in Figure~\ref{fig:thmC}.  The image of
$\td{\td{\gamma}}$ in $\pi_{t-1}(W_{k_i}^{k_{i+1}-1}(S^{-N_i}))$ projects
to an
element of $R^{[k_i]}_E(\alpha)$.  
The image of $\td{\td{\gamma}}$ in 
$\pi_{t-1}(W_{k_{i+1}}^{k_{i+2}-1}(S^{-N_{i+1}}))$ projects to an element of
$R^{[k_{i+1}]}_E(\alpha)$.
Let $\beta$ be the image of $\td{\td{\gamma}}$ in 
$\pi_{t-2}(S^{-M})$.  Following the outside of
Figure~\ref{fig:thmC} from the
lower left hand corner clockwise reveals that
$$ \beta \in \br{\bra{P_{-M}^{-N_{i+1}}}(R^{[k_{i+1}]}_E(\alpha))}.$$
Following the outside of Figure~\ref{fig:thmC} from the lower 
right-hand corner
counterclockwise gives
$$ \beta \in \bra{P_{-M}^{-N_i}}(\br{R}^{[k_i]}_E(\alpha)). $$
Thus $\beta$ is a common element of the two groups.
\end{proof}

\begin{proof}[Proof of Proposition~\ref{prop:E^Eroot}]
Suppose that the first filtered root invariant of $\alpha$ is carried by
the $-N$-cell of $P_{-\infty}$.  If the first filtered root invariant is
trivial, then $N = \infty$.
Consider the following diagram.
$$\xymatrix{
\pi_t(tS^0) \ar@/^2pc/@{~>}[rr]^{R^{[1]}_{E}(-)} &
\pi_{t-1}(W_I(P^J)) \ar[l] \ar[dl] \ar[dd] \ar[r] &
\pi_{t-1}(W_1^1(S^{-N})) \ar@{=}[dd]
\\
\pi_{t-1}(W_1(P)) \ar[d] \ar[u]
\\
\pi_{t-1}(W_1^1(P)) \ar@/_2pc/@{~>}[rr]_{R_{E \wedge \br{E}}(-)} &
\pi_{t-1}(W_1^1(P^{-N})) \ar[l] \ar[r] &
\pi_{t-1}(W_1^1(S^{-N}))
} $$
Let $f^\alpha$ be a lift of the image of $\alpha$ in $\pi_{t-1}(tS^0)$ to 
$\pi_{t-1}(W_I(P^J))$.  Then $f^\alpha$ has compatible images in every
group in the diagram.  The detecting element $\td{\alpha}$ is the image of
$f^\alpha$ in $\pi_{t-1}(W_1^1(P_{-\infty})) = \pi_{t-1}(\Sigma^{-1}
E\wedge \br{E}^{\wedge}_p)$.
\end{proof}

\begin{proof}[Proof of Corollary~\ref{cor:E^Eroot}]
First observe that if the first filtered root invariant is carried by the
$-M$-cell, then $M \ge N$.
Let $\td{\alpha}'$ be a `preferred' detecting element such as the one
described in Proposition~\ref{prop:E^Eroot}.  It suffices to show that if
$\beta \in R_{E \wedge \br{E}}(\td{\alpha})$ and $\beta' \in R_{E \wedge
\br{E}}(\td{\alpha}')$ (or $\beta' = 0$ if 
$R_{E \wedge \br{E}}(\td{\alpha}')$ lies in a higher degree than 
$R_{E \wedge \br{E}}(\td{\alpha})$), then $\beta - \beta' \in A$.  
The elements $\beta$
and $\beta'$ are the images of elements $\gamma$ and $\gamma'$,
respectively, under the collapse map
$$ \pi_{t-1}(W_1^1(P^{-N} )) \rightarrow
\pi_{t-1}(W_1^1(S^{-N} )). $$
where $\gamma$ and $\gamma'$ both map to the image of $\alpha$ under the
homomorphism
$$ \pi_{t-1}(W_1^1(P^{-N} )) \xrightarrow{\iota} 
\pi_{t-1}(W_0^1(P )). $$
Therefore there is an element $\mu \in \pi_{t}(W_0^{(0,1)}(P_{(-\infty,
-N+1)}))$ such that $\gamma - \gamma' =
\partial_1(\mu)$, where $\partial_1$ is the connecting homomorphism of the
long exact sequence associated to the cofiber sequence
$$ W_1^1(P^{-N}) \rightarrow W_0^1(P) \rightarrow W_0^{(0,1)}(P_{(-\infty,
-N+1)}). $$
A comparison of cofiber sequences shows that there is a commutative diagram
$$ \xymatrix{
\pi_t(W_0^{(0,1)}(P_{(-\infty,-N+1)})) \ar[r]^-{\partial_1} \ar[d]_{\nu} &
\pi_{t-1}(W_1^1(P^{-N})) \ar[d]^{\nu}
\\
\pi_t(W_0^{(0,1)}(P_{(-N,-N+1)})) \ar[r]_-{\partial} &
\pi_{t-1}(W_1^1(S^{-N}))
} $$
Thus we see
that 
$$ \beta - \beta' = \nu(\gamma - \gamma') = \nu \circ \partial_1(\mu) =
\partial\circ\nu(\mu)$$ 
so $\beta - \beta'$ is contained in $A$.
\end{proof}

\begin{proof}[Proof of Theorem~\ref{thmD}]
We shall refer to the appropriate $\ext$ groups as as the $E_2$ terms of the
ASS.  Observe that Lemma~\ref{filtlem} implies that $k_1 \ge k$.
We first will prove that $k_1 = k$.  Suppose that $k_1 > k$.
Let $\td{\td{\alpha}}$ be a lift of $\td{\alpha}$ to
$\pi_t(W_k(S^0))$.
Consider the following diagram.
$$ \xymatrix{
E^{k,t+k-1}_2(S^{-1}) \ar[d]_f &
\pi_{t-1}(W_k(S^{-1})) \ar[l] \ar[d]
\\
E^{k,t+k-1}_2(P_{-\infty}) &
\pi_{t-1}(W_k(P)_{-\infty}) \ar[l] &
\pi_{t-1}(W_{k_1}(P)_{-\infty}) \ar[l]
}$$
The element $\td{\td{\alpha}}$ maps to $\td{\alpha} \in
E^{k,t+k-1}_2(S^{-1})$.  By the definition of $k_1$, the image of
$\td{\td{\alpha}}$ in $\pi_{t-1}(W_k(P)_{-\infty})$ lifts to an element of 
$\pi_{t-1}(W_{k_1}(P)_{-\infty})$.  This implies that $f(\td{\alpha}) = 0$.
However, $\td{\alpha}$ is non-zero, and the algebraic Segal conjecture 
implies that $f$ is an isomorphism.  We conclude that $k_1 = k$.

The filtered root invariant 
$R^{[1]}_{H}(\alpha)$ is a subset of $E_1^{k,t+N_1+k-1}(S^0)$.  All of
the attaching maps of $P_{-\infty}$ are of positive Adams filtration, so
Theorem~\ref{thmB} implies that $R^{[1]}_E(\alpha)$ consists of 
$d_1$-cycles.  Let $f^\alpha \in \pi_{t-1}(W_I(P^J))$ be a lift of
$\alpha$, and let $\beta \in \pi_{t-1}(W_k^k(S^{-N_1}))$ 
be the corresponding element in
$R^{[k]}_{H}(\alpha)$.  Let $f_1^\alpha \in
\pi_{t-1}(W_k(P)_{-\infty})$ be the image of $f^\alpha$ under the natural
map
$$ W_I(P^J) \rightarrow W_k(P)_{-\infty}. $$
Consulting the diagram,
$$ \xymatrix{
& 
\pi_{t-1}(W_k(P)_{-\infty}) \ar[d] 
\\
E_2^{k,t+k-1}(P_{-N_1}) &
\pi_{t-1}(W_k^{k+1}(P_{-N_1})) \ar[d] \ar[l] &
\pi_{t-1}(W_k^{k+1}(S^{-N_1})) \ar[d] \ar[l]
\\
&
\pi_{t-1}(W_k^k(P_{-N_1})) &
\pi_{t-1}(W_k^k(S^{-N_1})) \ar[l]
}$$
both $f_1^\alpha$ and $\beta$ have the same image in 
$\pi_{t-1}(W_k^k(P_{-N_1}))$.  Let 
$\gamma$ be the image of $f_1^\alpha$ in $\pi_{t-1}(W_k^{k+1}(P_{-N_1}))$.
Since $\beta$ is a $d_1$-cycle, it lifts to $\td{\beta} \in 
\pi_{t-1}(W_k^{k+1}(S^{-N_1}))$.  Let $\td{\gamma}$ be the image of
$\td{\beta}$ in $\pi_{t-1}(W_k^{k+1}(P_{-N_1}))$.  It not necessarily the case
that $\gamma = \td{\gamma}$, but it is the case that the difference $\gamma
- \td{\gamma}$ vanishes in $E_2^{k,t+k-1}(P_{-N_1})$.  

Consider the following diagram.
$$\xymatrix@C-1.5em{ 
\pi_{t}(tS^0) \ar@{~>}[rrr]^{\br{R}^{[k]}_{H}(-)} &&&
\pi_{t-1}(W_k^{k+1}(S^{-N_1})) \ar[ddd] \ar[dl] \ar `r[ddddd] `[ddddd] [ddddd]
\\
& \pi_{t-1}(W_I(P^J)) \ar[ul] \ar[ddl] \ar[r] \ar[dr] &
\pi_{t-1}(W_k^k(S^{-N_1})) \ar[d]
\\
&& \pi_{t-1}(W_k^k(P_{-N_1})) 
\\
\pi_{t-1}(W_k(P)_{-\infty}) \ar[uuu] \ar[rrr] \ar[d] &&&
\pi_{t-1}(W_k^{k+1}(P_{-N_1})) \ar[ul] \ar[d] 
\\
E_2^{k,t+k-1}(P_{-\infty}) \ar[rrr]_{\nu_{N_1}} &&&
E_2^{k,t+k-1}(P_{-N_1})
\\
E_2^{k,t+k-1}(S^{-1}) \ar[u]^{\cong}_{f} \ar@{~>}[rrr]_{R_{alg}(-)} &&&
E_2^{k,t+k-1}(S^{-N_1}) \ar[u]
}$$
The element $\td{\beta}$ maps to an element of
$\br{R}^{[k]}_{H}(\alpha)$ in $E_2^{k,t+k-1}(S^{-N_1})$.  The map $f$
is an isomorphism, so the image of $f_1^\alpha$ lifts to an element
$\td{\alpha} \in E_2^{k,t+k-1}(S^{-1})$.  This is the choice of $\td{\alpha}$
which detects $\alpha$ that we appeal to in the statement of
Theorem~\ref{thmB}.  Recall that the image of $f_1^\alpha$ in 
$\pi_{t-1}(W_k^{k+1}(P_{-N_1}))$
was $\gamma$ while the image of $\beta$ is $\td{\gamma}$.  Since the
difference $\gamma - \td{\gamma}$ maps to zero in
$E_2^{k,t+k-1}(P_{-N_1})$, we may conclude that the image of $\td{\beta}$
in $E_2^{k,t+k-1}(S^{-N_1})$ is also an element of the algebraic root
invariant,
\emph{provided the algebraic
root invariant does not lie in a higher stem.}

We have proven that either what is claimed in Theorem~\ref{thmD} holds, 
or the 
algebraic root invariant lives in a larger stem than 
$R^{[1]}_{H}(\alpha)$.
We shall now show that this cannot happen.  Let $M$ be maximal such that 
the image of $\td{\alpha}$ in $E_2^{k, t+k-1}(S^{-1})$ maps to zero under
the composite
$$ E_2^{k, t+k-1}(S^{-1}) \xrightarrow{\nu_{M-1} \circ f} 
E_2^{k, t+k-1} (P_{-M+1}). $$
By the algebraic Segal conjecture, such a finite $M$ exists.
We wish to show that $N_1 = M$.  So far we know that $N_1 \le M$.  
In light of the definition of $N_1$, we 
simply must show that the $\alpha$ is sent to zero under the
composition
$$ \pi_{t-1}(S^{-1}) \rightarrow \pi_{t-1}
(W_0^{(k-1,k)}(P_{(-M+1,\infty)})).
$$
To this end consider the following diagram.
$$\xymatrix{
& \pi_t(W_0^{k-1}(P)_{-\infty}) \ar[d]_\partial &
\pi_t(W_{k-1}^{k-1}(P)_{-\infty}) \ar[d]^\nu \ar[l]
\\ 
\pi_{t-1}(W_k(P)_{-\infty}) \ar[r] \ar[d] &
\pi_{t-1}(W_k^k(P_{-M+1})) \ar[d] &
\pi_t(W_{k-1}^{k-1}(P_{-M+1})) \ar[l]^{d_1}
\\
\pi_t(tS^0) \ar[r] &
\pi_{t-1}(W_0^{(k-1,k)}(P_{(-M+1,\infty)}))
}$$
The central 
vertical column corresponds to part of the long exact sequence for a
cofibration.
The element $f^\alpha_1 \in \pi_{t-1}(W_k(P)_{-\infty})$ maps to the image
of $\alpha$ in $\pi_t(tS^0)$.  Let $g$ be the image of $f^\alpha_1$ in
$\pi_{t-1}(W_k^k(P_{-M+1}))$. By our choice of $M$, $g$ must vanish in
$E_2$, that is to say, there must be an element $h \in
\pi_{t}(W_{k-1}^{k-1}(P_{-M+1}))$ such that $d_1(h) = g$.  
The map $\nu$ is surjective, because we are
dealing with ordinary homology, and the map
$$ H_*(P)_{-\infty} \rightarrow H_*(P_{N})$$
is surjective for every $N$.  Thus $h$ lifts to an element $\td{h} \in 
\pi_t(W_{k-1}^{k-1}(P)_{-\infty})$.  Applying $\partial$ to the image of
$\td{h}$ in $\pi_t(W_0^{k-1}(P)_{-\infty})$, we get $g$, so $g$ is in the
image of $\partial$.  Thus, by exactness of the vertical column, the image
of $g$ in $\pi_{t-1}(W_0^{(k-1,k)}(P_{(-M+1,\infty)}))$ is trivial.  Thus
the image of $\alpha$ in $\pi_t(tS^0)$ maps to zero in 
$\pi_{t-1}(W_0^{(k-1,k)}(P_{(-M+1,\infty)}))$, which is what we were trying
to prove.
\end{proof}


\section{$bo$ resolutions}\label{sec:bo}

The $bo$ resolution for the sphere was the motivating example for the
theorems in Section~\ref{sec:results}.
The $d_1$'s on the $v_1$-periodic summand of the $0$ and $1$-lines of the 
$bo$-resolution reflect the root invariants of the elements $2^i$ at the
prime $2$.  We do not rederive these root invariants, but satisfy ourselves
in explaining how the theorems of Section~\ref{sec:results} 
play out in this context.  
Alternatively,
one might interpret this section as explaining how to use $bo$ to compute
the root invariants in the $K(1)$-local category.  This section is
furthermore meant to
explain how the techniques of this paper, when applied to the small
resolution of $L_{K(2)}(S^0)$ described in \cite{GoerssHennMahowaldRezk},
could be used to compute root invariants in the $K(2)$-local category. 

Recall that $\Sigma^4 bsp$ is a summand of $bo \wedge \Sigma \br{bo}$, 
hence its 
homotopy is a stable summand of the $1$-line of the $E_1$-term of the 
$bo$-ASS.  The composite of the $d_1$ of this spectral sequence
with the projection onto the summand
$$ \pi_t(bo) \xrightarrow{d_1} \pi_{t-1} (bo \wedge \br{bo}) \rightarrow
\pi_{t-1}(\Sigma^4 bsp) $$
is the map $\psi^3 - 1$ (up to a unit in $\ZZ_{(2)}$).
The fiber of $\psi^3 - 1$ is the $2$-primary $J$ spectrum.
\begin{figure}
\begin{center}
  \includegraphics[width=4.9in]{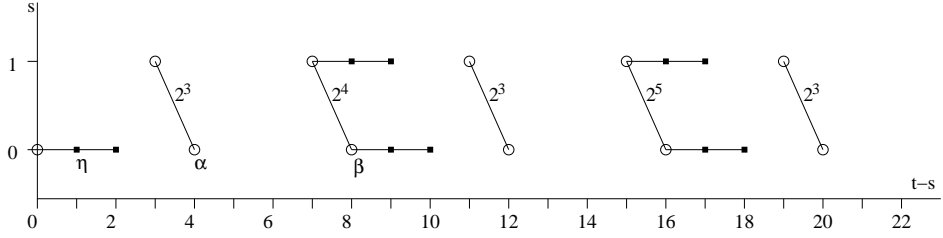}
\end{center}
\caption{The $v_1$-periodic summand of the $bo$-ASS.}
\label{bofig}
\end{figure}
Figure~\ref{bofig} shows the $v_1$-periodic summand of the 
$bo$-ASS 
with differentials.  In this figure, dots represent copies of $\FF_2$ and 
circles represent copies of $\ZZ$.  The $d_1$ differentials are (up to a
$2$-adic unit) multiplication
by the power of $2$ indicated.  What survives are the $v_1$-periodic elements
in $\pi_*^S$ at $p=2$.

Our theorems explain the relationship between $R_{bo} (2^k)$ and $R(2^k)$ for
all $k > 0$.  It is quite straightforward to calculate $R_{bo}(2^k)$
(see \cite{MahowaldRavenel}).  Let $\eta$, $\alpha$, and $\beta$ be the
multiplicative generators of $\pi_*(bo)$ in dimensions $1$, $4$, and $8$,
respectively.  The we have the following $bo$-root invariants.
$$
R_{bo}(2^k) = \begin{cases}
\beta^i & k = 4i \\
\eta \beta^i &  k = 4i+1 \\
\eta^2 \beta^i & k = 4i+2 \\
\alpha \beta^i & k = 4i+3
\end{cases}
$$
These are the filtered root invariants $R_{bo}^{[0]}(2^k)$, by
Proposition~\ref{prop:Eroot}.    
For $k \equiv 1,2 \pmod 4$, the elements $R_{bo}^{[0]}(2^k)$ are 
in the Hurewicz image
of $bo$, hence they are permanent cycles.  They detect
elements in $R(2^k)$.  

For $ k \equiv 3,4 \pmod 4$, the elements $R_{bo}(2^k)$ are not in the 
$bo$-Hurewicz
image.  Therefore, the elements $R^{[0]}_{bo} (2^k)$ support
differentials.  We have
$$ d_1 R^{[0]}_{bo} (2^k) = 2 \cdot \td{\alpha}_k $$ 
where $\td{\alpha}_k$ 
survives to a $v_1$-periodic element of order $2$ in dimension $4k-1$.
In the $bo$-bifiltration of $2^k$, the $-N_1$ cell carries the zeroth
filtered root invariant, where $N_1$ is given by
$$ N_1 = \begin{cases}
8i+5 & k = 4i+3 \\
8i+1 & k = 4i.
\end{cases}
$$
The first cell to attach to the $-N_1$ cell in
$P_{-N_1}$ is the $-N_1+1$ cell, and the attaching map is the degree $2$
map.
Therefore, we may deduce from Theorem~\ref{thmB} that $\td{\alpha}_k$
is an element of $R^{[1]}_{bo}(2^k)$.  The element $\td{\alpha}_k$ 
detects the
homotopy root invariant $R(2^k)$.  


\section{The algebraic Atiyah-Hirzebruch spectral sequence}\label{sec:aahss}

As discussed after the statement of Corollary~\ref{corA}, it is often the
case that one may know $R^{[k]}_E (\alpha)$ contains a permanent cycle in
the $E$-Adams spectral sequence, but nothing else.  One would want
to conclude that this permanent cycle detects the homotopy root invariant
$R(\alpha)$, but Theorem~\ref{thmA} says this
is only true modulo obstructions in higher $E$-Adams filtration.  
However, it is sometimes the case that brute
force computations in the $E_2$-term of the $E$-Adams spectral sequence 
for $P_{N}$ will 
yield enough data to eliminate such obstructions.  This section is
concerned with computations of the $E_2$-term $E_2(P_{N})$.
Specifically, 
we shall assume that $E = BP$
and describe the algebraic Atiyah-Hirzebruch spectral sequence (AAHSS), 
which computes $\ext(BP_*P_{N})$ from $\ext(BP_*)$.  In
Section~\ref{sec:proc}, we describe
a procedure (Procedure~\ref{procA})
that explains how a limited knowledge of differentials in the AAHSS
can be exploited to compute homotopy root invariants from filtered root
invariants.

The AAHSS is the spectral sequence obtained by applying $\ext$ to the
skeletal filtration of a complex.
The AAHSS has appeared in various forms in the literature.  In
\cite{Mahowald}, Mahowald uses ordinary cohomology, and uses this spectral
sequence to compute $\ext$ of various stunted projective spaces.  In
\cite{Sadofsky}, Sadofsky uses the $BP$ version that is the subject of this 
section.  There is some overlap with this section and the computations in
\cite{Sadofsky}.

We will first describe the AAHSS for 
computing the Adams-Novikov $E_2$-term $E_2(P_{N})$, 
where $N \equiv -1 \pmod q$.  
Let $F$ be the universal $p$-typical formal group law associated 
to $BP$, and write the $p$-series of $F$ as
$$ [p]_F(x) = \sum^F_{i \ge 0} v_i x^{p^i} = \sum_j c_j x^{(p-1)j + 1} $$
where $v_i$ are the Araki generators.  Then we can say the following about 
the coefficients $c_i$.
\begin{align*}
c_0 & = p \\
c_1 & = v_1 \\
c_j & \equiv 0 \pmod p \text{ for $1 < j < p + 1$} \\
c_{p+1} & \equiv v_2 \pmod p
\end{align*}
There are short exact sequences (compare with \cite[2.3]{Sadofsky})
\begin{gather}\label{eq:BPSES1}
0 \rightarrow BP_*(S^{kq-1}) \xrightarrow{\phi} BP_*(P^{kq-1}_N)
\rightarrow BP_*(P^{kq}_N) \rightarrow 0 \\
0 \rightarrow BP_*(P^{(k-1)q}_N) \rightarrow BP_*(P^{kq-1}_N) 
\xrightarrow{\nu} BP_*(S^{kq-1}) \rightarrow 0
\end{gather}
and these give rise to long exact sequences upon applying $\ext$.  They
have the following boundary homomorphisms.
\begin{align*}
\delta_1 & : \ext^{s,t}(BP_*P_N^{kq}) \rightarrow 
\ext^{s+1,t}(BP_*S^{kq-1}) \\
\delta_2 & : \ext^{s,t}(BP_*S^{kq-1}) \rightarrow 
\ext^{s+1,t}(BP_*P^{(k-1)q}_N)
\end{align*}
For $l \le k$, $BP_*(P_{lq - 1}^{kq})$ is generated as a $BP_*$-module by 
elements
$e_{jq-1}$ in dimension $jq-1$.
The map $\phi$ in the short exact sequence~\ref{eq:BPSES1} is given by
$$ \phi (\iota_{kq-1}) \mapsto \sum_i c_i e_{q(k - i) - 1}. $$
Splicing these sequences together creates an exact couple,
and the resulting spectral sequence is the AAHSS.
We shall index it just as one indexes the stable EHP spectral sequence
\cite[1.5]{Ravenel}.  In fact, if $N = q-1$ then this is precisely  
the $BP$-algebraic stable EHP spectral sequence.
Thus we have a spectral sequence
$$ E_1^{k,n,s} \Rightarrow \ext^{s, s+k} (BP_*P_N) $$
where the $E_1$ term is described below.
\begin{align*}
E_1^{k, 2m, s} & = \ext^{s, s+k} (BP_*S^{mq - 1}) \\
E_1^{k, 2m+1, s} & = \ext^{s+1, s+k} (BP_*S^{mq - 1}) 
\end{align*}
The indexing works out so that
$$ d_r: E_r^{k,n,s} \rightarrow E_r^{k-1, n-r, s+1}. $$

If we wish to compute $\ext(BP_*P_N)$ for $N = lq$, simply 
truncate the AAHSS for $\ext(BP_*P_{N-1})$
by setting $E_1^{k, 2l, s} = 0$ for all $k$ and $s$.

We shall refer to 
an element in the $E_1$-term of the AAHSS 
by its name in $\ext(BP_*)$ and the cell that is borne on.
Thus if $\gamma \in \ext(BP_*)$ is in 
$E_1^{k, 2m, s}$, then we shall refer to it as
$\gamma [mq - 1]$.  Likewise, if $\gamma$ is in $E_1^{k, 2m+1, s}$,
we shall refer to it as $\gamma[mq]$.

In order to implement Procedure~\ref{procA} we shall 
need to know how to explicitly
compute differentials in the filtered spectral sequence.
It is useful to use the diagram below. \emph{In this diagram and the
remainder of this section we will drop $BP_*$ from
our $\ext$ notation for compactness.}
$$
\xymatrix@R-1.7em{
&  \ext^{s+1,s+k}(P_N^{lq-1}) \ar[r]_-{\nu_*} \ar[dd] &
\ext^{s+1,s+k}(S^{lq-1})  
\\
& & E_1^{k-1, 2l, s+1} \ar@{=}[u]
\\
&  \ext^{s+1, s+k}(P_N^{lq}) 
\ar[r]_-{\delta_1} \ar[dd] & \ext^{s+2,
s+k}(S^{lq-1})  
\\
& & E_1^{k-1,2l+1,s+1} \ar@{=}[u]
\\
&  \vdots \ar[dd] 
\\
E_1^{k, 2m, s} \ar@{=}[d]
\\
\ext^{s, s+k}(S^{mq-1}) 
\ar[r]_-{\delta_2} &
\ext^{s+1,s+k}(P_N^{(m-1)q}) \ar[dd] 
\\
E_1^{k, 2m+1, s} \ar@{=}[d]
\\
\ext^{s+1,s+k}(S^{mq-1}) \ar[r]_-{\phi_*} &
\ext^{s+1, s+k}(P_N^{mq-1})
} $$
Suppose $\gamma[n]$ is an element of  $E_1$, where 
$n = mq - \epsilon$, $\epsilon =
0,1$.  Let $\gamma_1$ be the image of $\gamma$ in $\ext^{s+1,
s+k}(P_N^{mq-1})$ or $\ext^{s+1,
s+k}(P_N^{(m-1)q})$, depending on the value of $\epsilon$.  
Then lift $\gamma_1$ as
far as possible up the tower in the center of the above diagram.  
Suppose that $\gamma_1$ lifts to $\gamma_2 \in
\ext^{s+1,s+k}(P_N^{lq-\epsilon_1})$, where $\epsilon_1 = 0,1$.  Then 
if $\gamma_3$ is the image of $\gamma_2$ under $\nu_*$ or $\delta_1$
(depending on the value of $\epsilon_1$), there is an AAHSS differential
$$ d_r(\gamma[n]) = \gamma_3[n'] $$
where $n'=lq-\epsilon$.

The above description makes the following proposition clear. 
In the statement of this proposition, and in what follows, when 
we say that differentials are 
computed modulo $I_n$, we mean is that we are considering their
images in $ \ext(BP_*(-)/I_n)$.

\begin{prop}\label{prop:d_r, r odd}
Modulo $I_n$, we have the differential 
$$ d_{2p^n - 1} (\alpha[mq]) = v_n \cdot \alpha [mq-2p^n+1]. $$
\end{prop}

\begin{proof}
Take $\alpha \in \ext^{s+1, s+k}(S^{mq-1})$.  
Clearly, it is the case that
$$ \phi_*(\alpha) \underset{\pmod {I_n}}{\equiv} v_n \alpha
[mq-2p^n+1] + \cdots $$
lifts to $\ext^{s+1, s+k}(P^{mq-2p^n+1})$, and the image of this
lift under the map $\nu_*$ is 
$$ v_n \cdot \alpha \in \ext^{s+1,s+k}(S^{mq-2p^n+1}). $$
\end{proof}

In order to compute $\delta_1$ and $\delta_2$, we must know something of
the $BP_*BP$ comodule structure of $BP_*(P_N)$.  It suffices to
understand the $BP_*BP$ comodule structure of $BP_*(B\ZZ/p)$, 
since $P_0$ is a stable summand of $B\ZZ/p_+$.  The $BP_*BP$ comodule
structure of $BP_*(P_N)$ may then be deduced from James periodicity.
In \cite{Sadofsky}, it is proven that there is a cofiber sequence
$$ \CC P^\infty \rightarrow (\CC P^\infty)^{\xi^p} \rightarrow
\Sigma B\ZZ/p
$$
(where $\xi$ is the canonical complex line bundle on $\CC P^\infty$)
which induces the short exact sequence below.
$$ 0 \rightarrow BP_*\CC P^\infty \rightarrow BP_*(\CC P^\infty)^{\xi^p} 
\rightarrow BP_* \Sigma B\ZZ/p \rightarrow 0
$$
Since $BP_* (\CC P^\infty)^{\xi^p}$ surjects onto 
$BP_* \Sigma B\ZZ/p$, it suffices to understand the
$BP_*BP$-comodule structure of the former.

Recall that $BP^*(\CC P^\infty) = BP^* [[x]]$ where $x$ is the Euler
class of $\xi$.  Therefore, the image of the inclusion
$$ BP^*(\CC P^\infty)^{\xi^p} \rightarrow BP^*(\CC P^\infty) $$
is the ideal
$$ [p]_F(x) \cdot BP^* [[x]]  \subseteq BP^*(\CC P^\infty). $$
We shall identify $BP^*(\CC P^\infty)^{\xi^p}$ with this ideal.

Let $y_{2k}$ be the generator in $BP_{2k}((\CC P^\infty)^{\xi^p})$ dual
to $[p]_F(x) \cdot x^{k-1}$.  
Let $f(x)$ be the formal power series over $BP_*BP$ whose inverse is given by
$$ f^{-1}(x) = \sum^F_{i \ge o} t_i x^{p^i}.  $$
The power series $f(x)$ is the universal strict isomorphism of 
$p$-typical formal group
laws.

\begin{prop}
The $BP_* BP$ coaction on $y_{2k}$ is given below.
\begin{equation}\label{eq:coaction}
 \psi(y_{2k}) = \sum_{i+j = k} \left(\frac{f([p]_F(x))}{[p]_F(x)}
(f(x))^{j-1} \right)_{k-1} \otimes y_{2j}
\end{equation}
Here the subscript $k-1$ indicates the coefficient of $x^{k-1}$.
\end{prop}

\begin{proof}
We digress for a moment on cooperations.
For a spectrum $X$, the left unit 
$\eta_L : BP \rightarrow BP \wedge BP$ induces the $BP_*BP$ coaction
$$ \psi : BP_*(X) \xrightarrow{(\eta_L)_*} (BP \wedge BP)_*(X) \cong 
BP_*BP \otimes_{BP_*} BP_*(X). $$
It is sometimes convenient to consider the dual coaction on cohomology
$$ \psi^* : BP^*(X) \xrightarrow{(\eta_L)_*} (BP \wedge BP)^*(X). $$

In Theorem~11.3 of Part~II of \cite{Adams}, Adams describes the $MU_*MU$ 
coaction on $MU_*(\CC P^\infty)$.  When translated into $BP$-theory, the
coaction formula reads
$$ \psi(e_{2k}) = \sum_{i+j = k} \left( (f(x))^j \right)_k 
\otimes e_{2j}.
$$
Here, $e_{2k}$ is the generator of $BP_{2k}(\CC P^\infty)$ dual to $x^k \in
BP^{2k}(\CC P^\infty)$.  Upon dualization, we get the following formula for
the dual coaction on an element $h(x)$ in $BP^*(\CC P^\infty)$.
$$ \psi^*(h(x)) = f_*h(f(x)) $$
The polynomial $f_*h(x)$ is the polynomial obtained from applying the
right unit to the coefficients of $h(x)$.
We deduce the dual coaction on an element
$[p]_F(x)\cdot h(x)$ in $BP^*((\CC P^\infty)^{\xi^p})$, regarded as the
ideal $([p]_F(x))$ contained in $BP^*[[x]]$.
\begin{align*}
\psi^* ([p]_F(x)\cdot h(x)) & = f_*[p]_F(f(x)) \cdot f_*h(f(x)) \\
& = [p]_{f_*F}(f(x))\cdot f_*h(f(x)) \\
& = f([p]_F(x)) \cdot f_*h(f(x)) \\
& = [p]_F(x) \cdot \frac{f([p]_F(x))}{[p]_F(x)} \cdot f_*h(f(x))
\end{align*}
Here, $f_*F$ is the pushforward of the formal group law $F$ under the map
$f$.  
Letting $h(x) = x^{j-1}$ and dualizing, we have the desired formula for the
coaction on $BP^*((\CC P^\infty)^{\xi^p})$.
\end{proof}

We wish to use the above coaction formula to compute differentials in
the AAHSS with $p \ge 3$.
The first few terms of the relevant power series are listed below.
Everything in what follows is written in terms of Hazewinkel generators.
\begin{gather}
[p]_F(x)  = px + (1 - p^{p-1}) v_1x^p + \mathcal{O}(x^p+1) 
\label{pseries} \\
f(x)  = x - t_1 x^p + (pt_1^2 + v_1 t_1)x^{2p-1} + \mathcal{O}(x^{2p})
\label{fseries} 
\end{gather}
\begin{multline}
\label{qseries}
\frac{f([p]_F(x))}{[p]_F(x)}  = 1 - p^{p-1} t_1 x^{p-1} + \\
p^{p-2}(p^{p+1} t_1^2 + (1-p-p^{p-1}+2p)v_1 t_1) x^{2p-2} +
\mathcal{O}(x^{2p-1})
\end{multline}
We may compute these series further if we work modulo the ideal $I_2 =
(p, v_1)$.
\begin{gather}
[p]_F(x) \equiv v_2 x^{p^2} + \mathcal{O}(x^{p^2}+1) \pmod {I_2}
\label{pseries_modI2} \\
f(x) \equiv x - t_1 x^p + \mathcal{O}(x^{p^2}) \pmod {I_2}
\label{fseries_modI2} \\
\frac{f([p]_F(x))}{[p]_F(x)} \equiv 1 + \mathcal{O}(x^{p^2-1}) \pmod {I_2}
\label{qseries_modI2}
\end{gather}
These formulas give rise to the following proposition.
\begin{prop}\label{prop:odd_diffs}
In the filtered spectral sequence, the following formulas for
differentials hold up to multiplication by a unit in $\ZZ_{(p)}$.
\begin{align*}
  & d_2 (\alpha[kq-1]) \doteq
  \begin{cases}
    \alpha_1 \cdot \alpha[(k-1)q-1] & k \not \equiv 0 \pmod p \\
    0 & k \equiv 0 \pmod p
  \end{cases}
  \\
  & d_2 (\alpha[kq]) \underset{\pmod{I_1}}{\dotequiv}
  \begin{cases}
    \alpha_1 \cdot \alpha[(k-1)q] & k \not \equiv 1 \pmod p \\
    0 & k \equiv 1 \pmod p
  \end{cases}
  \\
  & d_q (\alpha[kq-1])  \underset{\pmod {I_2}}{\dotin}
  \begin{cases}
    \langle \alpha, \underbrace{\alpha_1, \ldots, \alpha_1}_{p-1}
    \rangle[(k-p+1)q-1] & k \equiv -1 \pmod p \\
    0 & k \not \equiv -1 \pmod p
  \end{cases}
  \\
  & d_q (\alpha[kq])  \underset{\pmod {I_2}}{\dotin}
  \begin{cases}
    \langle \alpha, \underbrace{\alpha_1, \ldots, \alpha_1}_{p-1}
    \rangle[(k-p+1)q] & k \equiv 0 \pmod p \\
    0 & k \not \equiv 0 \pmod p
  \end{cases}
\end{align*}
\end{prop}

\begin{proof}
We shall begin by computing $d_2(\alpha[kq-1])$.  In light of the
coaction equation \ref{eq:coaction}, we compute from \ref{fseries} and
\ref{qseries}
$$
\frac{f([p]_F(x))}{[p]_F(x)}\left( f(x) \right)^r  = 
x^r - (r+p^{p-1})t_1 x^{r+p-1} + \mathcal{O}(x^{r+p}).$$
We may conclude that the coaction formula on the generator $e_{kq-1}$ of
$BP_*(P^\infty_0)$ is given by
\begin{equation}\label{eq:t1coaction}
\psi(e_{kq-1}) = 1 \otimes e_{kq-1} + ((k-1)(p-1)-1) t_1 \otimes
e_{(k-1)q-1} + \cdots 
\end{equation}
Now $t_1$ represents $\alpha_1$ in the cobar complex, and $p \cdot
\alpha_1 = 0$.  

We now compute $d_2$.  
$$\delta_2 (\alpha[kq-1]) = ((k-1)(p-1)-1) \alpha | t_1 [(k-1)q-1] +
\cdots $$
and it follows immediately that 
\begin{align*}
d_2(\alpha[kq-1]) & = ((k-1)(p-1)-1) \alpha_1 \cdot \alpha [(k-1)q-1] \\
 & \equiv -k \alpha_1 \cdot \alpha[(k-1)q-1] \pmod p
\end{align*}
Of course, since $p \cdot \alpha_1 = 0$, one needs only to compute this
differential modulo $p$.

Now we shall deal with $d_2(\alpha[kq])$.  We shall work modulo $I_1$.
$$ \phi_* \alpha[kq] \equiv v_1 \cdot \alpha v_1 [(k-1)q-1] + \cdots $$
Therefore $\phi_* \alpha[kq]$ lifts to $\ext(P^{(k-1)q}_0)$.
So 
$$d_2(\alpha[kq]) = \delta_1 (v_1 \cdot \alpha[(k-1)q-1] + \cdots).$$
The boundary homomorphism $\delta_1$ is computed below on the level of
cobar complexes.
$$
\xymatrix@C-1em{
C^*(P^{(k-1)q-1}_0) \ar[r]^{\iota_*} \ar[dr]^d & C^*(P^{(k-1)q}) 
\\
C^*(S^{(k-1)q-1}) \ar[r]^{\phi_*} & C^*(P^{(k-1)q-1}_0)
\\
v_1 \cdot \alpha [(k-1)q-1] + \cdots \ar@{|->}[r] \ar@{|->}[dr]  &
v_1 \cdot \alpha [(k-1)q-1] + \cdots
\\
-(k-1) \alpha | t_1 [(k-1)q-1] \ar@{|->}[r] & -(k-1) v_1 \cdot
\alpha | t_1 [(k-2)q-1] + \cdots
} $$
Here, $C^*(-)$ represents the cobar complex.  We conclude that 
$$ \delta_1 (v_1 \cdot \alpha [(k-1)q-1] + \cdots) \equiv (1-k) \alpha
| t_1 [(k-1)q-1]  \pmod {I_1} $$

The $d_q$'s are proved similarly.  The relevant information to be
gathered from Equations \ref{fseries_modI2} and \ref{qseries_modI2} is
\begin{multline*}
\frac{f([p]_F(x))}{[p]_F(x)} \left( f(x) \right)^r \underset{
\pmod {I_2}}{\equiv}
x^r + \binom{r}{1} (-t_1) x^{r + (p-1)}  + \cdots 
\\
+ \binom{r}{p-1}
(-t_1)^{p-1} x^{r + (p-1)^2} + \mathcal{O}(x^{r + (p-1)^2 + 1})
\end{multline*}
The Massey product 
$$\langle \alpha, 
\underbrace{\alpha_1, \ldots, \alpha_1}_{p-1}
\rangle$$
corresponds (up to multiplication by a unit) to 
$$\alpha | t_1^{p-1} + \cdots $$
in the cobar complex.
We have that
$$ d_q(\alpha[kq-1]) \underset{\pmod {I_2}}{\in}
 \binom{r}{p-1} \langle \alpha, \underbrace{\alpha_1, \ldots,
\alpha_1}_{p-1}
\rangle [(k-p+1)q-1]. $$
where $r = (k-p+1)(p-1)-1$.
Consider the coefficient
$$ \binom{r}{p-1} = \frac{r \cdot (r-1) \cdots (r-p+2)}{(p-1)!}. $$
Since $p \not | \, p-1$, in order for $p$ to not divide this binomial
coefficient, it must be the case that $r \equiv -1 \pmod p$.  In terms
of $k$, this translates to $k \equiv -1 \pmod p$.

Just like in the calculation of $d_2\alpha[kq]$, the map $\delta_1$
introduces an offset in coefficients which results in the claimed
formula for $d_q \alpha[kq]$.
\end{proof}
 
It is interesting to note that these differentials, combined with the
image of $J$ differentials,  are enough to compute
the AAHSS completely through a certain range at odd primes.

In $\ext(BP_*)$, many elements live in an $\alpha_1$-$\beta_1$ tower, that
is, a copy of $P[\beta_1] \otimes E(\alpha_1)$.
The following general observation says that elements of the $E_1$-term of
the AAHSS which lie within an $\alpha_1$-$\beta_1$ tower typically 
either support
differentials are the target of differentials.

\begin{prop}\label{prop:beta_tower}
Let $\gamma \in \ext(BP_*)$.
Suppose furthermore that there exist elements $\lambda$ and $\mu$ in the
$E_2$ term of the Adams-Novikov spectral sequence such that either
\begin{gather*}
\alpha_1 \lambda = \gamma \\
\langle \gamma, \underbrace{\alpha_1, \ldots,  \alpha_1}_{p-1} \rangle = 
\mu
\end{gather*}
or
\begin{gather*}
\langle \lambda, \underbrace{\alpha_1, \ldots,  \alpha_1}_{p-1} \rangle
= \gamma \\
\alpha_1 \gamma = \mu.
\end{gather*}
(Here we insist that the Massey products have no indeterminacy.)
Assume furthermore that $p$ and $v_1$ do not divide any of $\lambda$,
$\mu$, or $\gamma$.
Then, in the AAHSS,  $\gamma[n]$ is either killed by a
differential, or supports a non-trivial differential.  
\end{prop}

\begin{proof}
The condition on $p$ and $v_1$ not dividing any of the elements 
$\lambda$, $\mu$, or $\gamma$ is just to ensure that we may use the
formulas in Proposition~\ref{prop:odd_diffs}.
Let us suppose that $\gamma$ satisfies the first set of conditions.
\begin{description}
\item[\rm{Case 1: $n \equiv -1 \pmod q$}] 
Write $n = mq-1$.  
  \begin{description}
    \item[\rm{Subcase (a): $m \equiv -1 \pmod p$}]
    Then $d_q \gamma [mq-1] = \mu[(m-p+1)q-1]$ 
    or $d_i \gamma[mq-1] \ne 0$ for some $i < q$.  

    \item[\rm{Subcase (b): $m \not \equiv -1 \pmod p$}]
    Then $d_2 \lambda[(m+1)q-1] = \gamma[mq-1]$.
  \end{description}

\item[\rm{Case 2: $n \equiv 0 \pmod q$}]
Write $n = mq$.
  \begin{description}
    \item[\rm{Subcase (a): $m \equiv 0 \pmod p$}]
    Then either $d_q \gamma [mq] = \mu[(m-p+1)q]$
    or we have $d_i \gamma[mq] \ne 0$ for some $i < q$.

    \item[\rm{Subcase (b): $m \not \equiv 0 \pmod p$}]
    Then $d_2 \lambda[(m+1)q] = \gamma[mq]$.
  \end{description}
\end{description}
Suppose now that $\gamma$ satisfies the second set of conditions.
\begin{description}
\item[\rm{Case 1: $n \equiv -1 \pmod q$}]
Write $n = mq-1$.
  \begin{description}
    \item[\rm{Subcase (a): $m \equiv 0 \pmod p$}]
    Then $d_q \lambda[(m+p-1)q-1] = \gamma[mq-1]$.

    \item[\rm{Subcase (b): $m \not \equiv 0 \pmod p$}]
    Then $d_2 \gamma[mq-1] = \mu[(m-1)q-1]$.
  \end{description}

\item[\rm{Case 2: $n \equiv 0 \pmod q$}]
Write $n = mq$.
  \begin{description}
    \item[\rm{Subcase (a): $m \equiv 1 \pmod p$}]
    Then $d_q \lambda[(m+p-1)q] = \gamma[mq]$.

    \item[\rm{Subcase (b): $m \not \equiv 1 \pmod p$}]
    Then $d_2 \gamma[mq] = \mu[(m-1)q]$
    or $d_1 \gamma[mq] \ne 0$. (This is the case of $p \cdot \gamma \ne
    0$.)
  \end{description}
\end{description}
\end{proof}

It is also useful to record all of the differentials supported by
the $v_1$-periodic elements of $\ext(BP_*)$ for $p > 2$ in the AAHSS.
These results are used in executing Step~$4$ of Procedure~\ref{procA}.

\begin{prop}\label{prop:imJdiffs_odd}
In the AAHSS, the differentials supported by
$Im J$ elements whose targets are $Im J$ elements are given by
\begin{align*}
d_1 (1[kq]) & = p[kq-1] \\
d_1 (\alpha_{i/j} [kq]) & = \alpha_{i/(j-1)} [kq-1] \text{ for $j >1$} \\
d_{2j}(1[kq-1]) & \doteq \tilde\alpha_{j}[(k-j)q-1] \text{ where 
$\nu_p(k) = j-1$.} \\ 
d_{2j+1} (\alpha_i[kq]) & \doteq \tilde\alpha_{i+j} [(k-j)q-1] \text{ where
$\nu_p(k+i) = j-1$.}
\end{align*}
(Here $\tilde\alpha_{l}$ is $\alpha_{l/m}$ for 
$m = \nu_p(l) + 1$.  It is the additive generator of $im J$ in the $lq-1$
stem.)
\end{prop}

\begin{proof}
These differentials detect the corresponding differentials in the
AHSS which converge to the stable
homotopy of $P_{N}$.  These differentials were computed for $P_0$ by
Thompson in
\cite{Thompson} and are summarized in \cite{Ravenel}.
\end{proof}

There are still $v_1$-periodic elements of $\ext(BP_*)$ 
which are eligible to support
differentials in the filtered spectral sequence (or subsequently in the
Adams-Novikov spectral sequence for $P_{N}$) 
whose targets are elements which are not $v_1$-periodic. 
These elements are described below.

\begin{cor}\label{cor:imJsurvivors}
The only elements in  $imJ$ in the AAHSS for
$P_{-\infty}$ which are neither targets or sources of the
differentials described in Proposition~\ref{prop:imJdiffs_odd} 
are
$$ \tilde\alpha_i[kq-1] $$
where $p^i | (k+i)$. 
\end{cor}

We obtain the following consequence which will become quite relevant
when we execute Step~$4$ of Procedure~\ref{procA}.

\begin{prop}\label{prop:beta_survive}
None of the elements $\beta_{i/j}[kq]$ are the target of differentials 
in the AAHSS or subsequently in the Adams-Novikov spectral sequence for
$\pi_*(P_N)$. 
\end{prop}

\begin{proof}
The only elements which could kill $\beta_{i/j}[kq]$ in the filtered
spectral sequence are of the form $\alpha_n [mq]$ and $1[mq-1]$.  
But these cannot by
Corollary~\ref{cor:imJsurvivors}. 
There are no elements in the $E_2$ term of the
Adams-Novikov spectral sequence which can kill $\beta_{i/j}[kq]$.  Such
elements would have to lie in Adams-Novikov filtration $0$.
\end{proof}


\section{Procedure for low dimensional calculations of root
invariants}\label{sec:proc}

In this section, we let $E=BP$, and concentrate on $BP$-filtered root
invariants.  
The following procedure is used in later sections to compute homotopy
root invariants from filtered root invariants using Theorem~\ref{thmA}.  
It only has a chance to
work through a finite range, and is very crude.  We state it mainly to
codify everything that must be checked in general to see that a filtered
root invariant detects a homotopy root invariant.

\begin{procedure}\label{procA}

Suppose we are in the situation where we know $\beta$ is an element
of $R^{[k]}_{BP}({\alpha})$, and we know that $\beta$ is a permanent cycle. 
Let $\iota_{-N}$ be the inclusion of the 
$-N$ cell of $P_{-N}$, where $-N = -N_i$ is the index of the bifiltration of
$\alpha$ that corresponds to the cell that bears the $k^\mathrm{th}$
filtered root invariant.  
Let $\nu_{-N} \alpha$ be the image of $\alpha$ in $\pi_*(P_{-N})$.
Then Theorem~\ref{thmA} tells us that $\nu_{-N} \alpha = \iota_{-N} \beta$
modulo elements in higher Adams-Novikov filtration.

\begin{description}

\item[Step 1]
Make a list
$$ \gamma_i [-n_i] $$
of additive generators in the $E_1$ term of the AAHSS for
$\ext(P_{-N})$, where $\gamma_i$ lies in Adams-Novikov 
filtration $s_i$ and
stem $k_i$, satisfying the following.
\begin{enumerate}

\item The homological degree $s_i$ satisfies
$$ s_i > \begin{cases}
k & n_i \equiv 1 \pmod q \\
k+1 & n_i \equiv 0 \pmod q
\end{cases} $$

\item The stem $k_i$ is greater than $p \cdot \abs{\alpha}$, and less than
$\abs{\beta}$.

\item We have $n_i = k_i - \abs{\alpha} + 1$, and $n_i \equiv 0,1 \pmod q$.

\end{enumerate}
These are precisely the conditions required for $\gamma_i[n_i]$ to be a  
candidate to survive in the AAHSS to an element in 
the same stem, but higher Adams-Novikov filtration, as $\iota_{-N} \beta$.
Condition (3) is a consequence of $P_{-N}$ having cells in dimensions
congruent to $0,-1$ modulo $q$.  We apply the first inequality in (2) because
Jones's theorem \cite{Miller} will preclude the possibility of these 
surviving to detect
the difference of $\iota_{-N} \br{\beta}$ and $\nu_{-N} \alpha$.
Thus we have
$$ \nu_{-N} \alpha = \beta[-N] + \sum_i a_i \gamma_i[-n_i] + (\text{ 
terms born on cells $< -(p-1) \abs{\alpha} + 1$}). $$

\item[Stem 2.]
Attempt to determine which $\gamma_i[n_i]$ are killed in the
filtration spectral sequence, or subsequently, in the Adams-Novikov 
spectral sequence for computing $\pi_*(P_{-N})$.  This will limit the
possibilities for what can detect the difference of $\iota_{-N} \br{\beta}$ 
and $\nu_{-N} \alpha$.

\item[Step 3]
First eliminate those $\gamma_i$ which are not permanent cycles
in the ANSS.  These cannot be root invariants.
Then, attempt to show that every remaining 
non-trivial linear combination of the 
elements $\gamma_i[n_i]$
which are not killed in the AAHSS actually
supports a non-trivial differential in the AHSS for 
computing $\pi_*(P_{-\infty})$.  
Suppose this is the case.
Then, if $\iota_{-N} \br{\beta} - 
\nu_{-N} \alpha$ is non-trivial, and is born on cells above the $-N$-cell,
it will project non-trivially to $\pi_*(P_{-M})$ for $M < N$, and will
represent the image of the homotopy root invariant.  
The homotopy root invariant is a permanent cycle in the 
AHSS for $\pi_*(P_{-\infty})$. 
Thus we may assume that the difference 
$\iota_{-N} \br{\beta} - \nu_{-N} \alpha$ is actually born on cells
above the $-(p-1) \abs{\alpha} + 1$ cell.  If this difference was 
non-trivial, we would violate Jones's theorem \cite{Miller}.
The difference is therefore trivial, and $\iota_{-N} \br{\beta}
= \nu_{-N} \alpha$.

\item[Step 4]
We will have shown that $\br{\beta} \in R(\alpha)$ if we know that 
$\iota_{-N} \br{\beta} \ne 0$.  It would suffice to show that 
$\beta[-N]$ is not killed in the AAHSS, and
also survives to a non-zero element in the ANSS
for computing $\pi_*(P_{-N})$.  Form a list of all of the elements
$$ \eta_i[n_i] $$
in the $E_1$-term of the AAHSS where $\eta_i$ is
in Adams-Novikov filtration $s_i$ and stem $k_i$ such that:
\begin{enumerate}
\item $s_i  \le \begin{cases}
k & -N \equiv -1 \pmod q \\
k-1 & -N \equiv 0 \pmod q
\end{cases}
$

\item $k_i \le \abs{\beta}$

\item $n_i = -N+\abs{\alpha}-k_i+1 \equiv 0, -1 \pmod q$

\end{enumerate}
These are precisely the elements
which have a 
chance of killing $\beta[-N]$ in the AAHSS, 
or subsequently in the ANSS.  
Next, show these elements actually do not kill $\beta[-N]$.  
\end{description}

\end{procedure}

\begin{rmk}\label{rmk:fss_ahss}
By comparing the definition of differentials in the AAHSS with the 
differentials in the AHSS, 
it is clear that if $\gamma$ detects $\br{\gamma}$ in the
ANSS for computing $\pi_*^S$, then if
$\gamma[n]$ supports a non-trivial AAHSS differential
$$ d_r(\gamma[n]) = \eta[m]$$
where $\eta$ is a permanent cycle in the ANSS,
then either
$\br{\gamma}$ supports a non-trivial AHSS $d_i$ for $i < r$, or there is an
AHSS differential
\begin{equation}\label{eq:AHSSdiff}
d_r (\br{\gamma}[n]) = \br{\eta}[m].
\end{equation}
where $\eta$ detects $\br{\eta}$.  However, the differential given in
Equation~\ref{eq:AHSSdiff} may be trivial in the AHSS if $\br{\eta}$ is the
target of a shorter differential.

Since in Step~$3$ of Procedure~\ref{procA} we are only
concerned with whether the elements $\br{\gamma}_i [n_i]$ support
non-trivial differentials
in the AHSS, it suffices to show that
the elements $\gamma_i[n_i]$ support non-trivial differentials in the AAHSS.  
One must
then make sure that the targets of the AAHSS differentials are not the
targets of shorter differentials in the AHSS.
\end{rmk}


\section{$BP$-filtered root invariants of some Greek letter
elements}\label{sec:anss}

In this section we compute the first two $BP$-filtered root
invariants of some chromatic families.  If
$\alpha = \alpha^{(n)}_{i/j_1,\ldots, j_k}$ is the $n^{\mathrm{th}}$ 
algebraic Greek
letter element of the ANSS which survives to a non-trivial element of
$\pi_*^S$, then one might expect that it should be the case that
\begin{equation}\label{eq:greekfiltroot}
\alpha^{(n+1)}_{i/j_1,\ldots, j_k} \dotin 
R^{[n+1]}_{BP}(\alpha^{(n)}_{i/j_1,\ldots, j_k})
\end{equation}
The purpose of this section is to show that Equation~\ref{eq:greekfiltroot}
holds for $n = 0$ at all
primes and for $n=1$ at odd primes.  Modulo an indeterminacy group which we
do not compute, we also show that the $n=1$ case of
Equation~\ref{eq:greekfiltroot} is true at the prime $2$.
Throughout this section, for a ring spectrum $E$, we shall let $\td{E}
\simeq \Sigma \br{E}$
denote the cofiber of the unit.  Also, whereas in Section~\ref{sec:aahss}
we
used Hazewinkel generators, \emph{in this section we always use Araki 
generators.}  This is because the $p$-series is more naturally expressed in
the Araki generators.

For appropriate $i$ and $j$, the elements which generate the $1$-line of
the ANSS are $\alpha_{i/j} \in BP_*BP$ and are given by
$$ \alpha_{i/j} = \frac{\eta_R(v_1)^i - v_1^i}{p^j}. $$

\begin{prop}\label{p's}
The first two filtered root invariants of $p^i$ are given by
\begin{align*}
(-v_1)^i & \in R^{[0]}_{BP}(p^i) \\
(-1)^i \cdot \alpha_i & \in R^{[1]}_{BP} (p^i)  
\end{align*}
\end{prop}

Define, for appropriate $i$ and $j$,
$$ \td \beta_{i/j} = 
\frac{(v_2 + v_1 t_1^p - v_1^p t_1)^i - v_2^i}{v_1^j} \in BP_*BP. $$
Observe that the image of $\td{\beta}_{i/j}$ in $BP_*BP/I_1$ coincides
with 
$$ \frac{\eta_R(v_2)^i - v_2^i}{v_1^j}. $$
It follows from the definition of the algebraic Greek letter elements
that, in the cobar complex,
$d_1(\td{\beta}_{i/j}) = p \cdot \beta_{i/j}$ unless $p = 2$ and
$i = j = 1$, in which case $\td{\beta}_1 = \alpha_{2/2}$, which is a
permanent cycle that detects $\nu$.

In order to compute the first filtered root invariants of the elements
$\alpha_{i/j}$, we first compute the $BP\wedge \td{BP}$-root invariants,
and then invoke Corollary~\ref{cor:E^Eroot}.  Note that $\td{BP} = 
\Sigma\br{BP}$, so the $BP \wedge \td{BP}$-root invariant is just the
suspension of the $BP\wedge \br{BP}$-root invariant.  The following result
holds for any prime.

\begin{prop}\label{prop:BP^BPalpha}
The $BP \wedge \td{BP}$-root invariant of $\alpha_{i/j}$ is given by
$$
R_{BP \wedge \td{BP}}(\alpha_{i/j}) = (-1)^{i-j} \td{\beta}_{i/j} + 
p\td{\beta}_{i/j} BP_* \td{BP}
$$
\end{prop}

\begin{cor}\label{cor:alpha's}
Suppose the Greek letter element $\alpha_{i/j}$ exists and is a permanent
cycle in the ANSS.  Let $\br{\alpha}_{i/j} \in \pi_*^S$ be the element that
$\alpha_{i/j}$ detects.
Then the zeroth filtered root invariant of $\br{\alpha}_{i/j}$ is
trivial.  If $p$ is odd, the first filtered root invariant of 
$\br{\alpha}_{i/j}$ is given by
$$
(-1)^{i-j} \td{\beta}_{i/j} + c \cdot \td{\alpha}_{i(p+1)+j} 
\in R^{[1]}_{BP}(\br{\alpha}_{i/j})
$$
where $c$ is some constant.  The second filtered root invariant is given by
$$
(-1)^{i-j} \beta_{i/j}  
\in R^{[2]}_{BP}(\br{\alpha}_{i/j}).
$$
\end{cor}

Deducing the first filtered root invariant in 
Corollary~\ref{cor:alpha's} from Proposition~\ref{prop:BP^BPalpha}
amounts to computing the
indeterminacy group $A$ described in Corollary~\ref{cor:E^Eroot}.  
The second filtered root invariant then follows from Theorem~\ref{thmB}.
We
describe a spectral sequence which computes $A$ (\ref{eq:res_AHSS}).  We
fully compute $A$
for odd primes $p$ and describe some of the $2$-primary aspects of this
computation in Remark~\ref{rmk:2alpha's}.

We have shown that
$\alpha_k \dotin R^{[1]}_{BP}(p^k)$.  It is shown in
\cite{MahowaldRavenel}, \cite{Sadofsky}
that these elements survive in the Adams-Novikov spectral sequence to
elements of $R(p^k)$.
Similarly, we have shown that $\beta_k \dotin R^{[2]}_{BP}(\br{\alpha}_k)$.
Again, in \cite{MahowaldRavenel}, \cite{Sadofsky} 
it is shown that these elements survive
to elements of $R(\br{\alpha}_k)$ for $p \ge 5$.  Sadofsky goes further in
\cite{Sadofsky} to show that $\beta_{p/2} \in R(\br{\alpha}_{p/2})$ for $p
\ge 5$.  

In general, we have shown that if $p^{k-1} \vert s$, then
$\beta_{s/k} \dotin R^{[2]}_{BP}(\alpha_{s/k})$.  According to the summary
presented in \cite[5.5]{Ravenel}, the elements $\beta_{s/k}$ exist and
are permanent cycles in the ANSS for $p \ge
5$ for these values of $k$ and $s$.  
However, without any additional
information, we can only deduce the conclusions of Theorem~\ref{thmA}.
In Section \ref{sec:anss3} we indicate what these
filtered root invariant calculations mean for the computation of homotopy
root invariants in low dimensions at the prime $3$.

In proving Propositions~\ref{p's} and \ref{prop:BP^BPalpha}, we shall need to
make use of the following well known computation (see Lemma~$2.1$ of 
\cite{AndoMoravaSadofsky}).

\begin{lem}\label{lem:tE_*}
Let $E$ be a complex oriented spectrum whose associated $p$-series 
$[p]_E(x)$ is not a zero divisor in $E^*[[x]]$.  Then the coefficient ring
of the Tate spectrum is given by
$$ \pi_*(tE) = (E_*[[x]]/([p]_E(x)))[x^{-1}] = E_*((x))/([p]_E(x)) $$
where the degree of $x$ is $-2$.  Furthermore, the inclusions of the
projective spectra are described as the inclusion of the fractional ideal
$$ x^N E_*[[x]]/([p]_E(x)/x) = E_*(\Sigma (B\ZZ/p)^{-2N})
\rightarrow \pi_*(tE) = E_*((x))/([p]_E(x)). $$ 
\end{lem}

According to Appendix B of \cite{GreenleesMay}, the skeletal and coskeletal
filtrations of $(E \wedge P)_{-\infty}$ give rise to the same notions of the 
$E$-root invariant.
It follows that $E$-root invariants are quite easy to compute for complex
orientable spectra $E$ provided one has some knowledge of the $p$-series of
$E$.  The method is outlined in the following corollary.

\begin{cor}\label{cor:Eroot}
Suppose that $\alpha$ is an element of $E_*$.  
Viewing $\br{\alpha}$ as the image of the 
constant power series $\alpha$
in $E_*((x))/([p]_E(x))$, suppose that $n$ is maximal so that 
$$ \br{\alpha} = a_n x^n + \mathcal{O}(x^{n +1}) $$
with $a_n \in E_*$ nonzero.  Then $a_n$ is an element of $R_E(\alpha)$.
\end{cor}

\begin{proof}[Proof of Proposition~\ref{p's}]
We first compute $R_{BP}(p^i)$.  The $p$-series $[p]_F(x)$ 
of the
universal $p$-typical formal group $F$ is given by
$$ [p](x) = px +_F v_1 x^p +_F v_2 x^{p^2} +_F \cdots $$
where $v_i$ are the Araki generators of $BP_*$.  In $tBP_*$, we have the
relation
$$ p = -v_1 x^{p-1} + \mathcal{O}(x^{p}) $$
which gives, upon taking the $i^{\rm th}$ power,
$$ p^i = (-v_1)^i x^{(p-1)i} + \mathcal{O} (x^{(p-1)i+1}). $$
Any other expression of $p^i$ in terms of $x^{(p-1)i}$ and higher order
terms will have a leading coefficient that differs by an element of the
ideal $ pBP_* \subset BP_*$.
Using Corollary~\ref{cor:Eroot} and Proposition~\ref{prop:Eroot}, 
we may conclude that
$$ (-v_1)^i + pBP_* = R_{BP}(p^i) \supseteq R^{[0]}_{BP}(p^i). $$
We calculate the Adams-Novikov $d_1$ on this coset as
$$ d_1(v_1^i) = \eta_R(v_1^i) - v_1^i. $$
The algebraic Greek letters $\alpha_i$ are defined by
$$ \alpha_i = \frac{\eta_R(v_1^i) - v_1^i}{p} $$
so we are in the situation where
$$ d_1 (R^{[0]}_{BP}(p^i)) \subseteq p \cdot (-1)^i \alpha_i + 
p^2 \alpha_i BP_*BP. $$
By Theorem~\ref{thmB}, we may conclude that
the first filtered root invariant $R^{[1]}_{BP}(p^i)$ is contained in the
coset
$ (-1)^i \alpha_i + p \alpha_i BP_*BP. $
\end{proof}

\begin{proof}[Proof of Proposition~\ref{prop:BP^BPalpha}]
The proof consists of two parts.  In Part~1 we prove the proposition with
$v_1$ inverted, and in Part~2 we prove that the $v_1^{-1}BP \wedge
\td{BP}$-root invariant can be lifted to compute the $BP \wedge \td{BP}$-root
invariant.
\vspace{.2in}

\noindent
{\it Part 1: computing the $v_1^{-1}BP \wedge \td{BP}$-root invariant}
\vspace{.15in}

We may calculate 
$tBP \wedge \td{BP}_*$ from $tBP_*BP$
using the split short
exact sequence
$$ 0 \rightarrow tBP_* \rightarrow 
tBP \wedge BP_* \rightarrow tBP \wedge \td{BP}_* 
\rightarrow 0. $$
Applying Lemma~\ref{lem:tE_*} to the left and right complex orientations of
$BP_*BP$, we get 
\begin{align*}
tBP \wedge BP_* 
& = BP\wedge BP_* ((x_L))
/([p]_L(x_L))  \\
& BP \wedge BP_* ((x_R))
/([p]_R(x_R)). 
\end{align*}
Here, the $p$-series $[p]_L(x)$, $[p]_R(x)$, and the coordinate $x_R$ are
given by the formulas
\begin{align*}
[p]_L (x) & = px +_{F_L} v_1 x^p +_{F_L} v_2 x^{p^2} +_{F_L} \cdots  \\
[p]_R (x) & = px +_{F_R} \eta_R (v_1) x^p +_{F_R} \eta_R(v_2) x^{p^2} 
+_{F_R} \cdots \\
x_L & = x_R +_{F_L} t_1 x_R^p +_{F_L} t_2 x_R^{p^2} +_{F_L} \cdots.
\end{align*}
The formal group laws $F_L$ and $F_R$ are the p-typical formal group laws 
over $BP_*BP$
induced by left and right units $\eta_L$ and $\eta_R$, respectively.
Consider the following computation in 
$tBP\wedge \td{BP}_*$.
\begin{equation}\label{eq:px_R}
\begin{split}
px_R & = [-1]_{F_R} \left(\eta_R(v_1) x_R^p +_{F_R} \eta_R(v_2) x_R^{p^2}
+_{F_R} \cdots \right) \\
& = -\eta_R (v_1) x_R^p - \eta_R(v_2) x_R^{p^2} + \mathcal{O}(x_R^{p^2+1})
\end{split}
\end{equation}
The last equality for $p > 2$ follows from the fact that for any
$p$-typical formal group law $G$ where $p > 2$, the $-1$-series is given by
$$ [-1]_G(x) = -x. $$
For $p = 2$ this is not true, but it turns out that the last equality of
Equation~\ref{eq:px_R} still
holds up to the power indicated in $tBP_*\td{BP}$.  The two expressions we 
want to be equal are written out explicitly below. 
\begin{multline}\label{eq:FERMAT1}
[-1]_{F_R} \left(\eta_R(v_1) x_R^2 +_{F_R} \eta_R(v_2) x_R^{4}
+_{F_R} \cdots \right)
\\
= (2t_1 - v_1)x_R^2 +
(14t_2 - 4t_1^3 - v_1 t_1^2 - 3v_1^2 t_1 - v_2 + v_1^3)x_R^4 +
\mathcal{O}(x_R^5) 
\end{multline}
\begin{multline}\label{eq:FERMAT2}
-\eta_R (v_1) x_R^2 - \eta_R(v_2) x_R^{4} + \mathcal{O}(x_R^{5})
\\
= (2t_1 - v_1)x_R^2 +
(14t_2 + 4t_1^3 - 13v_1 t_1^2 + 3v_1^2 t_1 - v_2)x_R^4 + 
\mathcal{O}(x_R^{5}) 
\end{multline}
The coefficient of $x_R^2$ is the same in Equations~\ref{eq:FERMAT1} and
\ref{eq:FERMAT2}.  Since we are working in
$tBP_*\td{BP}$, we are working modulo the $2$-series, and thus we only need
the coefficients of $x_R^4$ to be equivalent modulo $2$.  However, the
coefficient of $x_R^4$ in Equation~\ref{eq:FERMAT1} has an extra $v_1^3$,
but since we are working in the reduced setting of $tBP_*\td{BP}$, we have
$v_1 x_L^4 = 0$.  We may switch to $x_R$ since 
$x_L^4 = x_R^4 + \mathcal{O}(x_R^5)$.

Returning to our manipulations of the $p$-series,
upon dividing Equation~\ref{eq:px_R} by $x_R$, we get 
$$ p = -\eta_R (v_1) x_R^{p-1} - \eta_R(v_2) x_R^{p^2-1} +
\mathcal{O}(x_R^{p^2})$$
which implies
$$ \eta_R (v_1) = -\eta_R(v_2) x_R^{p^2-p} +
\mathcal{O}(x_R^{p^2-p+1})$$
since the image of $p$ is zero in $BP_* \td{BP}$.  Taking the $i^{\rm th}$ 
power and
exploiting the fact that the image of $v_1^i$ in $BP_*\td{BP}$ 
is also zero,
we may write
\begin{equation}\label{eqn:p^j_alpha}
p^j \alpha_{i/j} = \eta_R(v_1)^i - v_1^i = (-\eta_R(v_2))^i x_R^{i(p^2-p)} +
\mathcal{O}(x_R^{i(p^2-p)+1}).
\end{equation}
Since we are modding out by the $p$-series in $x_L$, we have the
relation
\begin{align*}
p x_L & = [-1]_{F_L} \left( v_1 x_L^p +_{F_L} v_2 x_L^{p^2} +_{F_L} \cdots
\right) \\
& = -v_1 x_L^p + \mathcal{O}(x_L^{p^2})
\end{align*}
or, dividing by $x_L$ and taking the $j^{\rm th}$ power,
\begin{equation}\label{eqn:p^j}
p^j = (-v_1)^j x_L^{j(p-1)} + \mathcal{O}(x_L^{j(p-1)+p^2-p})
\end{equation}
Upon combining
Equation~\ref{eqn:p^j} with Equation~\ref{eqn:p^j_alpha}, and
using the fact that $x_R^k = x_L^k + \mathcal{O}(x_R^{k+1})$, 
our original expression for $\alpha_{i/j}$ becomes 
the following. 
\begin{multline*}
(-v_1)^j\alpha_{i/j} + \alpha_{i/j} \mathcal{O}(x_R^{p^2-p})
\\
= (-1)^i(\eta_R(v_2)^i - v_2^i) x_R^{i(p^2-p) -
j(p-1)} + \mathcal{O}(x_R^{i(p^2-p) -j(p-1)+1})
\end{multline*}
Explicit formulas (see, for instance, \cite[4.3.21]{Ravenel}) reveal that
$$ \eta_R(v_2) = v_2 + v_1 t_1^p - v_1^p t_1 + p y $$
for some $y \in BP_*BP$.  Using this formula and Equation~\ref{eqn:p^j}
to write elements which are divisible by $p$ in terms of higher order
elements, we have
\begin{multline}\label{eq:v_1^jalpha_i/j}
(-v_1)^j\alpha_{i/j} 
+ \alpha_{i/j} \mathcal{O}(x_R^{p^2-p})
\\
= (-1)^i v_1^j \tilde \beta_{i/j} x_R^{i(p^2-p) -
j(p-1)} + \mathcal{O}(x_R^{i(p^2-p) -j(p-1)+1})
\end{multline}
In order to solve for $\alpha_{i/j}$, we would like to divide by $v_1^j$.
Therefore, we shall finish our algebraic manipulations in $tv_1^{-1}BP\wedge 
\td{BP}_*$.  We will then show that we can pull back our results to results
in $tBP\wedge \td{BP}_*$.  
Taking the image of Equation~\ref{eq:v_1^jalpha_i/j} in 
$tv_1^{-1}BP\wedge BP_*$
and dividing by $(-v_1)^j$, we get the expression
$$ \alpha_{i/j} = (-1)^{i-j}\td{\beta}_{i/j} x_R^{i(p^2-p) -
j(p-1)} + \alpha_{i/j}\mathcal{O}(x_R^{p^2 - p}) 
+ \mathcal{O}(x_R^{i(p^2-p) -j(p-1)+1})$$
By successively substituting the left hand side of this expression into
the right hand side, we obtain the expression below.
$$ \alpha_{i/j} = (-1)^{i-j}\tilde \beta_{i/j} x_R^{i(p^2-p) -
j(p-1)} + \text{higher order terms}.$$
We may conclude that $R_{v_1^{-1}BP\wedge \td{BP}}(\alpha_{i/j}) = 
\td{\beta}_{i/j}$, or that the $v_1^{-1}BP\wedge \td{BP}$-root invariant
lives in a higher stem.  
\vspace{.2in}

\noindent
{\it Part 2: lifting to the $BP \wedge \td{BP}$-root invariant}
\vspace{.15in}

We claim that the localization map
$$ BP \wedge \td{BP}_*(P_{-N}) \rightarrow v_1^{-1} BP \wedge \td{BP}_* 
(P_{-N}) $$
is an inclusion.  
In order to see this, chase the following diagram.
$$ 
\xymatrix{
 BP_*(P_{-N}) \ar[r] \ar@{^{(}->}[d] &  BP\wedge BP_*
(P_{-N})\ar[r]
\ar@{^{(}->}[d] &  BP \wedge \td{BP}_* (P_{-N})\ar[d]  
\\
v_1^{-1}  BP_*(P_{-N}) \ar[r]  & v_1^{-1} BP\wedge BP_* 
(P_{-N})\ar[r]
 & v_1^{-1}  BP \wedge \td{BP}_* (P_{-N}) 
}
$$
The relevant observations are that the other two localization maps in
the diagram are inclusions (since there is no $v_1$ torsion in these
groups), and 
the top and bottom sequences are compatibly split cofiber
sequences. 

We wish to deduce something about the $BP \wedge \td{BP}$ root
invariant from the computation of the $v_1^{-1} BP \wedge \td{BP}$
root invariant.
We refer to the following diagram.
$$
\xymatrix@C-1.4em{
v_1^{-1} BP \wedge
\td{BP}_* \ar[ddd] \ar@{~>}[rrr]^{R_{v_1^{-1}BP \wedge \td{BP}}(-)}
& & &  v_1^{-1} BP \wedge \td{BP}_*(S^{-N+1}) \ar[ddd]^\iota
\\
& *+[l]{BP \wedge \td{BP}_*} \ar[d] \ar[ul]
\ar@{~>}[r]^{R_{BP\wedge \td{BP}}(-)} & 
*+[r]{BP \wedge \td{BP}_*(S^{-N+1})} \ar[d]^\iota \ar[ur]
\\
& *+[l]{t BP \wedge \td{BP}_*} \ar[dl] \ar[r]_\nu & 
*+[r]{BP \wedge \td{BP}_*(\Sigma P_{-N})} \ar[dr] 
\\
t v_1^{-1} BP \wedge 
\td{BP}_* \ar[rrr]_\nu & & & v_1^{-1} BP \wedge \td{BP}_*(\Sigma P_{-N})}
$$
Here $N$ equals $2i(p^2 - p) - 2j(p-1)+1$.  We have shown that $\iota
((-1)^{i-j}\td{\beta}_{i/j}) = \nu(\alpha_{i/j})$ after $v_1$ is inverted.  
Because the
$v_1$-localization map is an inclusion, we may conclude that 
$ \iota ((-1)^{i-j}\td{\beta}_{i/j}) = \nu(\alpha_{i/j})$ in 
$BP \wedge \td{BP}_*(\Sigma P_{-N})$.
Therefore, we have computed the $BP_*\td{BP}$-root invariant
$$ R_{BP \wedge \td{BP}}(\alpha_{i/j}) = (-1)^{i-j} \td{\beta}_{i/j} + 
p\td{\beta}_{i/j} BP_* \td{BP}$$
or it lies in a larger stem.  The latter cannot be the case, however, as 
$\iota (\td{\beta}_{i/j})$ is non-zero.
\end{proof}

We wish to apply Corollary~\ref{cor:E^Eroot} to
Proposition~\ref{prop:BP^BPalpha} and 
conclude that 
$$ R^{[1]}_{BP}(\br{\alpha}_{i/j}) \subseteq (-1)^{i-j}\td{\beta}_{i/j} + 
p \td{\beta}_{i/j} BP_* \td{BP} + A. $$
Here $A$ is the image of the boundary homomorphism
$$ \pi_{iq-1}(W_0^{(0,1)}(P_{(-(ip-j)q-1, -(ip-j)q)})) \xrightarrow{\partial}
\pi_{iq-2}(W_1^1(S^{-(ip-j)q-1})). $$
Let $N$ be $(ip-j)q+1$.

The group $\pi_*(W_0^{(0,1)}(P_{(-N,-N+1)}))$ may be computed from the
AHSS associated to the filtration
\begin{equation*}
W_0^0(S^{-N}) \subset W_0^{(0,1)}(P^{-N+1}_{(-N,-N+1)}) \subset
W_0^{(0,1)}(P^{-N+2}_{(-N,-N+1)}) \subset \cdots.
\end{equation*}
The resulting spectral sequence takes the form
\begin{equation}\label{eq:res_AHSS}
E_1^{k,l} = 
\begin{cases}
0 & l \not\equiv 0,-1 \pmod q \\
\pi_k(W_0^1(S^l)) & l \equiv 0, -1 \pmod q,\:  l > -N \\
\pi_k(W_0^0(S^{-N})) & l = -N 
\end{cases} 
\end{equation}
and converges to $\pi_k(W_0^{(0,1)}(P_{(-N, -N+1)}))$.
The groups
$ \pi_*(W_0^1(S^l)) $ in the $E_1$-term of the spectral sequence
(\ref{eq:res_AHSS})
may be computed by taking the $0$ and $1$-lines of the $E_1$-term of 
the ANSS for the
sphere, and taking the cohomology with respect to the 
$d_1$ from the $0$-line to the $1$-line.  The differentials in 
spectral sequence (\ref{eq:res_AHSS})
are a restriction of the differentials in the AAHSS.  

The image of the map $\partial$ is generated by the image of 
Adams-Novikov $d_1$'s
supported on the $0$-line, the subgroup $pW_1^1(S^{-N})$, and 
the images of the higher differentials in the AAHSS whose sources are
permanent cycles in spectral sequence (\ref{eq:res_AHSS}) and whose 
targets are elements in
Adams-Novikov
filtration $1$ that are carried by the $-N$-cell.

We remark that for $0 \ne x \in \pi_*(W_1^1(S^0))$ which is \emph{not} a
permanent cycle, we have a non-trivial differential
$$ d_1(x[kq]) = px[kq-1]. $$
The reason the differential must be non-trivial is that there is no torsion
in the $E_1$-term of the ANSS, and if $px$ were the target of a $d_1$ 
in the ANSS, then $x$ would have to be a $d_1$-cycle.  Thus, if we are looking
for longer differentials in spectral sequence (\ref{eq:res_AHSS})
supported by $kq$-cells, we may restrict our search to those which are
actually $d_1$-cycles.

We now use spectral sequence (\ref{eq:res_AHSS}) to compute the
indeterminacy group $A$ for odd primes $p$.

\begin{proof}[Proof of Corollary~\ref{cor:alpha's}]
Since $filt_{BP}(\br{\alpha}_{i/j}) = 1$,
Lemma~\ref{filtlem} implies the filtered root invariant 
$R^{[0]}_{BP}(\br{\alpha}_{i/j})$ is
trivial.  In Proposition~\ref{prop:BP^BPalpha}, we found the $BP \wedge
\td{BP}$-root invariants $R_{BP \wedge \td{BP}}(\alpha_{i/j})$.
We will now deduce the first filtered root invariant
$R^{[1]}_{BP}(\br{\alpha}_{i/j})$ through the application of
Corollary~\ref{prop:E^Eroot}.  We just need to compute $A$ using
spectral sequence (\ref{eq:res_AHSS}).

The indeterminacy group $A$ is generated
by the images of Adams-Novikov $d_1$'s, the subgroup $pW_1^1(S^{-N})$, 
and the higher AAHSS differentials.  We just need to compute the latter.  
The
only elements of spectral sequence (\ref{eq:res_AHSS}) 
that could contribute to $A$ are those of the form
$$ \alpha_{k/l}[(i-k)q], \qquad k \le ip-j. $$  
Proposition~\ref{prop:imJdiffs_odd} tells us that only one can
contribute to $\delta$ and that contribution is 
given by 
$$ \partial(\alpha_{i(p+1)-j-l}[-(ip-j-l)q]) = 
\td{\alpha}_{i(p+1)-j}[-(ip-j)q-1] $$
where $l = \nu_p(i) + 1$.  Thus $A$ is also spanned by the element 
$\td{\alpha}_{-i(p+1)+j}$.  We conclude that
$$ (-1)^{i-j}\td{\beta}_{i/j} + c \cdot \td{\alpha}_{-i(p+1)+j} \in
R^{[1]}_{BP}(\br{\alpha}_{i/j}) $$
where $c$ is some constant.

To prove the second part of the proposition  we appeal to
Theorem~\ref{thmB}.  There is an Adams-Novikov differential 
$d_1 (\td{\beta}_{i/j}) = p \cdot \beta_{i/j}$.  The filtered root
invariant $R^{[1]}_{BP}(\alpha_{i/j})$ 
is carried by the $-N_1$-cell, where
$N_1 = 2(i(p^2 - p) - j(p-1)) + 1 = (ip - j)q + 1$.  The first cell to
attach nontrivially to this cell is the $-(ip - j)q$-cell, and
this is by the degree $p$ map.  We may conclude that 
$$ (-1)^{i-j} \beta_{i/j} \in R^{[2]}_{BP} (\br{\alpha}_{i/j}). $$ 
\end{proof}

\begin{rmk}\label{rmk:2alpha's}
Computing the group $A$ at the prime $2$ requires a more careful analysis.
The AAHSS differentials don't follow immediately from the $J$-spectrum AHSS
differentials computed in \cite{MahowaldJ}, because the varying
Adams-Novikov filtrations of the $v_1$-periodic elements.  
For instance, it turns out that
$$
\alpha_{4/4} = \td{\beta}_{2/2} + x_7 \in R^{[1]}_{BP} (\br{\alpha}_{2/2})
$$
where 
$$ x_7 \equiv v_2 t_1 + v_1(t_2 + t_1^3) \pmod 2. $$
If $i = j = 1$, then $\td{\beta}_1$ is a
permanent cycle which represents the element $\alpha_{2/2}$.
If $i = j = 2$, then $\td{\beta}_{2/2} + x_7$ is a permanent
cycle which represents the element $\alpha_{4/4}$.
Low dimensional calculations seem to indicate that in all other
circumstances we have 
$$ \beta_{i/j} + c \cdot \td{\alpha}_{3i-j+1} \alpha_1 
\in R^{[2]}_{BP}(\br{\alpha}_{i/j}) $$
where $c$ is some constant which may or may not be zero and $\td{\alpha}_k$
is equal to $\alpha_{i/j}$ with $j$ maximal.
The anomalous cases with $i=j=1,2$ correspond to the existence of the
`extra' Hopf
invariant $1$ elements $\nu$ and $\sigma$, which are detected in the ANSS
by $\alpha_{2/2}$ and $\alpha_{4/4}$, respectively.  These filtered root
invariants are thus consistent with the homotopy root invariant, which
takes each Hopf invariant $1$ element to the next one, if it exists.
\end{rmk}


\section{Computation of $R(\beta_1)$ at odd primes}\label{sec:beta1}

In this section we will compute the root invariant of $\beta_1$ at  odd
primes.  This result was stated, but not proved, in
\cite{MahowaldRavenel}.  We do this by first computing the $BP$-filtered
root invariants, and then by executing Step~$4$ of Procedure~\ref{procA}.

\begin{prop}
For $p >2$, the top filtered root invariant of $\beta_1$ is given by
$ R^{[2p]}_{BP}(\beta_1) \doteq \beta_1^p$. 
\end{prop} 

\begin{proof}[Sketch of Proof]
We first wish to show that $\beta_{p/p} \in R^{[2]}_{BP} (\beta_1)$.  In 
\cite{Ravenel}, it is shown that modulo the ideal $I_1$, the 
representatives for $\beta_1$ and $\beta_{p/p}$ in the cobar complex are
given by

\begin{gather*}
\beta_1 \equiv -\frac{1}{p} \sum_{0 < i < p} \binom{p}{i}
t_1^i|t_1^{p-i} \pmod {I_1}\\
\beta_{p/p} \equiv -\frac{1}{p} \sum_{0 < i < p^2} \binom{p^2}{i} t_1^i
|t_1^{p^2-i} \pmod {I_1}
\end{gather*}

Now, $\nu_p \binom{p^2}{i} = 1$ if and only if $p | i$.  Therefore, the
expression for $\beta_{p/p}$ modulo $I_1$ may be simplified.
$$ \beta_{p/p} \equiv - \frac{1}{p} \sum_{0 < i < p} \binom{p^2}{ip}
t_1^{ip} | t_1^{p(p-i)} \pmod {I_1}$$
Since $\binom{p}{i} \equiv \binom{p^2}{pi} \pmod p$ for $0 < i < p$,
we see that in $\ext_{BP_*BP/I_1}(BP_*)$, the element $\beta_{p/p}$ 
is obtained from
$\beta_1$ by application of the $0^\mathrm{th}$ algebraic Steenrod
operation. 
By computing the reduction map
$$ \ext_{BP_*BP}(BP_*) \rightarrow \ext_{BP_*BP/I_1} (BP_*/I_1) $$
we may conclude that  
$$ P^0 (\beta_1) \equiv \beta_{p/p} \pmod p$$
in $\ext_{BP_*BP}(BP_*)$.
By an algebraic analog of Jones's Kahn-Priddy theorem, this 
coincides with the $BP$-algebraic root invariant, which corresponds to
the first non-trivial filtered root invariant.  
Thus we have the filtered root invariant
$$ \beta_{p/p} \in R^{[2]}_{BP}(\beta_1). $$
In the ANSS, the element $\beta_{p/p}$ supports the
Toda differential.
$$ d_{q+1} (\beta_{p/p}) \doteq \beta_1^p \alpha_1 $$
The $-(p^2-p-1)q-1$ cell attaches to the $-(p^2-p)q-1$ cell of
$P_{-\infty}$ with
attaching map $\alpha_1$. 
Therefore, Theorem~\ref{thmB} tells us that
$$ \beta_1^p \dotin R^{[2p]}_{BP}(\beta_1). $$
There is no room for indeterminacy in the ANSS.
\end{proof}

\begin{cor}
For $p > 2$, the homotopy root invariant of $\beta_1$ is given by
$$ R(\beta_1) \doteq \beta_1^p. $$
\end{cor}

\begin{proof}
The element $\beta_1^p$ lies on the Adams-Novikov vanishing line, so it
must be the top filtered root invariant.  Therefore, we may apply
Corollary~\ref{corA} to see that either the image of the element
$\beta_1^p$ under the inclusion of the $(-(p^2-p-1)q-1)$-cell of 
$P_{-(p^2-p-1)q-1}$ is null, or $\beta_1^p$ actually detects the homotopy
root invariant.  We must therefore show that the element 
$\beta_1^p[-(p^2-p-1)q-1]$ in the AHSS for $P_{-(p^2-p-1)q-1}$ is not the
target of a differential.  We will actually show that 
$\beta_1^p[-(p^2-p-1)q-1]$ is not the target of a differential in the
AAHSS, and that there are no possible sources of differentials in the ANSS
for $P_{-(p^2-p-1)q-1}$ with target $\beta_1^p[-(p^2-p-1)q-1]$.

According to the low dimensional computations of the ANSS at odd primes
given in \cite[Ch. 4]{Ravenel}, the only elements in the $E_1$-term of the
AAHSS which can kill $\beta_1^p[-(p^2-p-1)q-1]$ in either the AAHSS or the
ANSS are the elements
\begin{align*}
\beta_k [-(k-1)(p+1)q] & \qquad 1 \le k < p \\
\alpha_{k/l} [-(k-p)q-1] & \qquad 1 \le k < p^2, 1 \le l \le \nu_p(k)+1
\end{align*}
as well as elements in $\alpha_1$-$\beta_1$ towers, i.e. those that
satisfy the hypotheses of Proposition~\ref{prop:beta_tower}.  These latter
elements cannot kill $\beta_1^p[-(p^2-p-1)q-1]$ in the AAHSS and cannot
survive to kill anything in the ANSS by Proposition~\ref{prop:beta_tower}.
Some care must be taken at $p=3$, but in this low dimensional range there
are no deviations from this pattern.

By Proposition~\ref{prop:odd_diffs},
the elements $\beta_k [-(k-1)(p+1)q]$ support non-trivial AAHSS
differentials
$$ d_2(\beta_k [-(k-1)(p+1)q]) \doteq \alpha_1 \cdot \beta_k [-k(p+1)q]). $$
According to Proposition~\ref{prop:imJdiffs_odd}, 
for $l < \nu_p(k)+1$, we have differentials in the AAHSS
$$ d_1(\alpha_{k/l+1} [-(k-p)q]) \doteq \alpha_{k/l} [-(k-p)q-1] $$
whereas for $l = \nu_p(k)+1$ and $k \ge 3$ we have
$$ d_5(\alpha_{k-2}[-(k-p-2)q]) \doteq \alpha_{k/l}[-(k-p)q-1]. $$
For $k = 2$ we have \cite{Thompson}
$$ d_4(1[pq-1]) = \alpha_2[(p-2)q-1]. $$
Finally, Proposition~\ref{prop:odd_diffs} implies that 
$$ d_q(\alpha_1[(p-1)q-1]) \doteq \beta_1[-1]. $$
There are no elements left to kill $\beta_1^p[-(p^2-p-1)q-1]$.
\end{proof}


\section{Low dimensional computations of root invariants 
at $p = 3$}\label{sec:anss3}

The aim of this section is to use knowledge of the
ANSS for $\pi_*^S$ in the first $100$ stems to compute the
homotopy root invariants of some low dimensional 
Greek letter elements $\alpha_{i/j}$ at $p=3$.
These results are
summarized in the following proposition.

\begin{prop}\label{prop:low_alpha}
We have the following root invariants at $p=3$.
\begin{align*}
R(\alpha_1) & = \beta_1 \\
R(\alpha_2) & \doteq \beta_1^2 \alpha_1 \\
R(\alpha_{3/2}) & = -\beta_{3/2} \\
R(\alpha_3) & = \beta_3 \\
R(\alpha_4) & \doteq \beta_1^5 \\
R(\alpha_5) & = \beta_5 \\
R(\alpha_{6/2}) & = \beta_{6/2} \\
R(\alpha_{6}) & = -\beta_6
\end{align*}
\end{prop}

It is interesting to note that although $\beta_2$ exists, it fails to be
the root invariant of $\alpha_2$.  The element $\beta_4$ does not exist, so
it cannot be the root invariant of $\alpha_4$.
In Section~\ref{sec:modcomp} we will prove that $\beta_i \dotin R(\alpha_i)$
for $i \equiv 0,1,5 \pmod 9$.  The remainder of this section is devoted to
proving Proposition~\ref{prop:low_alpha}.

\begin{figure}[tp]
\begin{center}
   \rotatebox{180}{\includegraphics{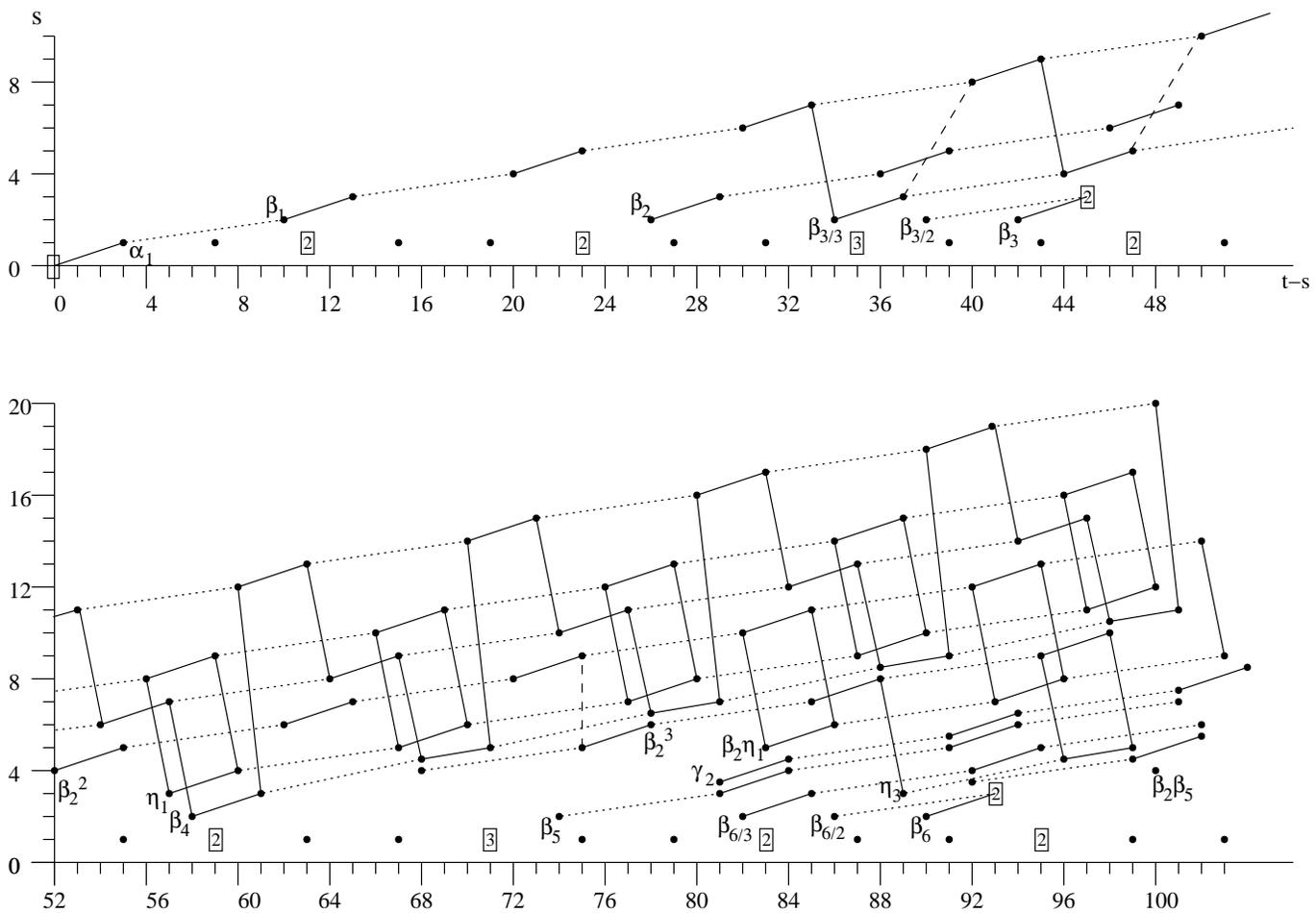}}
\end{center}
\caption{The Adams-Novikov spectral sequence at $p=3$}
\label{anss3E2fig}
\end{figure}

Figure~\ref{anss3E2fig} 
shows the Adams-Novikov $E_2$ term.
These
charts were created from the computations in \cite{Ravenel}.  Solid lines
represent multiplication by $\alpha_1$ and dotted lines represent the
Massey product $\bra{-, \alpha_1, \alpha_1}$.  Dashed lines represent hidden
extensions.

In the Section~\ref{sec:anss}, we proved in Corollary~\ref{cor:alpha's} 
that we have filtered root invariants
$$ (-1)^{i-j} \beta_{i/j} 
\in R^{[2]}(\alpha_{i/j}). $$ 
We supplement those results with two higher
filtered root invariants particular to $p=3$.

\begin{prop}
We have the following higher filtered root invariant of $\alpha_2$.
$$ \beta_1^2 \alpha_1 \dotin R^{[5]}_{BP}(\alpha_2) $$
\end{prop}

\begin{proof}
We know by Corollary~\ref{cor:alpha's} that
$$ -\beta_2 \in R^{[2]}_{BP}(\alpha_2). $$
The element $\beta_2$ is a permanent cycle in the ANSS, which detects an
element $\br{\beta}_2 \in \pi_*^S$.  We shall apply
Theorem~\ref{thmC} to determine the higher filtered root invariant from the
following hidden Toda bracket
$$ \beta_1^3 \doteq \bra{\alpha_1,3,\br{\beta}_2} $$ 
(see, for example, \cite{Ravenel}).  Considering the attaching map
structure of $P_{-\infty}$ we have the following equalities of Toda
brackets (the first is an equality of homotopy Toda brackets whereas the
second is an equality of Toda brackets in the ANSS).
\begin{gather*}
\bra{P^{-20}_{-25}} (\br{\beta}_2) 
\doteq \bra{\alpha_1, 3, \br{\beta}_2} 
\doteq \beta_1^3
\\
\bra{P^{-17}_{-25}} (\beta_1^2 \alpha_1) 
\doteq \bra{\alpha_1, \alpha_1, \beta_1^2\alpha_1} 
\doteq \beta_1^3
\end{gather*}
We conclude from Theorem~\ref{thmC}(3) that
$$ \beta_1^2 \alpha_1 \dotin R^{[5]}_{BP}(\alpha_2). $$
\end{proof}

\begin{prop}
We have the following higher filtered root invariants of $\alpha_4$.
\begin{align}
\beta_{3/3} \beta_1^2 & \dotin R^{[6]}_{BP}(\alpha_4) \\
\beta_1^5 & \dotin R^{[10]}_{BP}(\alpha_4)
\end{align}
\end{prop}

\begin{proof}
We know from Proposition~\ref{cor:alpha's} that 
$$ -\beta_4 \in R^{[2]}_{BP}(\alpha_4)$$
but 
$\beta_4$ is not a permanent cycle.  There
are Adams-Novikov differentials
\begin{align*}
d_5 (\beta_4) & \doteq \beta_{3/3}\beta_1^2 \alpha_1 \\
d_5 (\beta_{3/3} \beta_1^2) & \doteq \beta_1^5 \alpha_1.
\end{align*}
The higher filtered root invariants of $\alpha_4$ are derived from
two applications of Theorem~\ref{thmB}.  The relevant Toda brackets are
computed below.
\begin{gather*}
\bra{P^{-36}_{-40}}(\beta_{3/3} \beta_1^2) \doteq \alpha_1 \cdot
\beta_{3/3} \beta_1^2
\\
\bra{P^{-40}_{-44}}(\beta_1^5) \doteq \alpha_1 \cdot \beta_1^5
\end{gather*}
\end{proof}

We now apply Procedure~\ref{procA} to our filtered root invariants.
The results of executing Steps~1--3 in Procedure~\ref{procA} are
summarized in Table~\ref{tab:p=3}.  
The meanings of the contents of the columns are described below.
\begin{description}

\item[Element]
contains the
element $\alpha$ 
which we
wish to calculate the root invariant of.  

\item[$\mathbf{R^{[k]}}$]
indicates the $k^{\rm th}$ $BP$-filtered root invariant $\beta$ of
$\alpha$.  We would like $\beta$ to detect an element of
$R(\br{\alpha})$.  

\item[Cell]
indicates the cell of $P_{-\infty}$ that carries the filtered
root invariant.  

\item[$\mathbf{\{ \gamma_i\lbrack n_i \rbrack \} }$]
contains the list of elements of the
ANSS which
could survive to the difference between the filtered root invariant and
the homotopy root invariant.  This is the list described in Step~1 of
Procedure~\ref{procA}.  

\item[Diffs]
indicates the differentials in the AAHSS which kill
the elements $\gamma_i[n_i]$, as described in Step~2 of Procedure~\ref{procA}, 
or
the differentials supported by the $\gamma_i[n_i]$, as described
in Step~3 of Procedure~\ref{procA}.  If the element $\gamma_i$ supports a
non-trivial differential in the ANSS, we
shall place a ``(1)'' to indicate that it
is not a permanent cycle, and so cannot be a homotopy root invariant.
If the $\gamma_i[n_i]$ satisfies the conditions of
Proposition~\ref{prop:beta_tower}, then Steps~2 and 3 of 
Procedure~\ref{procA} are automatically satisfied, as explained in
Remark~\ref{rmk:fss_ahss}.  We shall indicate 
this with a ``(2)''.

\end{description}

\newpage

\tablecaption{Steps~1--3 of Procedure~\ref{procA} at
$p=3$}\label{tab:p=3}
\begin{center}
\tablefirsthead{\hline Element & $R^{[k]}$ & Cell & $\{ \gamma_i[n_i] \}$ & 
Diffs \\ \hline \hline}
\tablehead{
\multicolumn{5}{c}{Table~\ref{tab:p=3}, continued} \\
\hline Element & $R^{[k]}$ & Cell & $\{ \gamma_i[n_i] \}$ & 
Diffs \\ \hline \hline}
\tabletail{\hline}
\begin{supertabular}{c|c|c|c|c}
$\alpha_1$ & $\beta_1$ & $-2q$ & --- & --- 
\\
\hline
$\alpha_2$ & $\beta_1^2 \alpha_1$ & $-4q-1$ & --- & ---
\\
\hline
$\alpha_{3/2}$ & $\beta_{3/2}$ & $-7q$ & --- & --- 
\\
\hline
$\alpha_3$ & $\beta_3$ & $-8q$ & $\beta_2 \beta_1 \alpha_1[-7q-1]$ & $(2)$
\\
\hline
$\alpha_4$ & $\beta_1^5$ & $-9q$ & --- & --- 
\\
\hline
$\alpha_5$ & $\beta_5$ & $-14q$ & $\beta_4 \beta_1 \alpha_1[-13q-1]$ &
$(1)$
\\
& & & $ \eta_1 \beta_1[-12q-1]$ & 
$(1)$
\\
& & & $\eta_1 \beta_1 \alpha_1 [-13q]$ & 
$(1)$
\\
& & & $\beta_2^2 \beta_1[-11q]$ &
$(2)$
\\
& & & $\beta_{3/3} \beta_1^3 \alpha_1 [-12q-1]$ &
$(1)$
\\
& & & $\beta_2 \beta_1^3 \alpha_1[-10q-1]$ &
$(2)$
\\
& & & $\beta_2 \beta_1^4[-12q]$ &
$(2)$
\\
& & & $\beta_1^6 \alpha_1 [-11q-1]$ &
$(2)$
\\
& & & $\beta_1^7 [-13q]$ &
$(2)$
\\
\hline
$\alpha_{6/2}$ & $\beta_{6/2}$ & $-16q$ & $\beta_2 \eta_1[-15q-1]$ &
$(1)$
\\
& & & $\beta_2^3/\alpha_1[-13q-1]$ &
$(2)$
\\
& & & $\beta_4 \beta_1 \alpha_1[-12q-1]$ & 
$(1)$
\\
& & & $\beta_2^3[-14q]$ & $(2)$
\\
& & & $\eta_1 \beta_1 \alpha_1[-12q]$ & $(1)$
\\
& & & $\beta_2^2 \beta_1^2 \alpha_1[-13q-1]$ & $(2)$
\\
& & & $\beta_2^2 \beta_1^3 [-15q]$ & $(2)$
\\
& & & $\beta_{3/3} \beta_1^4 [-13q]$ & $(2)$
\\
& & & $\beta_2 \beta_1^5 \alpha_1[-14q-1]$ & $(2)$
\\
& & & $\beta_1^7[-12q]$ & $(2)$
\\
& & & $\beta_1^8 \alpha_1[-15q-1]$ & $(2)$
\\
\hline
$\alpha_6$ & $\beta_6$ & $-17q$ & $\eta_1 \beta_2[-15q-1]$ & $(1)$
\\
& & & $\beta_2^3 / \alpha_1 [-13q-1]$ & $(2)$
\\
& & & $\beta_4 \beta_1 \alpha_1 [-12q-1]$ & $(1)$
\\
& & & $\eta_1 \beta_2 \alpha_1 [-16q]$ & $(1)$
\\
& & & $\beta_2^3[-14q]$ & $(2)$
\\
& & & $\eta_1 \beta_1 \alpha_1[-12q]$ & $(1)$
\\
& & & $\eta_1 \beta_1^3[-16q-1]$ & $(1)$
\\
& & & $\beta_2^2 \beta_1^2 \alpha_1$ & $(2)$
\\
& & & $\beta_2^2 \beta_1^3[-15q]$ & $(2)$
\\
& & & $\beta_{3/3} \beta_1^4[-13q]$ & $(1)$
\\
& & & $\beta_{3/3} \beta_1^5 \alpha_1[-16q-1]$ & $(2)$
\\
& & & $\beta_2 \beta_1^5 \alpha_1[-14q-1]$ & $(2)$
\\
& & & $\beta_2 \beta_1^6 [-16q]$ & $(2)$
\\
& & & $\beta_1^8 \alpha_1[-15q-1]$ & $(2)$
\\
\end{supertabular}
\end{center}
\vspace{.25in}

Some care must be taken in the case of $\alpha_3$.  If $\bra{\alpha_1, 3,
\beta_3}$ is nonzero (i.e. contains $\pm \beta_2 \beta_1^2$), then the
AAHSS differential
$$ d_4(\beta_2 \beta_1 \alpha_1[-7q-1]) \doteq \beta_2 \beta_1^2 [-9q-1] $$
may be trivial in the AHSS with an intervening differential supported by
$\beta_3[-8q]$ (see Remark~\ref{rmk:fss_ahss}).  
This does not happen, though, as the Toda bracket 
$\bra{\alpha_1, 3,\beta_3}$ is actually zero.  It is easily seen to have
zero indeterminacy.  If the Toda bracket were non-trivial, it would have to
contain $\pm \beta_2 \beta_1^2$, and thus 
$\alpha_1 \cdot \bra{\alpha_1, 3,\beta_3}$ would be nonzero.  This cannot
be the case, since we have
$$ \alpha_1 \cdot \bra{\alpha_1, 3,\beta_3}  
= - \bra{\alpha_1, \alpha_1, 3} \cdot \beta_3
= \alpha_2 \cdot \beta_3 = 0. $$

Now we must complete Step~4 of Procedure~\ref{procA}.  We must determine
that the filtered root invariants survive to non-trivial elements in the
AHSS, and subsequently in the ANSS. 
For this we must determine which elements in the
AAHSS (or subsequently in the ANSS) 
can support differentials which kill the filtered
root invariants.  

Proposition~\ref{prop:beta_survive} implies that for all of the
$\alpha_{i/j}$ above whose filtered root invariants were $\beta_{i/j}$,
the image of $\beta_{i/j}$ in the appropriate stunted projective space
is non-trivial.  Thus we conclude that 
$$ \beta_{i/j} \dotin R(\alpha_{i/j}) $$
for $i \le 6$ and $i \ne 2,4$.  

Step~4 of Procedure~\ref{procA} is completed for $\alpha_2$ and $\alpha_4$
in the following lemmas.  In both of these lemmas, we denote by $\br{x}$
the element in homotopy detected by the ANSS element $x$.  There is no
ambiguity arising from higher Adams-Novikov filtration for 
the elements that we will be considering.

\begin{lem}\label{lem:alpha2_step_4}
The element $\br{\beta_1^2\alpha_1}[-4q-1]$ is not the target of a differential
in the AHSS for $P_{-4q-1}$.
\end{lem}

\begin{proof}
Differentials which could kill $\br{\beta_1^2\alpha_1}[-4q-1]$ must have
their source in the $7$-stem.
Proposition~\ref{prop:imJdiffs_odd} demonstrates that in the the targets of any
AHSS differentials supported by $ImJ$ elements in the $7$-stem 
are also $ImJ$ elements.
In the range we are considering, there is no room for shorter differentials.
The only element left which could kill $\br{\beta_1^2\alpha_1}[-4q-1]$ is
$\br{\beta_1^2}[-3q-1]$.  However, the complex $P^{-3q-1}_{-4q-1}$ is reducible,
so there can be no such differential.
\end{proof}

\begin{lem}
The element $\br{\beta_1^5}[-9q]$ is not the target of a differential
in the AHSS for $P_{-9q}$.
\end{lem}

\begin{proof}
Exactly as in the proof of Lemma~\ref{lem:alpha2_step_4},
the differentials supported by $ImJ$ elements in the $15$-stem can only hit
other $ImJ$ elements.  The only elements left which can support
differentials are given below.
\begin{gather*}
\br{\beta_1^2}[-q-1] \qquad 
\br{\beta_1^2 \alpha_1}[-2q] \qquad 
\br{\beta_2\beta_1}[-5q-1]
\\
\br{\beta_1^4}[-6q-1] \qquad 
\br{\beta_{3/3}\beta_1\alpha_1}[-8q]
\end{gather*}
All but $\br{\beta_{3/3}\beta_1\alpha_1}[-8q]$ are the target or source of
the AHSS $d_2$ differentials displayed below.
$$ \xymatrix@R-2em{
& \br{\beta_1^2}[-q-1] \ar@{|->}[r] & \br{\beta_1^2\alpha_1}[-2q-1]
\\
\br{\beta_1^2}[-q] \ar@{|->}[r] & \br{\beta_1^2 \alpha_1}[-2q]
\\
& \br{\beta_2\beta_1}[-5q-1] \ar@{|->}[r] & \br{\beta_2 \beta_1
\alpha_1}[-6q-1]
\\
\br{\beta_{3/3}\alpha_1}[-5q-1] \ar@{|->}[r] & \br{\beta_1^4}[-6q-1]
} $$
The last of these is the result of a hidden $\alpha_1$-extension in the
ANSS.  The remaining element $\br{\beta_{3/3}\beta_1\alpha_1}[-8q]$
cannot
support a differential whose target is $\br{\beta_1^5}[-9q]$, because 
$P^{-8q}_{-9q}$ is coreducible.
\end{proof}


\section{Algebraic filtered root invariants}~\label{sec:mss}

In this section we will describe the Mahowald spectral sequence (MSS),
which is a spectral sequence that computes $\ext(H_*X)$ by applying
$\ext(H_*-)$ to an $E$-Adams resolution of $X$.  This spectral sequence is
described in \cite{Miller81}.  We will briefly recall its construction for the
reader's convenience.  We will then define algebraic $E$-filtered root
invariants, and indicate how some of the results of Section~\ref{sec:results} 
carry
over to the algebraic setting to compute algebraic root invariants from
the algebraic filtered root invariants in the MSS. 

Let $E$ be a ring spectrum and suppose that $H$ is an $E$-ring spectrum in
the sense that there is a map of ring spectra $E \rightarrow H$.  Suppose
furthermore that $H_*H$ is a flat $H_*$ module.  Let $X$ be a (finite) 
spectrum.  
By applying $H_*$ to the
cofiber sequences that make up the $E$-Adams resolution, we obtain short exact
sequences
\begin{equation}\label{eq:MSSSES}
0 \rightarrow H_*(W_s(X)) \rightarrow H_*(W_s^s(X))
\rightarrow H_*(\Sigma W_{s+1}(X)) \rightarrow 0
\end{equation}
which are split by the $E$-action map
$$ H \wedge W^s_s(X) = H \wedge E \wedge \br{E}^{(s)} \wedge X \rightarrow
H \wedge \br{E}^{(s)} \wedge X = H \wedge W_s(X). $$
Thus the short exact sequences~\ref{eq:MSSSES} give rise to long
exact sequences when we apply $\ext_{H_*H}(-)$.  Therefore, upon applying
$\ext$ to the $E$-Adams resolution, we get a spectral sequence
$$ E_1^{s,t,k}(H_*X) = \ext^{s-k,t}(H_*W_k^k(X)) \Rightarrow \ext^{s,t}(H_*X)
$$
called the Mahowald spectral sequence (MSS).  For the remainder of this
section we shall assume $H = H\FF_p$.  

Many of the notions that we defined for the homotopy root invariant carry
over to the algebraic context.  We may define Tate comodules (over the dual
Steenrod algebra)
$$ t(H_*X) = \varprojlim_N (H_*(X \wedge \Sigma P_{-N})) = (tH\wedge X)_*. $$
Here, it is essential that the limit is taken \emph{after} taking homology.
We may use this to define algebraic
$E$-root invariants.

\begin{defn}[\it Algebraic $E$-root invariant \rm]
Let $\alpha$ be an element of the Ext group
$\ext^{s,t}(H_*X)$.
We have the following diagram of $\ext$ groups which defines the 
algebraic $E$-root invariant $R_{E,alg}(\alpha)$.
$$\xymatrix@C+1em{ 
\ext^{s,t}(H_*X) \ar[d]_{f} \ar@{~>}[r]^-{R_{E,alg}(-)} &
\ext^{s,t}(H_*E \wedge \Sigma^{-N+1}X) \ar[d]_{\iota_N}
\\
\ext^{s,t}(t(H_*(E \wedge X))) \ar[r]_-{\nu_{N}} &
\ext^{s,t}(H_*(E \wedge \Sigma P_{-N} \wedge X))
}$$
Here $f$ is induced by the inclusion of the $0$-cell of
$tS^0$, $\nu_N$ is the 
projection onto the $-N$-coskeleton,
$\iota_N$ is inclusion of the $-N$-cell,
and $N$ is minimal with respect to the property that
$\nu_{N} \circ f(\alpha)$ is non-zero.  Then the algebraic $E$-root 
invariant
$R_{E,alg}(\alpha)$ is defined to be the coset of lifts 
$\gamma \in \ext^{s, t+N-1}(H_*E \wedge X)$ of the
element $\nu_{N} \circ f (\alpha)$.  It could be the the case that
$f(\alpha) = 0$, in which case we say that the algebraic $E$-root
invariant is trivial.
\end{defn}

For simplicity, we now restrict our attention to the case $E = BP$, and we
work at an odd prime $p$.  Let
$\br{\Lambda}$ be the periodic lambda algebra \cite{Gray}.  The cohomology
of $\br{\Lambda}$ is $\ext(H_*)$.  The MSS has a concrete description in
terms of $\br{\Lambda}$.  Define a decreasing filtration on $\br{\Lambda}$
by
$$ \br{\Lambda} = F_0 \br{\Lambda} \subset F_1 \br{\Lambda} \subset F_2
\br{\Lambda} \subset \cdots $$
where $F_k \br{\Lambda}$ is the subcomplex of $\br{\Lambda}$ generated by
monomials containing $k$ or more $\lambda_i$'s.  The spectral sequence
associated to the filtered complex $\{F_k \br{\Lambda} \}$ is isomorphic to
the MSS \emph{starting with the $E_2$-term}.  The spectral sequence of the
filtered complex would agree with the MSS on the level of $E_1$-terms if
one were to apply $\ext(H_*-)$ to the correct $BP$-resolution (which
differs from the canonical $BP$-resolution).  Therefore, we shall also refer
to the spectral sequence associated to the filtered complex 
$\{F_k \br{\Lambda} \}$ as a MSS.

In analogy with the spaces $W_s(X)$, define subcomplexes of $\br{\Lambda}
\otimes H_*(X)$ (the complex which computes $\ext(H_*X)$) by
$$ W_k(H_*X) = F_k\br{\Lambda} \otimes H_*(X). $$
We define quotients $W_k^l(H_*X)$ by
$$ W_k^l(H_*X) = W_k(H_*X)/W_{l+1}(H_*X). $$
The cohomology $H^*(W_k^l(H_*X))$ is computed by restricting the MSS by
setting the $E_1^{*,*,i}$-term equal to zero for $i < k$ or $i > l$.
We define complexes $W_k^l(H_*P^N)$ to be the inverse limit
$$ W_k^l(H_*P^N) = \varprojlim_M W_k^l(H_*P^N_{-M}). $$

For sequences
\begin{align*}
I & = \{ k_1 < k_2 < \cdots < k_l \} \\
J & = \{ N_1 < N_2 < \cdots < N_l \}
\end{align*}
we can define 
filtered Tate complexes
$$ W_I(H_*P^J \wedge X) = \sum_i W_{k_i}(H_*P^{N_i}\wedge X). $$
Given another pair of sequences $(I',J') \le (I,J)$, we
define complexes
$$ W_I^{I'}(H_*P^J_{J'}\wedge X) = W_I(H_*P^J \wedge
X)/W_{I'+1}(H_*P^{J'-1}\wedge X) $$
where $I'+1$ (respectively $J'-1$) is the sequence obtained by increasing 
(decreasing)
every element of the sequence by $1$.

We shall now define algebraic filtered root invariants in analogy with
the definition of filtered root invariants given in
Section~\ref{sec:defs}.
Let $\alpha$ be an element of $\ext^{s,t}(H_*X)$.
We shall describe a pair of sequences 
\begin{align*}
I & = \{ k_1 < k_2 < \cdots < k_l \} \\
J & = \{ -N_1 < -N_2 < \cdots < -N_l \}
\end{align*}
associated to $\alpha$, which we define inductively.  
Let $k_1 \ge 0$ be
maximal such that the image of $\alpha$ under the composite
$$ \ext^{s,t-1}(H_*\Sigma^{-1}X) 
\rightarrow H^{s,t-1}(W_0^{k_1-1}(H_*P \wedge X)_{-\infty}) $$
is trivial.
Next, choose 
$N_1$ to be maximal such that the image of $\alpha$ under the composite
$$ \ext^{s,t-1}(H_*\Sigma^{-1} X) 
\rightarrow H^{s,t-1}(W_0^{(k_1-1, k_1)}(H_*P_{(-N_1+1, \infty)} \wedge X)) $$
is trivial.  Inductively, given 
\begin{align*}
I' & = (k_1, k_2, \ldots, k_i) \\
J' & = (-N_1, -N_2, \ldots, -N_i)
\end{align*}
let $k_{i+1}$ be maximal so that the image of $\alpha$ under the composite
$$ \ext^{s,t-1}(H_*\Sigma^{-1} X) 
\rightarrow H^{s,t-1}(W_0^{(I'-1,k_{i+1}-1)}(H_*P_{(J'+1,\infty)} \wedge
X)) $$
is trivial.  If there is no such maximal $k_{i+1}$, we declare that
$k_{i+1} = \infty$ and we are finished.  Otherwise, choose $N_{i+1}$ to be
maximal such that the composite
$$ \ext^{s,t-1}(H_*\Sigma^{-1} X) 
\rightarrow H^{s,t-1}(W_0^{(I'-1,k_{i+1}-1, k_{i+1})}(H_*P_{(J'+1,-N_{i+1}+1,
\infty)}\wedge X)) $$
is trivial, and continue the inductive procedure.  We shall refer to the
pair $(I,J)$ as the $BP$-bifiltration of $\alpha$.

Observe that there is an exact sequence
$$ H^{s,t-1}(W_I(H_*P^J \wedge X)) \rightarrow 
\ext^{s,t-1}(t(H_*\Sigma^{-1} X))
\rightarrow
H^{s,t-1}(W^{I-1}(H_*P_{J+1} \wedge X)). $$
Our choice of $(I,J)$ ensures that the image of $\alpha$ in 
$H^{s,t-1}(W^{I-1}(H_*P_{J+1} \wedge X))$ is trivial.  
Thus $\alpha$ lifts to an
element $f^\alpha \in H^{s,t-1}(W_I(H_*P^J \wedge X))$.

\begin{defn}[{\it Algebraic filtered root invariants} \rm]
Let $X$ be a finite complex and let $\alpha$ be an element of
$\ext^{s,t}(H_*X)$ of
$BP$-bifiltration $(I,J)$.  Given a lift 
$f^\alpha \in H^{s,t-1}(W_I(H_*P^J \wedge X))$, the $k^\mathrm{th}$
algebraic filtered
root invariant is said to be trivial if $k \ne k_i$ for any $k_i \in I$.
Otherwise, if $k=k_i$ for some $i$, 
we say that the image $\beta$ of $f^\alpha$ under the quotient map
$$ H^{s,t-1}(W_I(H_*P^J \wedge X)) \rightarrow
H^{s,t-1}(W_{k_i}^{k_i}(H_*\Sigma^{-N_{i}} X)) $$
is an element of the \emph{$k^\mathrm{th}$ algebraic 
filtered root invariant} of 
$\alpha$.
The $k^\mathrm{th}$ algebraic filtered root invariant is the coset 
$R^{[k]}_{E,alg}(\alpha)$ of the MSS $E_1$-term
$E_1^{s,t+N_i-1,k}(H_*X)$ of all such 
$\beta$ as we vary the 
lift $f^\alpha$.
\end{defn}

We wish to indicate how our filtered root invariant theorems carry over to
the algebraic context.  In order to do this we must produce algebraic
versions of $K$-Toda
brackets.  Suppose that $M$ is a finite $A_*$-comodule concentrated in
degrees $0$ through $n$ with a single
$\FF_p$ generator in degrees $0$ and $n$.
In what follows, we let $M^j$ be the sub-comodule of $M$ consisting of
elements of degree less than or equal to $j$, and let $M_i^j$
be the quotient $M^j/M^{i-1}$.  We shall omit the top index for the 
quotient $M_i = M/M^{i-1}$.

\begin{defn}[{\it Algebraic $M$-Toda bracket} \rm]
Let
$$ f: \ext^{s,t}(M_1) \rightarrow \ext^{s+1,t}(\FF_p) $$
be the connecting homomorphism associated to the short exact sequence
$$ 0 \rightarrow \FF_p \rightarrow M \rightarrow M_1 \rightarrow 0. $$
Let $\nu : M_1 \rightarrow \Sigma^n \FF_p$ be the projection onto the top
generator.
Suppose $\alpha$ is an element of $\ext^{s,t}(H_*X)$.
We have
$$ \ext^{s,t}(H_*X) \xleftarrow{\nu_*} \ext^{s,t+n}(H_*X \otimes M_1)
\xrightarrow{f} \ext^{s+1,t+n}(H_*X). $$
We say the algebraic $M$-Toda bracket 
$$ \bra{M}(\alpha) \subseteq \ext^{s+1,t+n}(H_*X) $$
is \emph{defined} if $\alpha$ is in the image of $\nu_*$.  Then the
algebraic $M$-Toda bracket is the 
collection of all $f(\gamma) \in \ext^{s+1,t+n}(H_*X)$ where $\gamma \in 
\ext^{s,t+n}(H_*X \otimes M_1)$ is any element satisfying $\nu_*(\gamma) =
\alpha$.
\end{defn}

We leave it to the reader to construct the dually defined algebraic
$M$-Toda bracket analogous to the dually defined $K$-Toda bracket, 
and MSS version of the algebraic $M$-Toda bracket 
analogous to the $E$-ASS version of the $K$-Toda bracket.  The definitions
of these Toda brackets given in Section~\ref{sec:toda} carry over verbatim
to the algebraic context.

With these definitions, our filtered root invariant results given in
Section~\ref{sec:results} have analogous algebraic statements.   
The first algebraic filtered root invariant is the algebraic $BP$-root
invariant.  If the $k^\mathrm{th}$ 
algebraic filtered root invariant is a permanent cycle in
the MSS, then it is the algebraic root invariant modulo $BP$-filtration
greater than $k$.  One may deduce higher algebraic filtered root invariants
from differentials and compositions in the MSS.  In other words, there are
algebraic versions of Theorem~\ref{thmA}, Corollary~\ref{corA},
Theorem~\ref{thmB}, Theorem~\ref{thmC}, Proposition~\ref{prop:Eroot},
Proposition~\ref{prop:E^Eroot}, and Corollary~\ref{cor:E^Eroot}, where one
makes the following replacements:
\begin{itemize}

\item Replace $\pi_*$ with $\ext$.
\item Replace root invariants with algebraic root invariants.
\item Replace $E$-root invariants with algebraic $E$-root invariants.
\item Replace filtered root invariants with algebraic filtered root
invariants.
\item Replace $K$-Toda brackets with algebraic $H_*(K)$-Toda brackets.
\item Specialize from $E$ to $BP$.
\item Replace the $E$-ASS with the MSS.

\end{itemize}


\section{Modified filtered root invariants}\label{sec:mod}

In this section we will recall the modified root invariant that Mahowald and
Ravenel define in \cite{MahowaldRavenel}.  
We will then explain how to define modified versions of filtered root
invariants, and how to adapt our theorems to compute modified root invariants.
We also discuss algebraic modified filtered root
invariants, which is the tool we will be using in
Section~\ref{sec:modcomp}.

Let $X$ be a finite complex.
We shall begin be recalling the definition of the modified root invariant
$R'(\alpha)$ for an element
$\alpha \in \pi_t(X)$.  
Let $V(0)$ be the mod $p$ Moore spectrum.  
We define modified versions of the stunted Tate spectra $P_M^N$ which will
have the same underlying space, but which we shall view as being built out
of $V(0)$ instead of the sphere spectrum.  To this end, define spectra
$$ {P'}^N = \begin{cases}
P^{N+1} & \text{if $N$ is of the form $kq-1$} \\
P^N & \text{otherwise.}
\end{cases}
$$
Likewise, for $M \le N$, define quotient spectra
$$ {P'}^N_M = {P'}^N/{P'}^{M-1}.  $$
Thus the spectrum ${P'}^N_M$ has a $V(0)$-cell in every dimension from $M$
to $N$ congruent to $-1$ modulo $q$.
There are cofiber sequences
$$ {P'}_{kq-1}^{(l-1)q-1} \rightarrow {P'}_{kq-1}^{lq-1} 
\rightarrow \Sigma^{lq-1}
V(0).$$

\begin{defn}[{\it Modified root invariant} \rm]
Let $X$ be a finite complex, and suppose we are given $\alpha \in
\pi_t(X)$.  The \emph{modified root invariant} 
of $\alpha$ is the coset 
of all dotted
arrows making the following diagram commute.
$$ \xymatrix{
S^{t} \ar@{.>}[r] \ar[d]^\alpha & \Sigma^{-N+1} V(0) \wedge X \ar[dd]  
\\
X \ar[d]
\\
tX \ar[r] & \Sigma P'_{-N} \wedge X 
}$$
This coset is denoted $R'(\alpha)$.  Here $N$ is chosen to be minimal
such that the composite $S^{t} \rightarrow \Sigma P'_{-N} \wedge X$ 
is non-trivial.
\end{defn}

Thus the definition of the modified root invariant differs from the
definition of the root invariant in that it takes values in
$\pi_*(V(0)\wedge X)$ instead of $\pi_*(X)$, and that $P$ has been replaced
with $P'$.  By making these replacements elsewhere, we may produce modified
versions of our other definitions, as summarized
below.  In the case of the algebraic filtered root invariants, we choose to
use the quotient $\br{\Lambda}_{(0)} = \br{\Lambda}/(v_0)$ described in
\cite{Gray}.  The cohomology of $\br{\Lambda}_{(0)}$ is given by
$$ H^*(\br{\Lambda}_{(0)}) = \ext(H_*V(0)). $$
One has modified versions of the subcomplexes $W_k(H_*X)$ given by
$$ W'_k(H_*X) 
= Im(W_k(H_*X) \hookrightarrow \br{\Lambda}\otimes H_*(X) \rightarrow
\br{\Lambda}_{(0)} \otimes H_*(X)). $$
By modifying our definitions, we produce:
\begin{itemize}
\item Modified $E$-root invariants
$$ R'_E: \pi_*(X) \rightsquigarrow \pi_*(E \wedge V(0) \wedge X) $$
\item Modified algebraic root invariants
$$ R'_{alg}: \ext(H_*X) \rightsquigarrow \ext(H_*V(0) \wedge X) $$
\item Modified filtered Tate spectra
$$W_I({P}'^J)$$
\item Modified filtered root invariants 
$$ {R_E^{[k]}}': \pi_*(X) \rightsquigarrow \pi_*(W_k^k(V(0) \wedge X)) $$
\item Modified algebraic $E$-root invariants
$$ R'_{E,alg} : \ext(H_*X) \rightsquigarrow \ext(H_*E \wedge V(0) \wedge X)
$$
\item Modified algebraic filtered root invariants
$$ {R^{[k]}_{E,}}'_{alg} : \ext(H_*X) \rightsquigarrow H^*({W'}_k^k(H_*X)) $$
\end{itemize}

We would like to reformulate the results of Section~\ref{sec:results} in
terms of modified root invariants.  We need a modified version of the Toda
brackets which appear in the statements of the main theorems.  Suppose that
$K$ is a finite complex built out of $V(0)$ with a single bottom $V(0)$-cell in
dimension $0$ and a single top $V(0)$-cell in dimension $n$.  
Let $K^j$ be the $j^\mathrm{th}$ $V(0)$-skeleton of $K$, and let $K_i^j$
be the quotient $K^j/K^{i-1}$.  We shall omit the top index  for the 
$i^\mathrm{th}$ $V(0)$-coskeleton $K_i = K/K^{i-1}$.

\begin{defn}[{\it Modified $K$-Toda bracket} \rm]
Let
$$ f: \Sigma^{-1} K_1 \rightarrow V(0) $$
be the attaching map of the first $V(0)$-coskeleton of $K$ to the zeroth
$V(0)$-cell, 
so that the cofiber of $f$ is $K$.
Let $\nu : K_1 \rightarrow \Sigma^nV(0)$ be the projection onto the top
$V(0)$-cell.
Suppose $\alpha$ is an element of $\pi_t(X \wedge V(0))$.
We have
$$ \pi_t(X \wedge V(0)) \xleftarrow{\nu_*} \pi_{t+n}(X \wedge K_1)
\xrightarrow{f_*} \pi_{t+n-1}(X \wedge V(0)). $$
We say the modified $K$-Toda bracket 
$$ \bra{K}'(\alpha) \subseteq \pi_{t+n-1}(X \wedge V(0)) $$
is \emph{defined} if $\alpha$ is in the image of $\nu_*$.  Then the
modified $K$-Toda bracket is the 
collection of all $f_*(\gamma) \in \pi_{t+n-1}(X \wedge V(0))$ 
where $\gamma \in 
\pi_{t+n}(X \wedge K_1)$ is any element satisfying $\nu_*(\gamma) =
\alpha$.
\end{defn}

Similarly we may provide a modified version of the dually defined and
$E$-ASS
$K$-Toda brackets given in Section~\ref{sec:toda}. 
Furthermore, for an $A_*$-comodule $M$ which is
cofree over the coalgebra $E[\tau_0]$, it is straightforward to define
the modified algebraic $M$-Toda bracket by varying the definition of the
$M$-Toda bracket given in Section~\ref{sec:mss}.

Modified versions of all of the results in Section~\ref{sec:results} hold,
with proofs that go through with only superficial changes.
Specifically, one must make the following adjustments.
\begin{itemize}
\item Replace $P$ with $P'$
\item Replace root invariants with modified root invariants.
\item Replace $E$-root invariants with modified $E$-root invariants.
\item Replace filtered root invariants with modified filtered root
invariants.
\item Replace $K$-Toda brackets with modified $K$-Toda brackets.
\item Replace the $E$-ASS for $X$ with the $E$-ASS for $X \wedge V(0)$.
\end{itemize}
Furthermore, 
there are
algebraic modified versions of Theorems~\ref{thmA}, 
\ref{thmB} and \ref{thmC}, Propositions~\ref{prop:Eroot} and 
\ref{prop:E^Eroot}, and Corollaries~\ref{corA} and \ref{cor:E^Eroot}, 
where one
makes the following replacements:
\begin{itemize}
\item Replace $P$ with $P'$.
\item Replace $\pi_*$ with $\ext$.
\item Replace root invariants with modified algebraic root invariants.
\item Replace $E$-root invariants with modified algebraic $E$-root invariants.
\item Replace filtered root invariants with modified algebraic filtered root
invariants.
\item Replace $K$-Toda brackets with modified algebraic $H_*(K)$-Toda brackets.
\item Specialize from $E$ to $BP$.
\item Replace the $E$-ASS for $X$ with the MSS for $H_*X \wedge V(0)$
\end{itemize}

We now modify our definitions even further, to give variations of the
second modified filtered root invariant $R''(\alpha)$ given in
\cite{MahowaldRavenel}.  
In what follows we shall always be working at a prime $p \ge 3$.

It is observed in
\cite{MahowaldRavenel} that there is a splitting
$$ P_{kq-1}^{(l+1)q} \wedge V(0) \simeq {P''}_{kq-1}^{lq-1} \vee
\Sigma^{kq}V(0) \vee \Sigma^{(l+1)q-1}V(0)
$$
In the notation of \cite{MahowaldRavenel}, we have
$$ {P''}_{kq-1}^{lq-1} = \br{P}_{kq-1}^{lq} $$
The spectra ${P''}_{kq-1}^{lq-1}$ are built out of the Smith-Toda complex 
$V(1)$. There is a $V(1)$-cell every in dimension from $kq-1$ to $lq-1$ 
congruent to $-1 \mod q$.
The decomposition of $P''_{-\infty}$ into
$V(1)$-cells is described on the level of cohomology in the following lemma,
which is a straightforward computation.

\begin{lem}\label{lem:QH^*P''}
Let $e_n^*$ be the generator of $H^*(P)_{-\infty}$ in dimension $n$, where $n
\equiv 0,-1 \pmod q$.  The cohomology of the Moore spectrum is given by
$$ H^*(V(0)) = E[Q_0] $$
as a module over the Steenrod algebra.
The cohomology $H^*(P'')_{-\infty} = H^*(P \wedge V(0))_{-\infty}$ is free
over the subalgebra $E[Q_0, Q_1]$ of the Steenrod algebra on the generators 
$e_{kq-1}^* \otimes 1$.
The action $E[Q_0,Q_1]$ on the free $E[Q_0,Q_1]$-submodule generated by
$e_{kq-1}^* \otimes 1$
is then given by the following formulas.
\begin{align*}
Q_0(e_{kq-1}^* \otimes 1) & = e_{kq}^* \otimes 1 - e_{kq-1}^* \otimes Q_0
\\
Q_1(e_{kq-1}^* \otimes 1) & = e_{(k+1)q}^* \otimes 1
\\
Q_1(e_{kq}^* \otimes 1 - e_{kq-1}^* \otimes Q_0) & = -e_{(k+1)q}^* \otimes
Q_0 
\\
Q_0(e_{(k+1)q}^* \otimes 1) & = e_{(k+1)q}^* \otimes Q_0
\end{align*}
\end{lem}

It is convenient to allow arbitrary subscripts and superscripts, so we
define for integers $M \le N$
$$ {P''}^N_M = {P''}_{kq-1}^{lq-1} $$
where $k$ (respectively $l$) is minimal (maximal) such that $M \le kq-1 \le
N$ ($M \le lq-1 \le N$).  If there is no such $k$ and $l$, then we have
${P''}^N_M = \ast$.

The following lemma is Lemma~3.7(e) of \cite{MahowaldRavenel}.
\begin{lem}\label{lem:H^*P''}
If $k$ and $l$ are congruent to $0 \pmod p$, then 
$H^*({P''}_{kq-1}^{lq-1})$ is free over the subalgebra $A(1)$ of the
Steenrod algebra.  
\end{lem}
Since there is no difference between the spectra $P''_{-\infty}$ and 
$P_{-\infty} \wedge V(0)$, and since $V(0)$ is
$p$-complete,
there is an analog of Lin's theorem \cite[Lemma 3.7(b)]{MahowaldRavenel}.
\begin{lem}\label{lem:lin''}
The map
$$ \Sigma^{-1}V(0) \rightarrow P''_{-\infty} $$
is an equivalence.
\end{lem}

In light of Lemma~\ref{lem:lin''}, Mahowald and Ravenel
\cite{MahowaldRavenel} define a second
modified root invariant.

\begin{defn}[{\it Second modified root invariant} \rm]
Let $X$ be a finite complex, and suppose we are given $\alpha \in
\pi_t(X \wedge V(0))$.  The \emph{second modified root invariant} 
of $\alpha$ is the coset 
of all dotted
arrows making the following diagram commute.
$$ \xymatrix{
S^{t} \ar@{.>}[r] \ar[d]^\alpha & \Sigma^{-N+1} V(1) \wedge X \ar[dd]  
\\
X \wedge V(0) \ar[d]
\\
tX \wedge V(0) \ar[r] & \Sigma P''_{-N} \wedge X 
}$$
This coset is denoted $R''(\alpha)$.  Here $N$ is chosen to be minimal
such that the composite $S^{t} \rightarrow \Sigma P''_{-N} \wedge X$ 
is non-trivial.
\end{defn}

The definition of the second modified root invariant may be obtained from
the definition of the (first) modified root invariant by smashing $X$ with
$V(0)$, replacing $P'$ with $P''$, and replacing $V(0)$ with $V(1)$.
Similarly, one may obtain second modified versions of all of the other
definitions we have been working with, as summarized below.  The 
second modified lambda complex
$W''_k(H_*X)$ is the appropriate subcomplex of $\br{\Lambda}_{(1)} \otimes
H_*(X)$, where $\br{\Lambda}_{(1)} = \br{\Lambda}/(v_0,v_1)$.
\begin{itemize}
\item Second modified $E$-root invariants
$$ R''_E: \pi_*(X\wedge V(0)) \rightsquigarrow \pi_*(E \wedge V(1) \wedge X) $$
\item Second modified algebraic root invariants
$$ R''_{alg}: \ext(H_*X \wedge V(0)) \rightsquigarrow \ext(H_*V(1) \wedge X) $$
\item Second modified filtered Tate spectra 
$$W_I({P''}^J)$$
\item Second modified filtered root invariants 
$$ {R_E^{[k]}}'': \pi_*(X \wedge V(0)) \rightsquigarrow 
\pi_*(W_k^k(V(1) \wedge X)) $$
\item Second modified algebraic $E$-root invariants
$$ R''_{E,alg} : \ext(H_*X \wedge V(0)) \rightsquigarrow 
\ext(H_*E \wedge V(1) \wedge X)
$$
\item Second modified algebraic filtered root invariants
$$ {R^{[k]}_{E,}}''_{alg} : \ext(H_*X \wedge V(0)) \rightsquigarrow 
H^*({W''}_k^k(H_*X)) $$
\end{itemize}
If $K$ is a finite complex built from $V(1)$ with single bottom and top
dimensional $V(1)$-cells, then we may define the second modified $K$-Toda
bracket 
$$ \bra{K}'': \pi_*(V(1) \wedge X) \rightsquigarrow \pi_*(V(1) \wedge X).
$$
Likewise, if $M$ is an $A_*$-comodule which is cofree over
$E[\tau_0,\tau_1]$, one can define a second modified algebraic $M$-Toda
bracket on $\ext(H_*V(1) \wedge X)$. 
Just as outlined for the case of modified root invariants in the first half
of this section, second modified
and second modified algebraic versions of the results of
Section~\ref{sec:results} hold.


\section{Computation of some infinite families of root invariants at $p =
3$}\label{sec:modcomp}

In this section we will extend the computation
$$ (-v_2)^i \in R''(v_1^i) $$
of \cite{MahowaldRavenel} for $p \ge 5$ to the prime
$3$ for $i \equiv 0,1,5 \pmod 9$.  
We will use our modified root invariant computations to deduce that, at the
prime $3$, the root invariant of the element $\alpha_i \in \pi_*^S$ 
is given by the element $\beta_i$ for $i \equiv 0,1,5 \pmod 9$.

Throughout this section we work at the prime $3$.  Low dimensional
computations indicate that there is a map
$$ v_2: S^{16} \rightarrow V(1).$$
Oka \cite{Oka} demonstrates that there is a map
$$ v_2^5: S^{80} \rightarrow V(1).$$
Composing these maps with iterates of the map
$$ v_2^9: \Sigma^{144}V(1) \rightarrow V(1) $$
given in \cite{BehrensPemmaraju}, we have maps
$$ v_2^i : S^{16i} \rightarrow V(1) $$
for $i \cong 0,1,5 \pmod 9$.  The computation of $\pi_*(L_2V(1))$ given in
\cite{GoerssHennMahowald} indicates that these are the only $i$ for
which $v_2^i$ can exist.

\begin{thm}\label{thm:R(v_1^i)}
At the prime $3$, for $i \equiv 0,1,5 \pmod 9$,
the second modified root invariant of $v_1^{i} 
\in \pi_*(V(0))$ is given by
$$ (-v_2)^i \in R''(v_1^{i}). $$
\end{thm}

We shall prove Theorem~\ref{thm:R(v_1^i)} with a sequence of lemmas.

\begin{lem}\label{lem:R''_BP}
For all $i$, we have the second modified $BP$-root invariant
$$ (-v_2)^i \in R''_{BP}(v_1^i). $$
\end{lem}

\begin{proof}
Modulo $I_1$, the $3$-series of the formal group associated to
$BP$ is given by
$$[3](x) \equiv v_1 x^3 + v_2 x^{9} + \cdots \pmod {I_1}. $$
Thus in $tBP \wedge V(0)_*$ (see Lemma~\ref{lem:tE_*}), we have
$$ v_1 = -v_2 x^6 + \cdots $$
and upon taking the $i^\mathrm{th}$ power, we have
$$ (v_1)^i = (-v_2)^i x^{6i} + \cdots $$
Using the second modified version of 
Corollary~\ref{cor:Eroot}, we may deduce the result.
\end{proof}

\begin{lem}\label{lem:R''_BP,alg}
For all $i$, we have the second modified algebraic $BP$-root invariant
$$ (-v_2)^i \in R''_{BP,alg}(v_1^i). $$
Here, $v_1$ and $v_2$ are viewed as elements in the cohomology of the 
periodic lambda
algebra.
\end{lem}

\begin{proof}
Apply $\ext(H_*-)$ to the diagram which realizes the second modified
$BP$-root invariant given by Lemma~\ref{lem:R''_BP}.
\end{proof}

In order to eliminate obstructions in higher Adams filtration, we prove the
following lemma.  This lemma is essentially contained in the proof of
Lemma~3.10 of \cite{MahowaldRavenel}.

\begin{lem}\label{lem:extP''}
We have the following $\ext$ calculation.
$$ \ext^{s,iq-1+s}(H_*P''_{-3iq-1}) = \begin{cases}
0 & s > i \\
\FF_p\{v_2^i\} & s = i
\end{cases}
$$ 
\end{lem}

\begin{proof}
In Lemma~\ref{lem:H^*P''}, we saw that $H^*(P''_{-3iq-1})$ is free
over the subalgebra $A(1)$ of the Steenrod algebra on generators in
dimensions congruent to $-1 \pmod {3q}$.  Therefore
$\ext_{A_*}(H_*P''_{-3iq-1})$ is built out of 
$\ext_{A_*}(\Sigma^{3kq-1} A(1)_*)$ for $k \ge -i$.  Now the May spectral
sequence $E_2$-term for $\ext_{A_*}(A(1)_*)$ is given by
$$ E_2 = E[h_{i,j} \: : \:  i \ge 1, j \ge 0, i+j \ge 2] \otimes P[b_{i,j}
\: : \: i \ge 1, j \ge 0, i+j \ge 2] \otimes P[v_i \: : \: i \ge 2]. $$
Analyzing the degrees that these elements live in establishes that the only
elements in $\ext_{A_*}(A(1)_*)$ which lie on or above the line of slope
$1/\abs{v_2}$ (in the $(t-s,s)$-plane) 
are those of the form $h_{1,1}^{\epsilon_1}
h_{2,0}^{\epsilon_2} v_2^k$ for $\epsilon_j = 0,1$.  Therefore, one easily
deduces that there are no
elements which lie above the line of slope $1/\abs{v_2}$ passing
through the point $(0,-3iq-1)$ which lie in the dimensions we are
considering, and the only element which lies in the $iq-1$ stem with Adams
filtration equal to $i$ is $v_2^i$.
\end{proof}

\begin{lem}\label{lem:R''_alg}
For all $i$, we have the second modified algebraic root invariant
$$ (-v_2)^i \in R''_{alg}(v_1^i) $$
\end{lem}

\begin{proof}
Applying the second modified algebraic version of
Proposition~\ref{prop:Eroot}, 
we may deduce
that the zeroth second modified algebraic filtered root invariant is given
by
$$ (-v_2^i) \in {R^{[0]}_{BP,}}''_{alg}(v_1^{i}). $$
The element $v_2^i$ is a permanent cycle in the MSS for $\ext(H_*V(1))$, so
we may use the second modified algebraic version of Theorem~\ref{thmA} to
deduce that the difference of the images of $(-v_2)^i$ and $v_1^i$ in 
$$ \ext^{i,iq-1+i}(H_*P''_{-3iq-1}) $$
is of $BP$-filtration greater than $0$.  Lemma~\ref{lem:extP''} says that there
are no non-zero elements of this $\ext$ group which could represent the
difference of these images, so we actually have
$$ (-v_2)^i \in R''_{alg}(v_1^i). $$
\end{proof}

We shall need a slightly different result for $i \equiv 5 \pmod 9$.

\begin{lem}\label{lem:v_2^5}
Modulo indeterminacy, the second modified algebraic root invariant of 
$v_1^{9t+5}$ is also 
given by
$$ (-v_2)^{9t+5} \pm v_2^{9t} v_3 g_0 b_0 \in R''_{alg}(v_1^{9t+5}) $$
\end{lem}

\begin{proof}
By Lemma~\ref{lem:extP''}, the element 
$ v_2^{9t} v_3 g_0 b_0 $ is in the kernel of the map
$$
\ext^{i, iq-1+i}(H_*\Sigma^{-3iq-1}V(1)) \rightarrow 
\ext^{i, iq-1+i}(H_*P''_{-3iq-1}) $$
where $i = 9t+5$.
Therefore it is in the indeterminacy of the second modified algebraic root
invariant.
\end{proof}

\begin{proof}[Proof of Theorem~\ref{thm:R(v_1^i)}]
The element $v_2^i$ is a permanent cycle in the ASS for $V(1)$ for $i = 0,1$.  
The existence of the self map $v_2^9$ of the 
complex $V(1)$ implies that $v_2^i$ is a permanent cycle in the ASS for $i
\equiv 0,1 \pmod 9$.  However, in the ASS for $V(1)$, the element $v_2^5$ is
\emph{not} a
permanent cycle; the element
$$ v_2^5 \pm v_3 g_0 b_0 $$
is the permanent cycle that detects the Oka element $v_2^5 \in
\pi_{80}(V(1))$ \cite[Remark 8.2]{BehrensPemmaraju}.  Thus the element
$v_2^{9t+5} \in \pi_*(V(1))$ is detected in the ASS by the permanent cycle  
$v_2^{9t+5} \pm v_2^{9t} v_3 g_0 b_0$.  For uniformity of notation,
let $\td{v_2^i}$ denote the ASS element
that detects $v_2^i \in \pi_*(V(1))$ for $i \equiv 0,1,5 \pmod 9$.

Lemmas~\ref{lem:R''_alg} and \ref{lem:v_2^5} imply that we have second
modified algebraic
root invariants 
$$ (-1)^i\td{v_2^i} \in R''_{alg}(v_1^i) $$
where $i \equiv 0,1,5 \pmod 9$.
By the second modified version of Theorem~\ref{thmD}, we have second
modified filtered root
invariants given by
$$ (-1)^i\td{v_2^i} \in {R^{[i]}_{H}}''(v_1^i). $$
The elements $\td{v_2^{i}}$ are permanent cycles in the ASS, 
so by the second modified version of Theorem~\ref{thmA}, we may conclude
that the difference between the images of $(-v_2)^i$ and $v_1^i$ in
$ \pi_{iq-1}(P''_{-3iq-1}) $
has Adams filtration greater than $i$.  According to
Lemma~\ref{lem:extP''}, there are no such elements, so the
images of $(-v_2)^i$
and $v_1^i$ in $\pi_{iq-1}(P''_{-3iq-1})$ are actually equal, and we have
$$ (-v_2)^i \in R''(v_1^i). $$
\end{proof}

We will now use our second modified root invariant computations to deduce the
following.

\begin{thm}
We have root invariants
$$ (-1)^{i+1} \beta_i \in R(\alpha_i)$$ 
for $i \equiv 0,1,5 \pmod 9$ at $p = 3$.  
\end{thm}

\begin{proof}
Let $\nu$ denote the map which is given by projection onto the top cell of
$V(0)$.  Since $P''_{jq-1}$ is a summand of $P_{jq-1} \wedge V(0)$, 
we have an induced
map
$$ \nu': P''_{jq-1} \rightarrow \Sigma P_{(j+1)q}. $$
We have the following diagram for $i \equiv 0,1,5 \pmod 9$.
$$
\xymatrix{
S^{qi-1} \ar[d]_{\alpha_i} \ar[dr]^{v_1^i} \ar[rr]_{(-v_2)^i}
\ar@/^2pc/[rrr]^{(-1)^i \beta_i} &&
\Sigma^{-3iq-1} V(1) \ar[dd] \ar[r]_{\nu''} &
S^{(-3i+1)q+1} \ar[ddd]_{\iota} \ar@{}[ddl]|{(-1)}
\\
S^0 \ar[dd] &
\Sigma^{-1}V(0) \ar[l]^{\nu} \ar[d] 
\\
& P''_{-\infty} \ar[dl]_\nu \ar[r] \ar@{}[uur]|{(\ast)} &
P''_{-3iq-1} \ar[dr]^{\nu'} 
\\
\Sigma P_{-\infty} \ar[rrr] &&&
\Sigma P_{(-3i+1)q}
}$$
The commutivity of the inner portion $(\ast)$ is our 
modified root
invariant computations of $v_1^i$.  
The portion of the diagram marked with a ``$(-1)$'' commutes with the
introduction of a minus sign.  This is due to the fact that the elements
$\beta_i$ are given by $\nu''_*(v_2^i)$ where
$$ \nu'' : V(1) \rightarrow S^6 $$
is the projection onto the cell corresponding to the element $Q_1 Q_0$ in 
the cohomology group
$$ H^*(V(1)) = E[Q_0,Q_1]. $$
According to the computations of Lemma~\ref{lem:QH^*P''}, this is the
negative of the cell in $P''_{-i3q-1}$ which corresponds to the cohomology
element
$$ e_{(-3i+1)q}^*\otimes Q_0 = Q_0 Q_1(e_{-3iq-1}^* \otimes 1). $$ 
Following the outer portion of the diagram, and introducing the minus sign,
we see that we have
$$ (-1)^{i+1} \beta_i \in R(\alpha_i) $$
unless
$\iota_*(\beta_i) = 0$.  That is not possible by
Proposition~\ref{prop:beta_survive}.
\end{proof}

\begin{rmk}
If one believes that the root invariant takes $v_1^9$ multiplication to
$v_2^9$ multiplication at $p=3$, analogous to the results of Sadofsky
\cite{Sadofsky} and
Johnson \cite{Johnson},
then it is seems likely that $\beta_i \dotin R(\alpha_i)$ for 
$i \equiv 0,1,3,5,6 
\pmod
9$. The present methods will not extend because $v_2^i$ is not a
permanent cycle in the ANSS for $V(1)$ for $i \not \equiv 0,1,5 \pmod 9$.  In
fact, at the time of writing we still do not know that the elements
$\beta_{9t+3}$ exist.  The elements $\beta_{9t+1} \beta_1 \alpha_1$ and 
$\beta_{9t+1} \beta_1^4$ do exist, and are
non-trivial by calculations of Shimomura and Wang \cite{ShimomuraWang}.  
We conjecture that we have root invariants
\begin{align*}
\beta_{9t+1} \beta_1 \alpha_1 & \dotin R(\alpha_{9t+2}) \\
\beta_{9t+1} \beta_1^4 & \dotin R(\alpha_{9t+4})
\end{align*}
The case $t=0$ was already
handled in Section~\ref{sec:anss3}.
\end{rmk}


\end{document}